\documentclass[oneside,a4paper,12pt]{amsart}
\usepackage{amsmath,amsthm,amsfonts,amssymb}
\usepackage[english]{babel}
\usepackage{geometry,graphicx,color}
\usepackage{varioref}
\input diagxy
\xyoption{all}
\xyoption{ps}
\newcommand{\xar}[1]{\ensuremath{\xrightarrow{#1}}}

\newcommand{\hra}[1]{\ensuremath{\stackrel{#1}{\hookrightarrow}}}

\newcommand{\sar}[1]{\ensuremath{\stackrel{#1}{\rightsquigarrow}}}
\newcommand{\bc}[1]{\boldsymbol{\mathcal{#1}}}
\newcommand{\wh}[1]{\widehat{#1}}
\newcommand{\wt}[1]{\widetilde{#1}}
\newcommand{\ol}[1]{\overline{#1}}
\newcommand{\mr}[1]{\mathrm{#1}}
\newcommand{\mb}[1]{\mathbf{#1}}
\newcommand{\bs}[1]{\boldsymbol{#1}}
\newcommand{\bg}[1]{\boldsymbol{\mathfrak{#1}}}

\newcommand{\mf}[1]{\mathfrak{#1}}
\newcommand{\ms}[1]{\mathsf{#1}}

\newcommand{\bem}[1]{\textbf{\emph{#1}}}
\newcommand{\us}[2]{\underset{#1}{#2}}

\newcommand{\br}[1]{[\![#1]\!]}
\newcommand{\bcc}[1]{\{\!\!\{#1\}\!\!\}}
\newcommand{\mc}[1]{\mathcal{#1}}

\def\PL{\mathrm{PL}}

\def\op{\mathrm{op}}
 
\def\R{\mathbb{R}}

\def\N{\mathbb{N}}

\def\Ob{\mathop{\mathrm{Ob}}}
\def\Ord{\mathop{\mathrm{Ord}}}
\def\Mor{\mathop{\mathrm{Mor}}}
\def\int{\mathop{\mathrm{int}}}

\def\relint{\mathop{\mathrm{relint}}}
\def\op{\mathrm{op}}
\def\id{\mathrm{id}}
\def\dom{\mathop{\mathrm{dom}}}
\def\codom{\mathop{\mathrm{codom}}}
\def\supp{\mathop{\mathrm{supp}}}
\def\ra{\rangle}
\def\la{\langle}
\def\sd{\mathop{\mathrm{sd}}}

\def\st{\mathop{\mathrm{star}}}

\def\St{\mathop{\mc St}}
\def\lk{\mathop{\mathrm{link}}}
\def\P{\mathbf{P}}
\def\OPL{\mathrm{OPL}}
\def\CPL{\mathrm{CPL}}
\def\conv{\mathop{\mathrm{conv}}}
\newtheorem{thm}{Theorem}
\newtheorem*{thm*}{Theorem}
\newtheorem*{thm_A}{Theorem A}
\newtheorem*{thm_A'}{Theorem A'}
\newtheorem*{thm_B}{Theorem B}
\newtheorem*{thm_C}{Theorem C}
\newtheorem*{thm_D}{Theorem D}

\newtheorem{lemma}{Lemma}

\newtheorem{prop}{Proposition}
\theoremstyle{remark}

\newtheorem*{a_lemma}{$\pmb \frown$-Lemma}
\newtheorem*{s_lemma}{$\pmb \sqcap$-Lemma}
\title{Combinatorial fiber bundles and fragmentation of a fiberwise
PL-homeomorphism}
\author{}

\begin{document}

\maketitle
\begin{center}{\textbf{\large{N. Mn\"ev}}}\\
\vspace{.5cm} \small St.~Petersburg Department of the \\
Steklov Mathematical Institute,
\\St.~Petersburg, Russia. \\
\texttt{mnev@pdmi.ras.ru} \\
\end{center}

\begin{abstract}

With a compact $\PL$ manifold $X$ we associate a category $\bg T(X)$.
The objects of $\bg T(X)$ are all combinatorial manifolds of type $X$, and
morphisms are combinatorial assemblies. We prove that the
homotopy equivalence
$$B\bg T (X) \approx B\PL(X)$$ holds, where $\PL(X)$ is the simplicial
group of $\PL$-homeomorphisms.
Thus the space $ B\bg T (X)$ is a canonical countable (as a
CW-complex) model of $B\PL(X)$. As a result, we obtain functorial
pure combinatorial models for $\PL$ fiber bundles with fiber $X$
and a $\PL$ polyhedron $B$ as the base. Such a model looks like a $\bg
T (X)$-coloring of some triangulation $K$ of $B$. The vertices of
$K$ are colored by objects of $\bg T (X)$ and the arcs are colored
by morphisms in such a way that the diagram arising from the 2-skeleton
of $K$ is commutative.
Comparing
with the classical results of
geometric topology, we obtain combinatorial models of the real Grassmannian
in small dimensions: $B\bg T(S^{n-1}) \approx B\mr O(n)$ for
$n=1,2,3,4$. The result is proved in a sequence of results on similar
models of $B\PL(X)$. Special attention is  paid to the main
noncompact case $X=\R^n$ and to the tangent bundle and Gauss functor
of a combinatorial manifold. The trick that makes the proof possible
is a collection of lemmas on ``fragmentation of a fiberwise
homeomorphism,'' a generalization of the folklore lemma on
fragmentation of an isotopy.
\end{abstract}

\newpage \mbox{}
\newpage
\pagestyle{myheadings}
\markboth{}{N. MN\"EV, COMBINATORIAL FIBER BUNDLES}

\tableofcontents

\section{Introduction}
\subsection{}
\begin{figure}[b]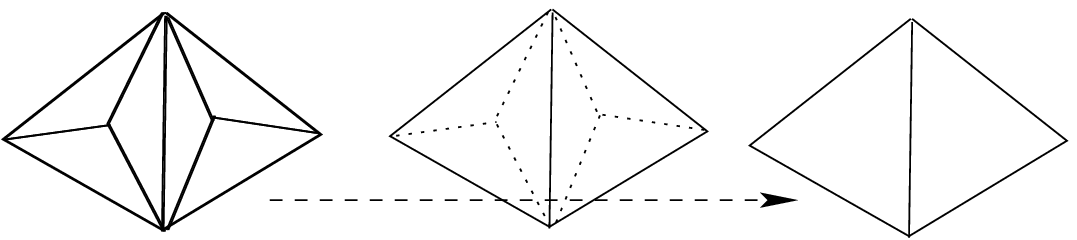 \caption{}\label{pic1}
\end{figure}

\begin{figure}[h] 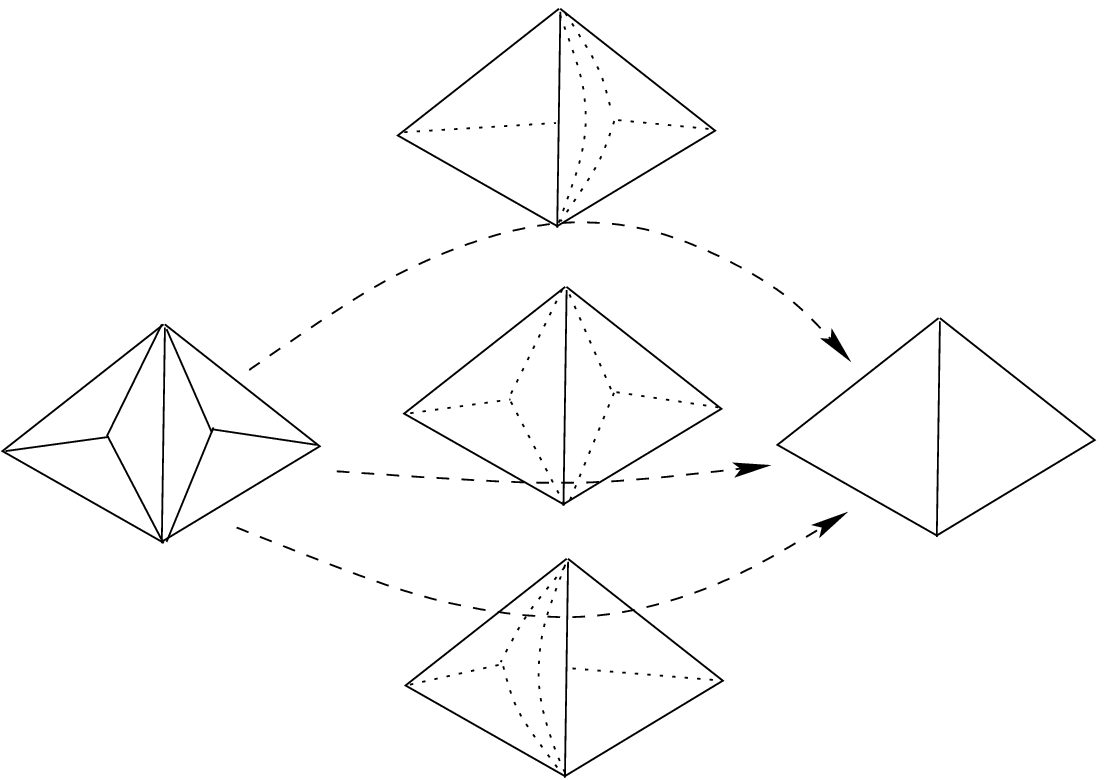 \caption{}\label{pic2}
\end{figure}

Let $X$ be a compact $\PL$ manifold. There is a natural
generalization of  piecewise linear triangulations of $X$, namely,
the structures of  piecewise linear regular cell (or ``ball'')
complexes\footnote{For exact definitions, see Sec.\vref{ball}. One
may imagine something like the boundary complex of a convex 3-polytope
as a ``ball complex'' and a planar 3-connected graph as an ``abstract
ball complex.'' Or one may simply think about geometric triangulations
instead of ``ball complexes'' and about combinatorial manifolds
instead of ``abstract ball complexes.''} on $X$.  The set of all
regular $\PL$ ball complexes on $X$ is partially ordered by
subdivision. We denote this poset by $\mb{R}(X)$. It is convenient
to consider a subdivision $\ms Q_0 \trianglelefteq \ms Q_1 $ of ball complexes
as a morphism of ``geometric assembly'' with source $\ms Q_0$ and
 target $\ms Q_1$. By forgetting the geometry,
to a geometric $\PL$ ball complex $\ms Q$ we can associate an abstract $\PL$ ball
complex $\mb P(\ms Q)$ (an ``abstract $\PL$ ball complex'' is a
natural generalization of the notion of an abstract simplicial
complex). The correspondence $\mb P$ sends the poset $\mb{R}(X)$ to some
new category $\bg{R}(X)$ whose objects are abstract
$\PL$ ball complexes and morphisms are ``abstract assemblies.'' One
may imagine an abstract assembly $\mf Q_0 \rightsquigarrow \mf Q_1$
of abstract  ball complexes as a way of gluing together the abstract
balls of $\mf Q_0$ into larger balls so as to obtain the
complex $\mf Q_1$. This way of gluing may be not unique.
Figure~\ref{pic1} should give an idea of a unique geometric assembly
of two particular geometric ball complexes, and
Fig.~\ref{pic2} should give an idea of three possible
combinatorial assemblies $$\mb P (\ms Q_0)=\mf Q_0\rightsquigarrow
\mf Q_1 = \mb P (\ms Q_1).$$
With the functor $\mb P$ we associate a cellular map of classifying
spaces $$B\mb R(X)\xar{B\mb P}B\bg R(X).$$ The map $B \mb P$ has a
description in terms of an action of the group of $\PL$ homeomorphisms
on $\PL$ ball complexes on $X$. Namely, the natural action of $\PL$
homeomorphisms on the set $\mb  R(X)$ can be extended to a cellular
action of a discrete group $\PL^\delta (X)$ on $B\mb R(X)$:
\begin{equation}\label{eq1} \PL(X)^\delta \times B\mb R(X)\xar{} B \mb R(X).\end{equation}
Then the cellular space of orbits $B \mb R(X)/\PL(X)^\delta$
coincides with $B \bg{R}(X)$, and $B\mb R(X)\xar{ B\mb
P}B \bg R(X)$ is a projection to the space of orbits. We should
mention that the action~(\ref{eq1}) is highly nonfree.

The category $\bg{R} (X)$ is an object of the classical combinatorial
topology of the manifold $X$. For example, Alexander's theorem on
combinatorial manifolds~\cite{Al} is the assertion that the space
$B\bg{R}(X)$ will remain connected if we restrict the class of all
morphisms to the more tame class of ``stellar assemblies.''

Denote by $\PL(X)$ the simplicial group of $\PL$ homeomorphisms of
$X$. Denote by $|\PL(X)|$ the cellular topological group that is the
geometric realization of $\PL(X)$. In  statistical models of
topological quantum field theory,
the simple fact is known that
\begin{equation}\label{eq2}\pi_1 B\bg{R}(X)\approx \pi_0 |\PL(X)|.\end{equation}
The group $\pi_0 |\PL(X)|$ is the mapping class group of the manifold
$X$. We prove the following generalization of
(\ref{eq2}).
\begin{thm_A}
The spaces $B\bg R (X)$ and  $B\PL(X)$ are homotopy equivalent.
\end{thm_A}
Thus the category $\bg R(X)$ is a discrete category that represents a
delooping of the simplicial group $\PL(X)$.  Generally, for any
topological (or simplicial) group there exists a delooping using
a discrete category or even a discrete
monoid~\cite{McDuff:1979}. The source of problems about discrete
categories representing deloopings of classical spaces is
the algebraic K-theory of topological spaces, starting with the famous
Hatcher's paper \cite{hatcher:1975}. The closest assertion to our
theorems is probably Steinberger's theorem \cite[Theorem~1,
p.~12]{steinberger:1986}, which is a refinement of Hatcher's conjecture
\cite[Proposition~2.5, p.~109]{hatcher:1975}. Steinberger's theorem says that
the discrete category of ordered simplicial complexes whose morphisms
are monotone maps with contractible preimages
of simplices
classifies Serre $\PL$ bundles.

Let $E\PL(X)$ be the contractible total space of the universal
principal bundle for the group $|\PL(X)|$, let
\begin{equation}\label{e3}|\PL(X)|\times E\PL(X)\xar{}E\PL(X)\end{equation} be
the canonical free action, and let $E\PL(X)\xar{}B\PL(X)$ be the
projection
to the space of orbits. Informally speaking, we prove that the
nonfree action~(\ref{eq1}) of the discrete group $\PL^\delta(X)$ on the
contractible space $B\mb R(X)$ can be deformed to the
canonical free action~(\ref{e3}) by a homotopy. In this form, our results are
relatives of Levitt's models for $B\PL$ (see \cite{Le}) presenting
$B\PL$ as orbit spaces. But in our case we are able to
eliminate geometry completely.

The first situations where Theorem A can be regarded as already known
appear when $\dim X = 1$ and $X$ is the
interval $I$ or the
circle
$S^1$. Here $\bg R(I)^\op$ is the category whose objects are all
finite ordinals and morphisms are generated by injective monotone maps
and the additional map of inverting the order;
$\bg R(S^1)^\op$
is the category whose objects are all cyclically ordered finite sets and
morphisms are generated by injective monotone maps and the additional map
of inverting the order.
The category $\bg R(S^1)^\op$ is closely related
to Connes' cyclic category. In these cases,
$$B\bg R(I)\approx B\PL(I)\approx BO(1) $$ and
$$B \bg  R(S^1) \approx B\PL(S^1) \approx BO(2). $$
The last assertion is a close relative of the theorem on the homotopy
type of the cyclic category (see \cite{Loday:1998}).

The assertion that $B\bg R (X) \approx B\PL(X)$ can be rephrased with the
help of the theory of representable homotopy functors. An
equivalent theorem states that there is a one-to-one correspondence
between isomorphism classes of $\PL$ fiber bundles with fiber $X$
on a polyhedron $B$ and concordance classes of $\bg R (X)$-colored
triangulations of $B$ (the vertices of triangulations are colored by
objects and the 1-simplices are colored by morphisms in such a way that
the 2-skeleton turns into a commutative diagram).
Thus a functorial
combinatorics of $\PL$ fiber bundle shows up,
and it is natural to pose
a question about generalizations of the known
combinatorial formulas for characteristic classes to $\PL$ fiber bundles.

We prove  an analog of  Theorem A for the subcategory $\bg T(X)$ of
$\bg R (X)$. The subcategory $\bg T(X)$ is formed by abstract
simplicial complexes (i.e., combinatorial manifolds of type $X$) and
their abstract assemblies.
\begin{thm_B}
The spaces $B\bg T (X)$ and  $B\PL(X)$ are homotopy equivalent.
\end{thm_B}

We may hope that there exist interesting subcategories of $\bg
T(X)$ for which the theorem is still valid. Probably, it is possible to
refine
our constructions so that they will work for the category of Brouwer
manifolds and linear representable assemblies. In the original
project (see \cite{AM}) it was supposed that the minimal
subcategory of $\bg T(X)$  modeling
$\PL(X)$ is the category of locally stellar
manifolds and stellar assemblies. This would be true in the case of
the positive solution of the famous problem (see~\cite[p.~14]{Hu}) of the
existence of a common geometric stellar subdivision for any two linear
triangulations of the simplex. The positive solution of this problem
would follow from the ``strong Oda conjecture'' on decomposition of
a birational isomorphism of smooth toric varieties. However,
a serious flaw was discovered in the proof \cite{Mo} of the strong Oda
conjecture (see \cite{AKMW}), and now the situation with the problem
of the existence of a
common stellar subdivision looks pessimistic. But in fact we need the following
weaker assertion, for which we
still have a hope.

\medskip\noindent
\textbf{Conjecture.}
{\em The poset of all linear triangulations of the
simplex ordered by
stellar subdivision is homotopy trivial.}
\medskip

Probably, with some cosmetic changes in the proofs, analogs of
Theorems~A,~B hold for any compact stratified polyhedron. We
analyze only the most important case of this kind: that
of the sphere with a fixed
point: $X=(S^n, \{*\})$.

\begin{figure}[h]
\includegraphics{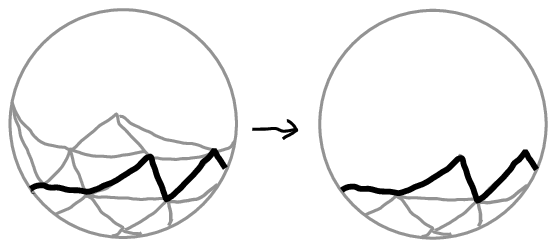}
\caption{}\label{aggr}
\end{figure}

According to Kuiper and Lashof \cite{KLI}, there is a one-to-one
correspondence between isomorphism classes of $\PL$ Milnor
$n$-microbundles, $\PL$ fiber bundles with fiber $\R^n$, and $\PL$
fiber bundles whose fiber is the pair $(S^n, \{*\})$. This
correspondence produces the homotopy equivalences $B\PL_n\approx
B\PL(\R^n) \approx B\PL(S^n, \{*\})$. We apply our scheme to
$B\PL(S^n, \{*\})$. As a result, we are able to build a combinatorial
model of $B\PL_n$ that has a remarkable property: the tangent bundle and
the Gauss map of a
combinatorial manifold obtain
a canonical
combinatorial form. Let $\bg R_n$ be the category whose objects are
abstract $n$-dimensional spherical $\PL$ ball complexes with a
marked  $n$-ball and morphisms are combinatorial
assemblies sending a marked ball to a marked ball (see Fig.~\ref{aggr}).
\begin{thm_C} The spaces
$B\bg R_n$ and $B\PL_n$ are homotopy equivalent.
\end{thm_C}

Let $\ms M^n$ be a combinatorial manifold. Let $(\mb P \ms M^n)^\op$
be the poset of simplices of $\ms M^n$ with the reversed order. There is a
functor $(\mb P \ms M^n)^\op \xar{\mb G} \bg R_n$. The value of $\mb
G$ at a simplex of $\ms M^n$ is the star of this simplex with an
extra $n$-ball attached by the link of the simplex so as to
obtain an $n$-sphere. This ball is the marked ball.

\begin{figure}\includegraphics[scale=0.5]{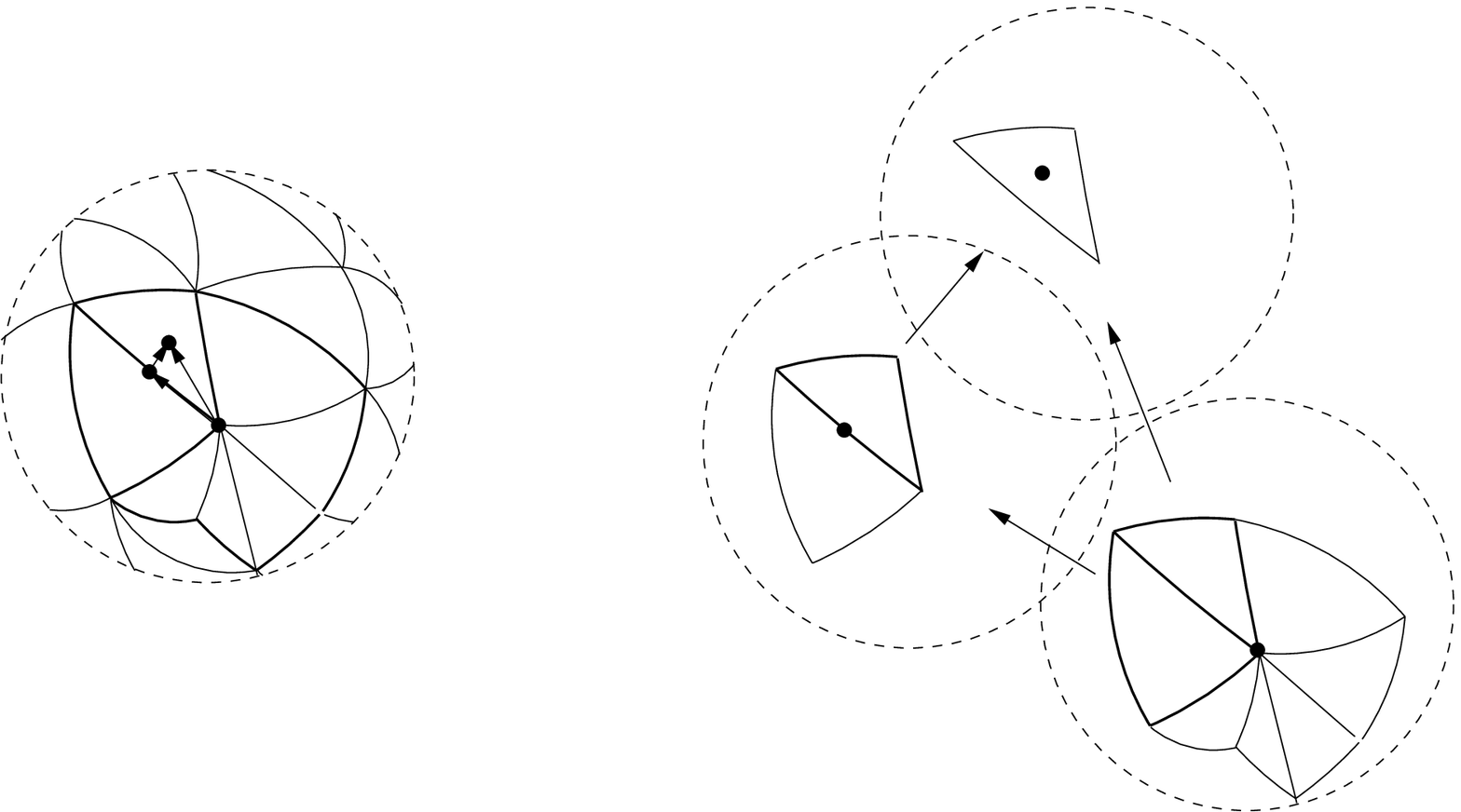}
\caption{} \label{tangent} \end{figure}

It is easy to imagine $\mb G$ for the combinatorial sphere $\ms S^n$.
For a simplex $s\in \ms S^n$, the sphere $\mb G(s)$ is the
$n$-dimensional abstract  ball complex that is the result of assemblying
all the simplices of $\ms S^n$ that
are not in the star of $s$ into one marked $n$-ball. If $s_0 \subset s_1 \in \ms S^n$, then
$\st s_0 \supset \st s_1$ and all the extra simplices of $\mb
G(s_0)$ are dissolved in the marked ball of $\mb G(s_1)$. This
defines the assembly morphism $\mb{G}(s_0 \subset s_1)$ (see Fig.~\ref{tangent}).

The space $B (\mb P \ms M^n)^\op$ is homeomorphic to $|\ms M^n|$.
\begin{thm_D} Two maps $|\ms M^n| \xar{} B\PL_n$ are homotopy equivalent:
the Gauss map for the tangent bundle of the manifold $|\ms M^n|$ and
the map $|\ms M^n|\xar{B \mb G} B\bg R_n\approx B\PL_n$.
\end{thm_D}

\begin{figure}
\includegraphics{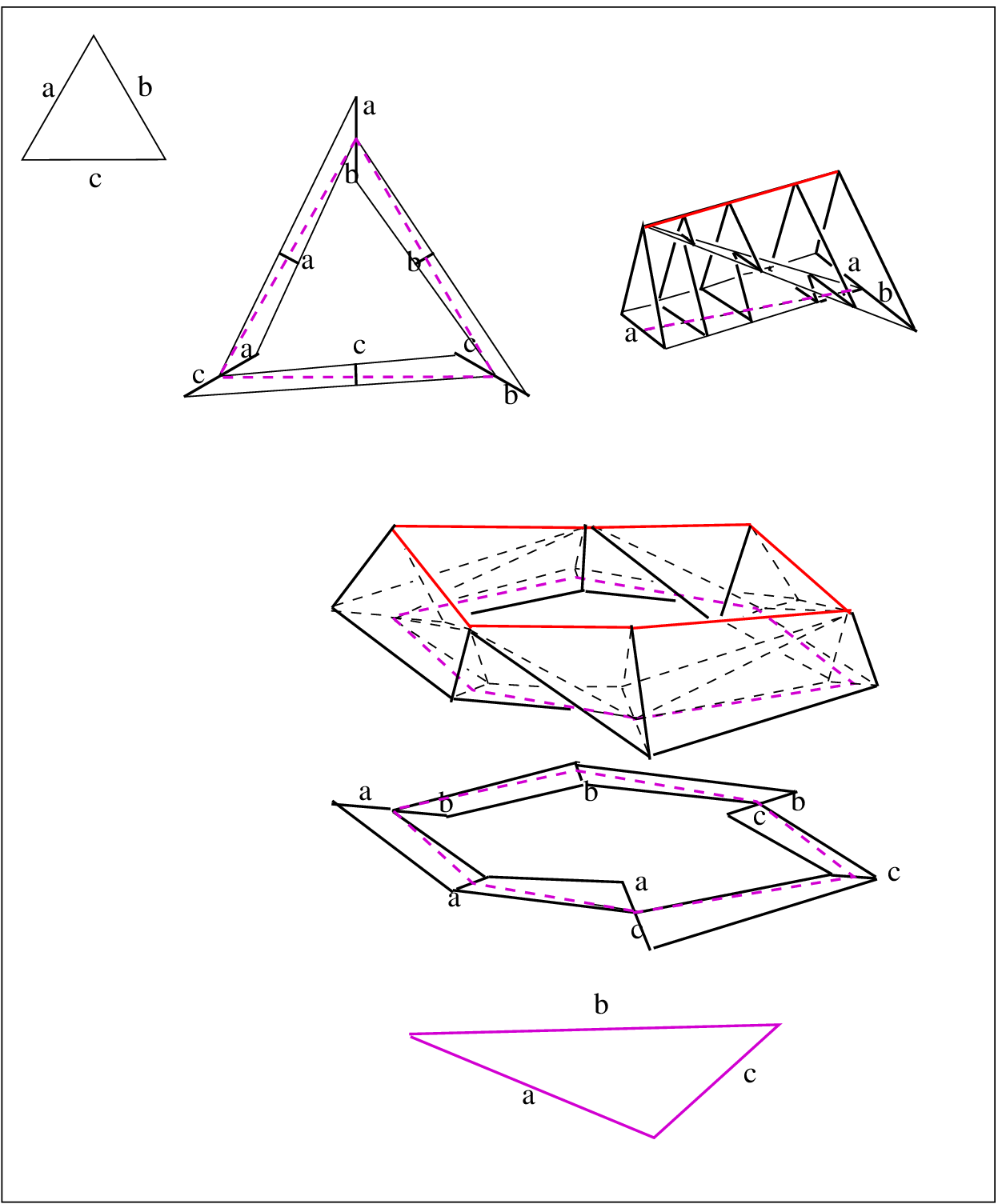} \caption{\label{e5}}
\end{figure}

The construction has the following remarkable property.  The
simplicial bundle $\mr{hocolim\,} \mb G \xar{} \mc N \mb P M^n$ is a
spherical bundle with zero and infinity sections. It is the
Kuiper--Lashof model of the tangent bundle $T|\ms M^n|$. The simplicial complex
$\mr{hocolim\,} \mb G$ is again a combinatorial manifold, and we
can iterate the construction. Already in the case of convex
polytopes, this situation looks interesting. Figure~\vref{e5} is
our attempt to draw  the canonical cellular structure on the tangent
bundle of a triangle.
\subsection{}
The main technical result of the paper is the proof of Theorem~A. Let us
describe how we would like to, but cannot,  prove such a theorem.
This speculation is borrowed  from \cite{steinberger:1986}.

There
are two simple functorial constructions, which are probably
originated from \cite{Whitehead:1939}. We will denote these
constructions by
$\mr{hocolim}$ and $\mr{hocolim}^{-1}$.
\begin{itemize}
\item The construction $\mr{hocolim}$ takes an $\bg R(X)$-coloring of a
triangulation of a polyhedron $B$ and produces a $\PL$ fiber bundle on $B$.
\item The construction $\mr{hocolim}^{-1}$ takes a $\PL$ fiber bundle on
$B$ with fiber $X$ and
produces an $\bg R(X)$-coloring of some triangulation of $B$. \\
\end{itemize}
The construction $\mr{hocolim}$ has many different names:
``iterated mapping cone'' \cite{Cohen:1967}, ``homotopy colimit,''
``Grothendieck construction'' \cite{GJ}, ``double bar-construction''
\cite[\S\,12]{May:1975}. The inverse construction,
$\mr{hocolim}^{-1}$, appears explicitly in \cite{hatcher:1975}
and uses triangulations of bundles.
One can triangulate a bundle.
Let us fix such a triangulation. We obtain a diagram
of combinatorial morphisms of the ball structures on $X$ in the fibers
over the vertices of the first barycentric subdivision of the base.
Then the dual morphisms form an $\bg R(X)$-coloring of the base.
The composite $ \mr{hocolim}\circ\mr{hocolim}^{-1}$, given a
bundle, produces an isomorphic bundle. We will obtain a short proof of
Theorem~A if we present a construction of a canonical concordance
between an $\bg R(X)$-coloring $\xi$ of a simplicial complex $K$ and
the coloring $\mr{hocolim}^{-1}\circ\mr{hocolim}\, \xi$ of the
complex $\sd K$. Unfortunately, there is no guarantee of the existence
of such a canonical construction. For the case of Serre bundles,
Steinberger \cite{steinberger:1986} constructed such a functorial
concordance, but his construction does not work for fiber bundles.
It follows from our results that $\xi$ and
$\mr{hocolim}^{-1}\circ\mr{hocolim}\, \xi$ are nevertheless
concordant, but our concordance is  absolutely transcendental.
\subsection{}
The core of our proof of  Theorem~A is Lemma~\vref{lem_frag_on_qube}
``on simultaneous fragmentation of fiberwise $\PL$ homeomorphisms of
the trivial fiber bundle over the cube.'' This lemma is a
straightforward generalization of  the lemma on fragmentation of an
isotopy. Let us roughly describe our scheme of reasoning.

Assume that we wish to proof the fact of the homotopy equivalence $B\bg R
(X)\approx B\PL(X)$ in its reformulation in terms of Brown's theory of
representable homotopy functors. Thus we wish to find a functorial
way to relate a fiber bundle with an $\bg R (X)$-coloring of a
polyhedron $B$, and {\it vice versa}. Let us describe
some process of constructing a bundle from a coloring. We are unable
to fight with the combinatorics of the $\mr{hocolim}$-construction, so we
will replace it by the traditional construction of
pasting trivializations
with the help of structure homeomorphisms.

\begin{figure}\caption{\label{polycone00}}
$$\includegraphics[scale=0.7]{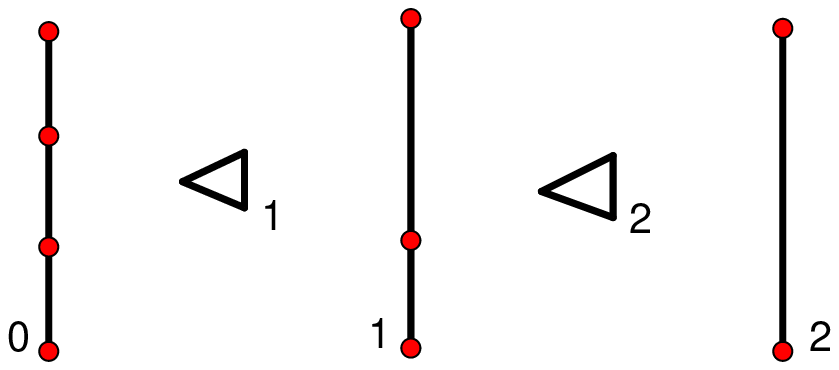}$$
\end{figure}

\begin{figure}\caption{\label{polycone11}}
$$\includegraphics[scale=0.7]{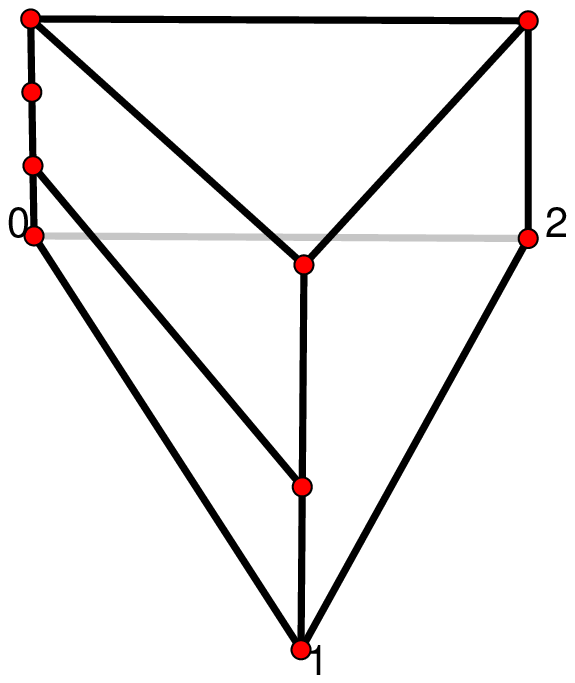}$$
\end{figure}

Let $K$ be an $\bg R (X)$-colored simplicial complex. The coloring of $K$
induces a coloring of each $k$-simplex of $K$ by a chain
$$\bg Q_0\rightsquigarrow \bg Q_1\rightsquigarrow{\ldots} \rightsquigarrow \bg Q_k$$
of abstract assemblies. This chain can be realized by a chain
\begin{equation}\label{qqqqq} Q=( \ms Q_0 \trianglelefteq \ms Q_1 \trianglelefteq {\ldots}  \trianglelefteq \ms Q_k)\end{equation}
of geometric assemblies of geometric $\PL$ ball complexes.  With the
chain $Q$ we can associate the ball decomposition
 of the trivial bundle
$X\times \Delta^k\xar{\pi}\Delta^k$ into the  horizontal
``prisms,'' i.e., the trivial subbundles
whose fiber is a ball.
Figures~\ref{polycone00} and~\ref{polycone11} illustrate
the construction of the
prismatic decomposition of
$\pi$ from a chain of
geometric assemblies.

\begin{figure}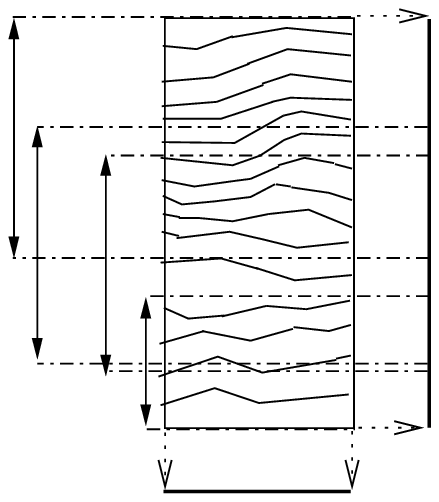\caption{\label{f4}}\end{figure}

\begin{figure}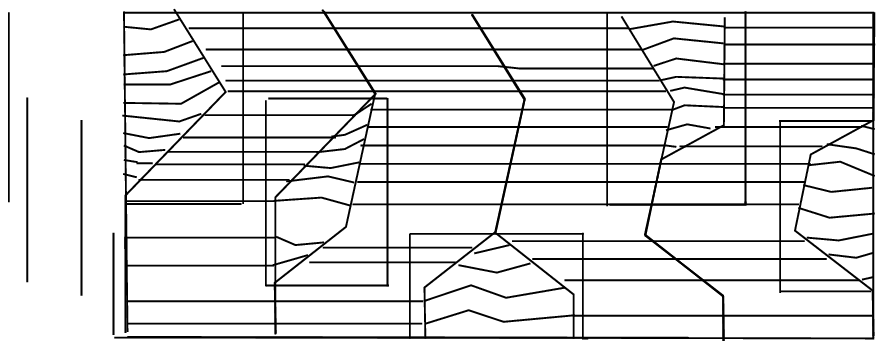\caption{\label{f5}}\end{figure}

The combinatorics of the coloring associates to any pair of
simplices $s_0 \subset s_1$ in $K$ a combinatorial isomorphism of
two prismatic structures on the trivial bundle over $s$. By  Alexander's
trick, all these combinatorial isomorphisms can be represented by
fiberwise structure $\PL$ homeomorphisms of the fiber bundle with
base $L$
and fiber $X$. All these structure homeomorphisms map
prisms to prisms. As a result, from the $\bg R(X)$-colorings we obtain
the class of fiber bundles with unusual, ``prismatic,'' structure homeomorphisms.
In this setup, the inverse problem is to learn
how to deform the structure homeomorphisms of an arbitrary $\PL$ fiber
bundle into the ``prismatic'' form and construct a consistent
coloring of the base
in a controllable way. At this point it is useful to recall the
proof of  the lemma on fragmentation of an isotopy.
In the $\PL$ case, this lemma was proved
by Hudson  \cite{Hu}. It states that for any cover
$U =\{ U_i \}_i$ of a manifold $X$ by open balls and for any
$\PL$ homeomorphism $X \xar{f}X$ isotopic to the identity there
exists a finite decomposition $f=f_1\circ{\ldots} \circ f_m$ such that
for every $i$ there exists $j$ with $\supp f_i \subset U_j$. The proof of the
fragmentation lemma contains more information than its statement. In
the proof we pick an arbitrary $\PL$ isotopy $F$ connecting $f$ and the
identity. Then we deform $F$ in the class of isotopies with fixed
ends to an isotopy $F'$ of a special form. The isotopy $F'$
corresponds to a chain of isotopies that are fixed on the complements
of the open balls from $U$. The isotopy  $F$ is a fiberwise
homeomorphism
$$\bfig \Vtriangle[X\times {[0,1]}`X\times {[0,1]}`{[0,1]};F`\pi_2`\pi_2]\efig$$
such that $F_0 = \id$ and $F_1=f$. The homeomorphism $F$
is in turn the same thing as a
one-dimensional foliation $\mc F$ on $X\times [0,1]$
transversal to the fibers of the projection $\pi_2$ (see Fig.~\ref{f4}).
The homeomorphism $F'$ corresponds to a foliation $\mc F'$ with the
following property: for any point $b\in [0,1]$, all points  $x
\in X$ such that the leaf of $\mc F$ ``is not horizontal'' at
$(x,b)$ are contained in an element of $ U$ (see \ref{f5}).
Inspecting  the figure of  $\mc F'$, it is easy to see that we can
subdivide the base $[0,1]$ into intervals $u_1,{\ldots} ,u_m$ and
introduce a prismatic structure on all subbundles $X\times u_i
\xar{\pi_2} u_i$ such that the induced homeomorphisms $F'\lfloor_{u_i}$
are prismatic. Thus the construction of the  fragmentation lemma
allows us to deform a fiberwise homeomorphism of the trivial bundle
over the interval into a system of prismatic homeomorphisms over a
subdivision of the interval. The deformation  $F\rightsquigarrow F'$ has
a canonical form; it has a coordinate generalization to
homeomorphisms of the trivial bundle over the cube. Our main task is to
formulate and analyze this generalization.
\subsection{The plan of the paper}
Section~2 contains a detailed definition of combinatorial assemblies of
 abstract ball complexes.

Section~3 contains a universal construction of the tautological
 $\PL(X)$ fiber bundle on $B\bg R(X)$. We emphasize the special role
 of ``prismatic''
 homeomorphisms, which are used for
 constructing the classifying map \begin{equation}\label{eqq}B\bg R(X)\xar{}B\PL(X).\end{equation}
 We introduce the simplicial groupoid of prismatic homeomorphisms $\mr{Prism}(X)$ and then
 construct a map of simplicial sets
 $\ol{\mc W}\mr{Prism}(X)\xar{}\mc W\mr{Prism}(X)$ that is a formal analog of the
$\mc W$-construction for the universal principal bundle. It is easy
to compare the  $\mc W$-construction for $\mr{Prism}(X)$ with the
standard $\mc W$-construction for $\PL(X)$ and with the map $B\mb
R(X)\xar{B\mb P}B\bg R(X)$. Thus  Theorem~A is reduced to Lemma~\vref{CombFrag}
``on prismatic trivialization.''

In Sec.~4 we reduce Lemma~\ref{CombFrag} ``on prismatic trivialization''
to the pure geometric Lemma~\vref{CombFrag1} ``on a common $\mb
R(X)$-triangulation of a family of fiberwise homeomorphisms.'' A
descriptive formulation of Lemma~\ref{CombFrag1} is as follows.
Consider the
simplicial set $\underline{\mr{Prism}}^m(X)$ whose
typical $k$-simplex is a collection $\la \ms Q, G_1,{\ldots} ,G_m\ra$,
where $\ms Q$ is a chain of geometric assemblies~(\ref{qqqqq}) and
$G_i\in \PL_k(X)$, $i=1,{\ldots} ,m$, are $\ms Q$-prismatic
homeomorphisms. That is, the homeomorphisms $G_i^{-1}$ send the prisms of
$\ms T(\ms Q)$ to prisms. The set $\underline{\mr{Prism}}^{m+1}(X)$
is naturally embedded into $\underline{\mr{Prism}}^m(X)\times \PL(X)$.
Lemma~\ref{CombFrag1} states that

\medskip\noindent
{\em the pair $(|\underline{\mr{Prism}}^m(X)\times \PL(X)|,
|\underline{\mr{Prism}}^{m+1}(X)|)$
is homotopy trivial.}

\medskip\noindent
The translation of this assertion into the common language reads as follows:

\medskip\noindent
{\em one can deform any new nonprismatic homeomorphism to the prismatic
form jointly with some family of homeomorphisms in such a way that
all prismatic homeomorphisms in the family will remain prismatic.}

\medskip\noindent
This fact is the central technical result of the paper. The plan of
the proof of Lemma~\ref{CombFrag1}  is contained in Sec.~5, and the proof
itself occupies Secs.~6--15.  We introduce the general notions of
``Alexandroff presheaf'' and  ``prismaticity
of a fiberwise
homeomorphism with respect to an Alexandroff presheaf.'' For these
generalized prismatic homeomorphisms we develop some surgery
centered around a generalized Hudson's construction of fragmentation
for a $\PL$ isotopy.

In Sec.~16 we describe how to tweak the general scheme in order
to obtain a proof of Theorem~B.

In Sec.~17  we describe how to tweak the general scheme in order to obtain a proof
of Theorem~C. Then we demonstrate that our combinatorial
construction of the tangent bundle represents the Milnor tangent
microbundle. This proves  Theorem~D.

\subsection{}
Our theorems appear as an answer to the natural question about
relations between geometric and abstract triangulations of a
manifold. This question arose at A.~M.~Vershik's seminar in the
author's student years. It was converted into conjectures
during the joint work with Peter Mani-Levitska and Laura Anderson.
The conjecture on a
combinatorial model of  $B\PL_n$ (which is now
Theorem~C) was a $\PL$ analog of the conjectures on
the MacPhersonian\footnote{The proof of these conjectures in \cite{B1}
contains
 a very serious flaw (see \cite{MnevBiss}).}
(a hypothetical combinatorial model for $B\mr
O(n)$ \cite{MnevZigler:1993}). We can mention  that our theorems
combined with the classical knowledge on the relations between $\PL(S^n)$,
$\mr{Diff}(S^n)$, and
 $\mr O(n+1)$ (see \cite{Hatcher:1983}) produce the following combinatorial
 models of the Grassmannians  $B\mr O(n)$, $n=1,2,3,4$:
$$B\bg R(S^{n-1})\approx B\bg T(S^{n-1})
\approx B \bg R_{n} \approx B\mr O(n).$$

As mentioned above, the first project \cite{AM} of proving our
theorems was based on the proof \cite{Mo} of the strong Oda's
conjecture. This proof appears to be wrong (see \cite{AKMW}). Our
current proofs are independent from Oda's conjecture.

The author is grateful to A.~M.~Vershik for permanent support and
wise advices, to the St.~Petersburg Department of Steklov
Institute of Mathematics for financial support and wonderful
atmosphere, to Peter Mani and Laura Anderson for countless
stimulative talks and fantastic hospitality
at Bern and College Station during the initial stage
of the project.

\section{Assemblies of ball complexes, the poset $\mb R(X)$, and the category $\bg R(X)$}
\label{ball} In this section we define geometric and combinatorial
assemblies of ball complexes.
\subsection{}
Our principal category is the category  $\PL$ of
piecewise linear Euclidean polyhedra and piecewise linear maps. The
foundations of $\PL$ topology can be found in the  books \cite{RS,
Hu} and the notes of Zeeman's seminar \cite{Ze1}.

\medskip
\noindent\textbf{Warning}: \emph{in this paper, all polyhedra, manifolds,
maps, etc.\ are assumed to be piecewise linear unless another category is
specified.}

\subsection{Ball complexes}
For general information on topological ball complexes, see
the book \cite{LW}. $\PL$ ball complexes appeared in $\PL$
geometric topology (e.g., in \cite{Casson:1967}). We recall
the standard definition of a ``topological ball complex''
or, equivalently, a ``finite regular $CW$-complex.''

A \bem{topological ball complex} is a finite cover $S$ of
a Hausdorff space $X$ by closed topological balls such that
\begin{itemize}
\item[(i)] the relative interiors of the balls from $S$ form a partition of $X$;
\item[(ii)] the boundary of every ball from $S$ is the union of balls of smaller
dimension.
\end{itemize}
A \bem{ $\PL$ ball complex} on a Euclidean polyhedron $X$ is a
finite cover $S$ of $X$ by closed $\PL$ balls such that the
conditions~(i), (ii) are satisfied.
The main example of a $\PL$ ball complex is a finite geometrical
simplicial complex. {\em In what follows, a ``ball complex''
means a ``$\PL$ ball complex.''}

\subsection{The category $\mb{PLball}$}
The category $\mb{PLball}$ is the category whose objects are ball
complexes and morphisms are maps that send the relative interiors of
balls into (not necessarily {\em onto}) the relative interiors of balls.
To be more  precise, a morphism  $( X_0,  S_0) \xar{} ( X_1,  S_1)$
is a pair $( h, \xi )$, where $X_0\xar{h}X_1$ is a $\PL$ map and $ S_0
\xar{\xi}  S_1$ is a map of ball sets such that for every $s\in
S_0$ the inclusion $h(\relint s)\subseteq \relint \xi(s)$ holds. The category
$\mb{PLball}$ is not very interesting, it merely contains our working
subcategories.

\subsection{Abstract ball complexes}
For every ball complex  $( X, S )$, the polyhedron $X$ is determined up to
homeomorphism by pure combinatorial data, namely, by the combinatorics of
the adjacency of balls of $S$. Let us formulate this assertion in detail.

Let $\ms D = (X, S)$ be a ball complex. Denote by $\P (\ms D)$ the
partial order by inclusion on $S$. Consider the abstract simplicial
complex $\Ord \P(\ms D)$, the order complex of the poset $\P (\ms D)$.
Let $|\Ord \P(\ms D)|$ be the
geometric realization of $\Ord
\P(\ms D)$. Let us introduce the standard notation. Given a poset $\mc P$
and an element $p\in \mc P$, denote by $p_\leq$ the subposet of $P$
formed by all elements that are less or equal to $p$ (the
``lower principal ideal generated by $p$''). Denote by $p_<$ the
ideal formed by all elements that are strictly less than $p$. The following
theorem holds (see \cite{LW}).

\begin{thm} \label{LW1} For every ball complex $\ms D=( X,{S})$ there is
a cellular homeomorphism $$( X,{S}) \approx (
\mathop{|\mathrm{Ord}}\P(\ms D)|,
 \{|\Ord p_\leq|\}_{p\in \P(\ms D)}).
$$
\end{thm}
 The polyhedra $\{|\Ord
p_<|\}_{p\in\P(\ms D)}$ from  Theorem~\ref{LW1} are automatically the boundary
spheres of cells. The last property allows us to define an
abstract ball complex. The following theorem holds (see \cite{Bj84}).

\begin{thm} \label{Bj84}
Let $\mc{P}$ be a finite poset satisfying the following property:
if $p \in \mc{P}$ and the rank of $p$
is equal to $k$, then $|\Ord p_<| {\approx} S^{k-1}$. Then
$$( | \Ord \mc{P}|,%
\{|\Ord p_\leq|\}_{p\in \mc{P}} ) $$ is a ball complex.
\end{thm}
A poset that satisfies the conditions of  Theorem~\ref{Bj84} is
called an \bem{abstract ball complex}. Theorems~\ref{LW1} and~\ref{Bj84}
were originally formulated for topological ball complexes,
but the proofs work in the $\PL$ category without any changes. Thus by an
\bem{abstract $\PL$ ball complex} (in this paper, it will be called
an ``abstract ball complex'') we mean a finite poset $\mc{P}$ such that
$|\Ord p_<|\us{\PL}{\approx} S^{k-1}$ for every $p \in \mc{P}$ of rank $k$.
The $\PL$ version of  Theorem~\ref{Bj84}
states that in this case $( | \Ord \mc{P}|, \{|\Ord p_\leq|\}_{p\in
\mc{P}} ) $ is a $\PL$ ball complex.

\subsection{The functor $\mb{PLball} \xar{\mb P }\mb{Posets}$}
Consider a $\mb{PLball}$-morphism $( X_0, S_0) \xar{( h,\xi)} ( X_1,
S_1)$. From the definition of a ball complex it follows that $\xi$ is
a morphism of the posets of balls. Therefore the correspondence $\ms
D \mapsto \mb P (\ms D)$ is a functor with values in the category
of posets.
\subsection{Geometric assemblies of ball complexes, the poset $\mb R(X)$}
\label{ss26} A \bem{geometric assembly} of ball complexes on $X$ is
a $\mb{PLball}$-morphism $( X, S_0) \xar{( h,\xi)} ( X,S_1)$ such
that $h$ is the identity map. Such a situation is possible only when
the relative interior of every ball $s\in S_0$ is contained in the
relative interior of the ball $\xi(s)\in S_1$. This means that the
partition of $X$ into the relative interiors of balls from $S_1$ is
subdivided by the partition into the relative interiors of balls from
$S_0$. Therefore a geometric assembly morphism $( X, S_0) \xar{} (
X,S_1)$ is unique if it exists. Thus the geometric assemblies form
a poset $\mb R (X) \hra{} \mb{PLball}$ on the set of all ball complexes
with the underlying polyhedron $X$. We denote a geometric assembly $\ms
Q_0 \xar{} \ms Q_1$ by $\ms Q_0 \trianglelefteq \ms Q_1$.
\subsection{The category $\bc{R}(X)$}\label{hom_sborki}
Consider the subcategory $\mb{\bc{R}}(X) \hra{} \mb{PLball}$ whose
objects are regular ball complexes on $X$ and morphisms are of the form
$( X, S_0) \xar{( h,\xi)} ( X,S_1)$, where $X\xar{h}X$ is a
homeomorphism. The poset $\mb{R}(X)$ sits in $\bc{R}(X)$ as a
subcategory.
The morphisms of $\bc{R}(X)$ are generated by two classes:
\begin{itemize}
\item[(i)] assemblies $\trianglelefteq$,
\item[(ii)] cellular homeomorphisms (i.e., homeomorphisms sending every
 ball onto a ball).
\end{itemize}
Obviously, the following proposition holds.
\begin{prop}\label{decomp}
Every morphism $\ms Q_0\xar{f}\ms Q_1$ of $\bc{R}(X)$ has two
canonical decompositions
\begin{equation} \label{h-r}
\bfig \Atriangle(0,0)/<-`>`>/<300,300>[{\ms Q_0'}`\ms Q_0`\ms Q_1;\trianglelefteq_1`h_1`f]
\Vtriangle(0,-300)/`>`{<-}/<300,300>[\ms Q_0`\ms Q_1`{\ms Q_1'};`h_2`\trianglelefteq_2] \efig
\end{equation}
into a homeomorphism and an assembly,
i.e., an assembly in the source and a homeomorphism in the target or an
assembly in the target and a homeomorphism in the source.
\end{prop}
\subsection{The category $\tilde{\bg R }(X)$
of combinatorial assemblies of geometric ball complexes and the
category $\bg{R}(X)$}~\label{2.7} Consider the category $\tilde{\bg R
}(X)$ whose objects are ball complexes on $X$ and an
$\tilde{\bg{R}}(X)$-morphism $\ms Q_0 \xar{\mf f} \ms Q_1$ is a
morphism of posets $\mb P \ms Q_0 \xar{\mf f} \mb P \ms Q_1 $ that
is representable by some $\bc R (X)$-morphism. That is, $\mf f$ is a
poset morphism such that there is a morphism $\ms Q_0 \xar{f} \ms Q_1$ for
which $\mb P (f) = \mf f$. Such poset morphisms will be called
\bem{abstract assemblies}. Using the decomposition from
Proposition~\ref{decomp}, we can easily establish the following fact.
\begin{prop}
The composition of abstract assemblies is an abstract assembly.
\end{prop}
Therefore $\tilde{\bg R}(X)$ is a well-defined category;
there is a forgetful functor $\bc R(X)\xar{\bc A} \tilde{\bg R}(X)$,
which is identical on objects and sends a morphism of $\bc
R(X)$ to the corresponding abstract assembly. The poset  $\mb R (X)$
is a subcategory of $\bc R (X)$. Denote by $\mb A $ the
composite
$$\mb A=(\mb R(X)\hra{} \bc R (X)\xar{\bc A} \tilde{\bg R}(X)).$$
Note that it is easy to give a transcendental but ``pure
combinatorial'' definition of an abstract assembly of abstract ball
complexes. Namely, the following proposition holds.
\begin{prop} \label{3}
A morphism
$\ms Q_0 \sar{\mf f} \ms Q_1 $ is an abstact assemly if and only if for every
$b \in \mb P \ms Q_1$ of rank $k$ the poset ${\mf
f}^{-1}(b_\leq) \subseteq \mb P \ms Q_0$ is an abstract ball complex
representing a
$k$-dimensional $\PL$ ball.
\end{prop}
From the $\PL$ variant of  Theorem~\vref{LW1} we know that any two
combinatorially isomorphic geometric ball complexes are isomorphic
in $\bc R (X)$. This allows us to identify, up to {\em equivalence of
categories}, the category $\tilde{\bg R}(X)$ with the category $\bg
R(X)$ whose objects are all abstract ball complexes representing $X$
by geometric realizations
and morphisms are abstract assemblies of
abstract ball complexes, i.e., poset morphisms $\mc P_0 \xar{\mu} \mc
P_1$, $\mc P_0, \mc P_1 \in \Ob \bg R(X)$, representable by $\bc R
(X)$-morphisms of geometric realizations. Thus we obtain the commutative
triangle of functors
\begin{equation}\label{tria}\bfig \dtriangle[\mb R(X)`\tilde{\bg R}(X)`\bg R (X),; \mb A`\mb P`\mb F] \efig\end{equation}
where $\mb F$ is the forgetful functor inducing an equivalence of the
categories $\tilde{\bg R}(X)$ and
 $\bg R(X)$.
 \subsection{$\mc N \bg R (X)$ as the orbit space of an action of $\PL^\delta (X)$
on a contractible space}
Now we switch to a more scientific point of view on the functor $\mb
R (X)\xar {\mb P}\bg R(X)$. Let us pass to the nerve $\mc N \mb R
(X)\xar{\mc N\mb P}\mc N\bg R(X)$ of this functor. There is an
action of a
discrete simplicial group $\PL^\delta (X)$ on $\mc N \mb
R (X)$. A typical $k$-simplex of $\PL^\delta (X)$ is an ordered
set $g=(g_0,{\ldots} ,g_k)$, $g_i \in \PL(X)$, $i=1,{\ldots} ,k$, $g_i=g_j$. If
$\ms Q \in \mc N_k \mb R (X)$, $\ms Q=(\ms Q_0 \trianglelefteq {\ldots} \trianglelefteq \ms
Q_k)$, then we set
$$g \ms Q = (g_0 \ms Q_0 \trianglelefteq g_1 \ms Q_1 \trianglelefteq {\ldots}  \trianglelefteq g_k \ms Q_k)\in \mc N \mb R (X).$$
This defines an action
$$ \PL^\delta(X)\times \mc N \mb R (X) \xar{} \mc N \mb R (X).$$
\begin{lemma}
$\mc N \bg R (X) = \mc N \mb R (X) /\PL^{\delta}(X)$, and
 $\mc N\mb R (X) \xar{\mc N \mb P} \mc N \bg R (X)$ is a
projection to the space of orbits.
\end{lemma}
\begin{proof}
By definition (see, for example, \cite[p.~303]{He}), the $k$-simplices of
$\mc N \mb R (X) /\PL^{\delta}(X)$ are the orbits of the action of
$\PL^\delta_k(X)$ on $\mc N_k \mb R(X)$ with the induced simplicial
structure: for $\ms Q \in \mc N_k \mb R(X)$, the orbits of $d_i \ms Q
$ and
 $s_i \ms Q$
are determined by the orbit of $\ms Q$. This exactly identifies $\mc N
\mb R (X) /\PL^{\delta}(X)$ and $\mc N \bg R (X)$.
\end{proof}
Let us also mention the following fundamental fact.
\begin{lemma} \label{lem2}
The space  $|\mc N \mb R(X)|$ is contractible.
\end{lemma}
\begin{proof}
This easily follows from the fact that $\mb R(X)$ is a poset and for
every pair $\ms Q_0, \ms Q_1 \in \mb R (X)$ there exists a common
subdivision, i.e., there is $\ms Q_3$  such that $\ms Q_3 \trianglelefteq \ms
Q_0$, $\ms Q_3 \trianglelefteq \ms Q_0$ (see \cite{Quillen:1978}).
\end{proof}
We wish to emphasize that only in the  category $\PL$ there exists a
common subdivision of ball complexes. Due to this fact, in other
categories the existence of analogs of  Lemma~\ref{lem2} is
problematic.

\section{Prismatic homeomorphisms, reduction of  Theorem~A
to the lemma on prismatic trivialization }\label{const1}

In this section, we present a construction of a bundle with base
$B$ and fiber $X$ from an $\mf{R}(X)$-coloring of a triangulation
of $B$. This construction is described universally, as a
simplicial morphism $\mc N\mf{R}(X)\xar{}\ol{\mc W} \PL(X)$ from the
nerve $\mf{R}(X)$ to the $\ol{\mc W}$-construction of $B\PL(X)$.

\subsection{Prismatic decompositions of trivial fiber bundles}\label{prism}
Let $\bs m$ be the totally ordered set $\{0<1<{\ldots} <m\}$, a ``finite
ordinal.'' Let $\mb N$ be the category of all finite ordinals and
monotone maps. Consider a chain of geometric assemblies of ball
complexes $\bs m \xar{\ms Q} \mb R(X)$. It this subsection,
with the chain $\ms Q=(\ms Q_0\trianglelefteq_1{\ldots} \trianglelefteq_m \ms Q_m)$
we will associate the
structure of a ball complex $\ms{T}(\ms Q)$ on the polyhedron $X
\times \Delta^m$ and a $\mb{PLball}$-morphism $\ms{T}(\ms Q)\xar{\ms
e ( \ms Q)}[\Delta^m] $. We denote by $[\Delta^m]$ the \bem{standard
``ball simplex,''} the ball complex formed by all faces of the
standard simplex $\Delta^m$. Let $\bs n \xar{\theta} \bs m $ be an
$\mb N$-morphism, and let $$\Delta^n\xar{\Delta(\theta)}\Delta^m$$ be
the corresponding cosimplicial morphism. By
$[\Delta^n]\xar{[\Delta(\theta)]}[\Delta^m]$ we denote the induced
$\mb{PLBall}$-morphism of the ball complexes of the standard simplices.
By $[d^i],[s^i]$ we denote the standard {\em cellular} cofaces and
codegenerations.

Note that in the category $\mb{PLball} $
there exist \bem{induced assemblies}, i.e., if a ball complex $\ms
B_0$ is a subcomplex of  $\ms B_1$ and $\ol{\ms{B}}_0 \trianglelefteq \ms
{B}_0$, then there exists a universal dashed arrow in the diagram
\begin{equation}
\bfig \label{ind_subd}
\square/^{(}..>`>`..>`^{(}->/[\ol{\ms B}_0`\ol{\ms B}_1`\ms B_0` \ms B_1.;`\trianglelefteq`\trianglelefteq`]
\efig
\end{equation}
The set of balls of the complex $\ol{\ms B}_1$ is the set of all
balls of the complex $\ms{B}_1$ with all the balls
of $\ms{B}_0$ deleted and all the balls of $\ol{\ms B}_0$ added.

Consider the flag of  faces of the simplex
\begin{equation*}\label{flag0} \Delta^0 \hra{d^1}
\Delta^1 \hra{d^2}\cdots\hra{d^{m}}\Delta^m.
\end{equation*}
Consider the induced flag of trivial fiber bundles
\begin{equation}\label{flag}
\xymatrix{{X\times \Delta^0} \ar[d] \ar@{^{(}->}^{\id \times d^1}[r] &
 {X\times \Delta^1} \ar[d] \ar@{^{(}->}^{\id \times d^2}[r] & \cdots \ar@{^{(}->}^{\id
\times d^m}[r] & {X\times \Delta^m} \ar[d] \\
 {\Delta^0}\ar@{^{(}->}[r]^{d^1}&
  {\Delta^1}\ar@{^{(}->}[r]^{d^2} & \cdots \ar@{^{(}->}^{d^m}[r]
  &{\Delta^m}.
}
\end{equation}
Consider the
following construction
on the total space of the flag~(\ref{flag}). Consider the  staircase of ball assemblies
and embeddings
\begin{equation}  \label{saw}
\xymatrix{ \ms Q_0\times [\Delta_0] \ar@{^{(}..>}[r] & \bullet  \ar@{^{(}..>}[r] & \bullet  \ar@{^{(}..>}[r] & \cdots  \ar@{^{(}..>}[r] &  \ms{T (Q)} \\
           \ms Q_1\times [\Delta^0] \ar@{^{(}->}[r]^{\id
\times [d^1]} \ar@{<-}[u]^{\trianglelefteq\times\id}& \ms Q_1
 \times [\Delta^1] \ar@{^{(}..>}[r] \ar@{<..}[u]^\trianglelefteq &
 \bullet \ar@{<..}[u]^\trianglelefteq \ar@{^{(}..>}[r] &
 \cdots \ar@{^{(}..>}[r] & \bullet  \ar@{<..}[u]^\trianglelefteq \\
            & \ms Q_2 \times [\Delta^1] \ar@{^{(}->}[r]^{\id \times
            [d^2]} \ar@{<-}[u]^{\trianglelefteq\times\id} & \ms Q_2\times [\Delta^2] \ar@{^{(}..>}[r]
             \ar@{<..}[u]^\trianglelefteq & \cdots \ar@{^{(}..>}[r] & \bullet \ar@{<..}[u]^\trianglelefteq \\
            &  & \vdots \ar@{<..}[u]^{\trianglelefteq \times \id} & \vdots &  \vdots \ar@{<..}[u]^\trianglelefteq\\
            &  & & \ms Q_m\times [\Delta^{m-1}] \ar@{<..}[u]^{\trianglelefteq \times \id}
 \ar@{^{(}->}[r]^{\id \times [d^{m}]} & \ms Q_m \times [\Delta^m]. \ar@{<..}[u]^\trianglelefteq}
\end{equation}
By applying (\ref{ind_subd}) we can fill the north-east corner of
(\ref{saw}). Thus $\ms{ T(Q)} \trianglelefteq \ms{Q}_m\times
[\Delta^m]$. In our situation, the construction of the induced assembly
commutes with the projection to the base. Thus we obtain a nontrivial
cellular structure $\ms{ T(Q)} \xar{\ms e (\ms Q)} [\Delta^m]$ on
the trivial fiber bundle $X\times \Delta^k \xar{\pi_2}\Delta^k$.
Figure~\vref{polycone00} presents an example of a 3-chain of geometric
assemblies on an interval. Figure~\vref{polycone11}  presents the
corresponding ball structure on $[0,1]\times \Delta^2$.

We can explicitly describe the ball structure $\ms{T}(\ms
Q)\xar{\ms e(Q)} [\Delta^m]$. Let $ S_i$ be the set of balls of
$\ms Q_i$.  For any subset $\mb k$ of the set $\bs m $, we
denote by $\Delta^{\mb{k}} \subseteq \Delta^m=\Delta^{\bs{m}}$ the
face of $\Delta^m$ generated by vertices with numbers from $\mb{k}$.
Then the ball complex $[\Delta^m]$ has the form $( \Delta^m,
\{\Delta^{\mb{k}}\}_{\mb k \in 2^{\bs{m}}})$. For any $\mb k \in
2^{\bs{m}}$, denote by $\max (\mb k) \in \bs m$ the maximal element of $\mb k$.
\begin{prop}\label{prizms}
The balls of $\ms{T}(\ms Q)$ are indexed by the pairs
$$( \mb k , B ) \text{ such that } \mb k \in 2^{\bs{m}} \text{ and } B \in S_{\max(\mb k)}$$
and have the form $T_{(\mb k , B)}=B\times \Delta^{\mb k} \subset
X\times \Delta^{m} $. The adjacency of balls is as follows:
$$ T_{( \mb k_0 , B_0 )} \us{X\times \Delta^m}{\subseteq} T_{( \mb k_1 , B_1 )} $$
if and only if
$$\mb k_0 \us{\bs m}{\subseteq} \mb k_1 \quad\text{and}\quad B_0 \us{X}{\subseteq} B_1.$$
The cellular morphism $\ms{e(Q)}$ projects the ball
 $T_{( \mb k , B )}$ onto the ball
$\Delta^{\mb k}$ of the base.
\end{prop}
The balls of $\ms{T (Q)}$ will be called \bem{prisms}. The
construction $\ms{T(-)}$ is functorial with respect to morphisms
of faces and degenerations of simplices of $\mc N \mb R(X)$.
Obviously, the correspondence $\ms Q \mapsto \ms{ e(Q)}$ is a
contravariant functor from the category $C\mc N \mb
R(X)$ of simplices of the simplicial set $\mc N \mb R(X)$ to the category of
cellular fiber bundles.
\subsection{The simplicial groupoid of prismatic
 homeomorphisms  $\mathrm{Prism}(X)$}
\label{prism_hom}\mbox{}

\subsubsection{} Let us fix some conventions and notations.
Let $B$ be a polyhedron. We will denote by $\PL_B(X)$ the group of
fiberwise homeomorphisms of the trivial bundle $X\times B
\xar{\pi_2}B$. If $A \xar{h} B$ and $f\in \PL_B(X)$, then we denote
by $f\lfloor_h \in \PL_A(X)$ the homeomorphism induced by $h$. If $h$
is an embedding, then we simply write $f\lfloor_A$. We regard the group $\PL(X)$
of $\PL$ homeomorphisms of $X$ as a simplicial group
with the group of $m$-simplices $\PL_m(X)=\PL_{\Delta^m}(X)$.
\subsubsection{}\label{ss322}
As above, let $\ms Q=(\ms Q_0\trianglelefteq_1{\ldots} \trianglelefteq_m \ms Q_m) $, and let
$\ms T(\ms Q)\xar{\ms{ e(Q)}} [\Delta^m]$ be the corresponding
prismatic decomposition of the trivial fiber bundle. We say that a
homeo\-mor\-phism $f \in \PL_m(X)$ is \bem{$\ms Q$-prismatic} if
for every ball  $T_{\mb k,B}$ of the complex $\ms {T(Q)}$
\begin{equation}\label{412}f^{-1}(T_{\mb (k,B)})=f^{-1}
\lfloor_{\max \mb k}(B)\times \Delta^{\mb k},\end{equation}
where by $X\xar{f^{-1}\lfloor_i} X$ we denote the fiber of the
homeomorphism $f^{-1}$ over the $i$th vertex of $\Delta^m$. Put
\begin{equation}\label{413} f^{-1} \ms Q=
(f^{-1}\lfloor_{ 0} (\ms Q_0)\trianglelefteq  f^{-1}\lfloor_{ 0} (\ms Q_1)
\trianglelefteq {\ldots}  \trianglelefteq f^{-1}\lfloor_{ 0} (\ms Q_m) ). \end{equation}
\begin{prop}\label{prop6}
A $\ms Q$-prismatic homeomorphism $f$ induces a fiberwise cellular
homeomorphism of the cellular bundles:
$$\bfig
\Vtriangle[\ms{T}(f^{-1} \ms Q)`\ms {T}(\ms Q)`{[\Delta^m]}.;f`\ms{e}(f^{-1} \ms Q)`\ms{e(Q)} ]
\efig$$
\end{prop}
\begin{proof}
From condition~(\ref{412}) it follows that
for every $i = 0,{\ldots} ,m$, $B \in \ms Q_i$,
$$f^{-1}\lfloor_i (B) = f^{-1}\lfloor_0(B).$$
Indeed, $f^{-1}(T_{(\{0,i\},B)}=B\times
\Delta^{\{0,i\}})=f^{-1}\lfloor_i(B)\times \Delta^{\{0,i\}}\equiv
f^{-1}\lfloor_0(B) \times \Delta^{\{0,i\}}$, where the last identity
is coordinatewise. Therefore
$$f^{-1}\lfloor_i (\ms Q_i)=f^{-1}\lfloor_0 (\ms Q_i)$$ and
$$f^{-1}\lfloor_i (\ms Q_i) = f^{-1}\lfloor_0 (\ms Q_i) \trianglelefteq f^{-1}\lfloor_0
 (\ms Q_j)=f^{-1}\lfloor_j (\ms Q_j) \text{ for } i \leq j.$$
Thus condition~(\ref{412}) can be rewritten as
  $$f^{-1}(T_{(B,\mb k)}(\ms Q))= T_{(f^{-1}\lfloor_0(B),\mb k)}(f^{-1} \ms Q),$$
i.e., as the assertion that
  $f^{-1}(T_{(B,\mb k)}(\ms Q))$ is a prism of $\ms T(f^{-1}\ms Q)$.
\end{proof}
\subsubsection{} \label{ss323}
Now we can define a simplicial groupoid
$$\mb{N}^{\op} \xar{\mr{Prism (X)}}\mb{Groupoids}.$$
The set $\Ob_m \mr{Prism} (X)$ of $m$-objects of this groupoid
is the
set of all $m$-simplices from $\mc N\mb R (X)$. The set of
$m$-morphisms
$$ (\mr{Prism}(X))(\ms{Q^0,Q^1})$$ is the set of all prismatic homeomorphisms
sending $\ms{T(Q^0)}$ to $\ms{T(Q^1)}$ according to construction from Proposition~\ref{prop6}
The faces and degenerations are defined in a natural way and
agree with the faces and degenerations in $\mc N \mb R(X)$. The
simplicial ``space'' of the groupoid $\mr{Prism}(X)$ is the simplicial
set $\Mor \mr{Prism} (X)$ of all morphisms of the groupoid
$\mr{Prism}(X)$. Thus the topological space of the groupoid is the
space $|\Mor \mr{Prism} (X)|$.

\begin{lemma}[\bem{On extension of a prismatic homeomorphism}]\label{Alex}\mbox{}\\
{\bf \rm (a)} Let  $\ms Q^0, \ms Q^1 \in  (\mc N \mb R(X))_0$,
and let $\mb P(\ms Q^0) \xar{\mu} \mb P (\ms Q^1)$ be an isomorphism. Then there exists
 a cellular isomorphism $\ms Q^0 \xar{f} \ms Q^1$ such that $\mb P(f)=\mu$.

\smallskip\noindent
{\bf \rm (b)} Let $\ms Q^0, \ms Q^1 \in (\mc N \mb R(X))_1$, and let
$\ms Q^0_0 \xar{f_0} \ms Q^1_0$ and $\ms Q^0_1 \xar{f_1} \ms Q^1_1$
be cellular ho\-meo\-mor\-phisms such that the diagram
$$\bfig \Square[\mb P (\ms Q_0^0)`\mb P(\ms Q^0_1)`\mb P (\ms Q_1^0)`\mb P
(\ms Q_1^1);\mb P(\trianglelefteq)`\mb P(f_0)`\mb P(f_1)`\mb P(\trianglelefteq)] \efig $$ is commutative.
Then there exists $f \in \Mor_2 \mathrm{Prism}(X)$, $\ms{Q}^0 \xar{f}
\ms{Q}^1$,
such that  $d_0 f=f_1$, $d_1 f=f_2$.

\smallskip\noindent
{\bf \rm (c)} Let $\ms Q^0, \ms Q^1 \in (\mc N \mb R(X))_m$ with $m
\geq 2$. Let  $f_0,{\ldots} ,f_m \in \Mor_{m-1} \mr {Prism} (X)$, $d_i \ms
Q^0 \xar{f_i} d_i \ms Q^1$, be a collection of prismatic
homeomorphisms such that $d_i f_j = d_{j-1}f_i$ whenever  $i<j$.
Then there exists
 $f \in \Mor_{m} \mr {Prism} (X)$,
$\ms Q^0 \xar{f} \ms Q^1$, such that $d_i f = f_i$ for $i=0,{\ldots} ,m$.
\end{lemma}
\begin{proof}
This lemma is a form of Alexander's trick. We will present a
detailed proof, because this lemma is important for our further
constructions.

1. We need the following observation, which is
standard in the theory of $\PL$ fiber bundles.

\medskip\noindent
\emph{Consider the space $\R^k\times \R^l$, the projection $\R^k\times
\R^l \xar{\pi_2}\R^l$,  and two geometric $d$-dimensional simplices
$S^0,S^1$ in $\R^k\times \R^l$: $S^0 = \conv (s^0_0,{\ldots} ,s^0_d)$, $S^1
= \conv (s^1_0,{\ldots} ,s^1_d)$. Assume that $\pi_2(s^0_i)=\pi_2(s^1)_i$,
$i=0,{\ldots} ,d$. Let $S^0\xar{A}S^1$ be an affine map such that
$A(s^0_i)=s^1_i$, $i=0,{\ldots} ,d$. Then $A$ is a $\pi_2$-fiberwise map,
i.e., the diagram
$$\bfig \Vtriangle[S^0`S^1`\R^l;A`\pi_2`\pi^2]\efig$$ is commutative.}

\medskip\noindent
This  follows from the fact that the map $A$ sends a point $x\in
S^0$, $x=\sum_{i=0}^d t_i s^0_i$, to the point  $A(x)=\sum_{i=0}^d t_i
s^1_i \in S^1$ and the calculation
$$\pi_2(A(x))=\sum_{i=0}^d t_i \pi_2 (s^1_i)=
\sum_{i=0}^d t_i \pi_2 (s^0_i) = \pi_2(x).$$

\smallskip
2. (Parametric  $\PL$ Alexander's trick.) \emph{Consider the projection
$\Delta^k \times \Delta^l \xar{\pi_2} \Delta_l$. Consider the
$(k+l-1)$-sphere $$\partial (\Delta^k \times \Delta^l) = (\partial
\Delta^k \times \Delta^l)\cup (\Delta^k \times \partial \Delta^l).$$
Assume that  a fiberwise homeomorphism
$$\bfig \Vtriangle[\partial ( \Delta^k \times \Delta^l)
`\partial (\Delta^k \times \Delta^l)`\Delta^l;F`\pi_2`\pi^2]\efig$$ is fixed.
Then there is a fiberwise homeomorphism
$$\bfig \Vtriangle[\Delta^k \times \Delta^l`\Delta^k \times \Delta^l`\Delta^l;
 G`\pi_2`\pi^2]\efig$$ such that
 $G|_{\partial ( \Delta^k \times \Delta^l)}=F$. }

First, we should triangulate the homeomorphism $F$. Second,
using the convexity ({\it sic}!) of the prism $\Delta^k \times \Delta^l$, we
can find stellar extentions of the triangulations in the preimage
and image to combinatorially isomorphic triangulations of
$\Delta^k
\times \Delta^l$. Such a stellar extention determines a fiberwise
extention of $G$ by Step~1 of the proof.

\smallskip
3. Step~2 of the proof implies a more general fact.
\emph{Assume that we have two balls $B^k, B^l$. Consider the projection $B^k
\times B^l \xar{\pi_2} B_l$. Consider the $(k+l-1)$-sphere
$$\partial (B^k \times B^l) =
(\partial B^k \times B^l)\cup (B^k \times \partial B^l).$$
Assume that a fiberwise homeomorphism $$\bfig \Vtriangle[\partial (
B^k \times B^l) `\partial (B^k \times B^l)`B^l;F`\pi_2`\pi^2]\efig$$
is fixed. Then there exists a fiberwise homeomorphism $$\bfig
\Vtriangle[B^k \times B^l`B^k \times B^l`B^l;
 G`\pi_2`\pi^2]\efig$$
such that $G|_{\partial ( B^k \times B^l)}=F$.}

For the proof, we choose two homeomorphisms $B^k \xar{H_1}
\Delta^k$ and $B^l \xar{H_2} \Delta^l$, then apply a homeomorphism $B^k \times B^l \xar{H_1 \times H_2} \Delta^k
\times \Delta^l$ to $B^k \times
B^l$, and refer to Step~2 of the proof.

\smallskip
4. Let $\partial [\Delta^m]$ be the ball complex of the boundary of the
standard simplex. For $(\bs m \xar{\ms Q} \mb R (X)) \in (\mc N \mb
R (X))_m $, denote by $\theta_{(-1)}(\ms Q)$ the cellular fiber bundle
on $\partial [\Delta^m]$  induced from $\ms e(\ms Q)$ by the embedding
$\partial [\Delta^m] \hra{} [\Delta^m]$:
$$\bfig
\square/^{(}..>`^..>`>`^{(}->/[\Theta_{(-1)}(\ms Q)`\ms T (\ms Q)`{\partial [\Delta^m]}`{[\Delta^m].};`\theta_{-1}(\ms Q)` \ms e(\ms Q)`]
\efig$$
The balls of $\Theta_{(-1)}$ are of the form $T_{(\mb k, B)} \in \ms
T (\ms Q)$, $\mb k \neq \bs m$. Denote by $\Theta_{(i)}(\ms Q)$ the
subcomplex of $\ms T (\ms Q)$ that consists of the balls
$\Theta_{(-1)}(\ms Q)$ together with all balls of the form
$T_{(\bs m, B)}$, $B \in \ms Q_m$, $\dim B \leq i $. Recall that $n=\dim
X$. We obtain the filtration
$$\xymatrix{ \Theta_{(-1)}(\ms Q) \ar@{^{(}->}[r] \ar[d]_{\theta_{(-1)}(\ms Q)}&
 \Theta_{(0)}{(\ms Q)} \ar[d]_{\theta_{(0)}({\ms Q})} \ar@{^{(}->}[r] & \ar@{..>}[dl] \cdots \ar@{^{(}->}[r]
   & \ar[dll]^{\theta_{(n)}(\ms Q)=\ms e (\ms Q)} \Theta_{(n)}(\ms Q)=\ms T (\ms Q) \\
             \partial [\Delta^m] \ar@{^{(}->}[r] & [\Delta^m].& &
}$$
The restriction of  $\ms e (\ms Q)$ to $\Theta_{(i)}(\ms Q)$ will be
denoted by $\theta_{(i)} (\ms Q)$.

\smallskip
5. Now, using Step~3 of the proof, we can build $f$ inductively. The
maps $f_0,{\ldots} ,f_m$ are pasted
together to form a homeomorphism
$F^{(-1)}$ that coincides with $f_i$ restricted to the bundle
over the $i$th face of the sphere $\partial [\Delta^m]$:
$$\bfig \Vtriangle[\Theta_{(-1)} \ms Q^0
` \Theta_{(-1)} \ms Q^1 `\partial {[\Delta^m].};F^{(-1)}%
`\theta_{(-1)}(\ms Q^0)`\theta_{(-1)} e(\ms Q^1)]\efig$$
Then necesarily
$f_i|_0=f_j|_0=F^{(-1)}|_0$ for
 $i,j = 1,{\ldots} ,m$ and
$\ms Q^1_i = F^{(-1)}|_0 (\ms Q^0_i)$. The homeomorphism $F^{(-1)}$
sends a ball $T_{(B, \mb k)}$ of the complex $\Theta_{(-1)} \ms Q^0$
to the ball $T_{(F^{(-1)}|_0(B),\mb k)}$ of the complex
$\Theta_{(-1)} \ms Q^1$.

{\em Inductive step.} Assume that we have a prismatic homeomorphism
$\Theta_{i-1} (\ms Q^0) \xar{F^{(i-1)}} \Theta_{i-1}(\ms Q^1)$
sending a ball $T_{(B, \mb k)}$ of the complex $\Theta_{(i-1)} \ms
Q^0$ to the ball $T_{(F^{(-1)}|_0(B),\mb k)}$ of the complex
$\Theta_{(i-1)} \ms Q^1$. Let us extend this homeomorphism to a
fiberwise cellular homeomorphism
$$\bfig \Vtriangle[\Theta_{i} \ms Q^0
` \Theta_{i} \ms Q^1 `{[\Delta^m].};F^i`\theta_{i}(\ms Q^0)%
`\theta_{i} e(\ms Q^1)]\efig$$
To achieve this aim, we must extend $F^{(i-1)}$ to the balls of the
form $T_{(\bs m, B)}$ with  $\dim B = i$. Such an extension
exists, because  $\partial T_{(\bs m, B)} \subset \Theta_{(i-1)}(\ms
Q^0)$ and
 $F^{(i-1)}|_{\partial T_{(\bs m, B)}}$ satisfies the conditions of Alexander's trick
 in the form given in Step~3 of the proof.
\end{proof}

\subsection{The simplicial sets $\mc W\mr{Prism}(X)$ and $\ol{\mc W} \mr{Prism} (X)$}
\label{wprism}
\subsubsection{}\label{w_const}
We recall the so-called $W$-construction, i.e., the Eilenberg--MacLane
construction of the principal simplicial bundle $EG\xar{u_G}BG$ for
a simplicial group $G$.

Put
\begin{equation*}
(\ol{\mc W}G)_m = \left\{
\begin{array}{ll}
G_{m-1}\times G_{m-2}\times\cdots\times G_0 & \text{ for } m>0,  \\
\{*\} \text{ (one-element set) }& \text{ for } m=0.
\end{array} \right.
\end{equation*}
For $m>0$, denote an element of
$(\ol{\mc W} G)_m$ by $[g_{m-1},{\ldots} ,g_0]$.
We define the faces of one-dimensional simplices and the
degeneracy of the 0-dimensional simplex as follows:
\begin{equation}\label{w_const1}
\begin{array}{ll}
d_i([g_0])=* & \text{ for  } i=0,1, \\
s_0(*)=[e_0] & (e_0\text{ is the identity of the group } G_0).
\end{array}
\end{equation}
For $m>1$ the faces are defined as follows:
\begin{equation}\label{w_const_2}
d_i[g_{m-1},{\ldots} ,g_0]= \left\{
\begin{array}{ll}
{[}g_{m-2},{\ldots} ,g_0 {]}&\text{ for } i=0, \\
{[}d_{i-1}g_{m-1},{\ldots} d_1g_{m-i+1},\us{m-i-1}
{(g_{m-i-1}\circ d_0g_{m-i})},g_{m-i-2},{\ldots}, g_0 {]} &\text{ for } i=1,{\ldots} ,m-1, \\
{[}d_{m-1}g_{m-1},{\ldots} , d_1g_1{]} &\text{ for } i=m;
\end{array} \right.
\end{equation}
the degeneracies are defined by
$$
s_i [g_{m-1},{\ldots} ,g_0] = \left\{
\begin{array}{ll}
{[}e_m,g_{m-1},{\ldots} ,g_0{]} & \text{ for } i=0, \\
{[}s_{i-1}g_{m-1},{\ldots} ,s_0g_{m-i},e_{m-i},g_{m-i-1},{\ldots} ,g_0{]} & \text{ for }i=1,{\ldots} ,m,
\end{array} \right.
$$
where $e_j$ is the identity of the group $G_j$.

Put $(\mc WG)_m=G_m\times\cdots \times G_0$. An element of
$(\mc WG)_m$ will be denoted by $\la
g_m,{\ldots} ,g_0\ra$. The faces are
$$
d_i\la g_{m},{\ldots} ,g_0 \ra= \left\{
\begin{array}{ll}
{\la}d_{i}g_{m},{\ldots}, d_1g_{m-i+1},\us{m-i-1}{(g_{m-i-1}\circ d_0g_{m-i})},g_{m-i-2},%
{\ldots} , g_0 {\ra} &\text{ for } i=0,{\ldots} ,m-1, \\
{\la}d_{m}g_{m},{\ldots} , d_1g_1{\ra} &\text{ for } i=m.
\end{array} \right.
$$
The degeneracies are
$$
s_i\la g_{m},{\ldots} ,g_0 \ra= \la s_i g_m, s_{i-1}g_{m-1},{\ldots}
,s_0g_{m-i},e_{m-i},g_{m-i-1},{\ldots} ,g_0\ra.
$$
On $\mc WG$ we define a free action  $G\times \mc WG \xar{}\mc WG $ of the
group  $G$ by the rule $(h_m,\la g_m,{\ldots} ,g_0\ra) \mapsto
\la h_m g_m,g_{m-1},{\ldots} , g_0\ra$. The principal fiber bundle
$$\mc W G \xar{u_G} \ol{\mc W}G: \la g_m,{\ldots} ,g_0
\ra \mapsto [g_{m-1},{\ldots} ,g_0] $$
corresponding to this action is a universal principal $G$-bundle
for $G$.

\subsubsection{}
Let us develop a version of the $\mc W$-construction for the simplicial
groupoid of prismatic homeomorphisms. Let $g \in \Mor_m
\mr{Prism}(X)$ be a prismatic homeomorphism. Denote by $\dom(g)\in
\mc N_k \mb R(X) $ its image and by
 $\codom (g) \in \mc N_k \mb R(X)$ its preimage.

For $m\geq 1$, the set of $m$-simplices of
$\ol{\mc W} \mr{Prism}(X)$ is the set of all pairs
$$(\ms Q,[ g_{m-1},{\ldots} ,g_0]),$$
where $\ms Q \in \mc N_m \mb R (X)$ and $[ g_{m-1},{\ldots} ,g_0]$ is a
sequence of prismatic homeomorphisms such that $g_i\in \Mor_i%
 \mr{ Prism}(X)$, $i=0,{\ldots} ,m-1$, $d_0 \ms Q = \dom g_{m-1}$,
 and $d_0 \codom g_i = \dom g_{i-1}$. For $m=0$  the simplices are
identified with the elements of $\mb R(X)$.

We define the faces
of 1-simplices and the
degeneracy of the 0-simplex in $\ol{\mc W} \mr{Prism}(X)$ by the formulas
\begin{equation}
\begin{array}{ll}
d_0(\ms Q, [g_0])=\codom g_0,  \\
d_1 (\ms Q, [g_0])= d_1 \ms  Q, \\
s_0(\ms Q)=(s_0\ms Q,[e_0]) & (e_0\text{ is the identity of the group } G_0).
\end{array}
\end{equation}

For $m>1$ the faces and degeneracies are defined as follows:
\begin{multline}\label{w_prism1}
d_i(\ms Q, [g_{m-1},{\ldots} ,g_0])\\
=\left\{
\begin{array}{ll}
(\codom g_{m-1}, {[}g_{m-2},{\ldots} ,g_0 {]})&\text{ for } i=0, \\
(d_i \ms Q, {[}d_{i-1}g_{m-1},{\ldots} ,d_1g_{m-i+1},\us{m-i-1}{(g_{m-i-1}\circ d_0g_{m-i})},g_{m-i-2},{\ldots} , g_0 {]}) &\text{ for } i=1,{\ldots} ,m-1, \\
(d_m \ms Q, {[}d_{m-1}g_{m-1},{\ldots} , d_1g_1{]}) &\text{ for } i=m
\end{array} \right.
\end{multline}
and
\begin{multline}\label{w_prism2}
s_i (\ms Q, [g_{m-1},{\ldots} ,g_0]) \\ = \left\{
\begin{array}{ll}
(s_0 \ms Q, {[}e_m,g_{m-1},{\ldots} ,g_0{]}) & \text{ for } i=0, \\
(s_i \ms Q, {[}s_{i-1}g_{m-1},{\ldots} ,s_0g_{m-i},e_{m-i},g_{m-i-1},{\ldots} ,g_0{]}) &
\text{ for }i=1,{\ldots} ,m,
\end{array} \right.
\end{multline} where $e_j$ is the identity of the group $G_j$.

Now we define a simplicial set $\mc W \mr{Prism} (X)$.

The elements of $(\mc W \mr{Prism} (X))_m$ are all pairs $(\ms Q,
\la g_m,{\ldots} ,g_0\ra)$, where $\ms Q \in \mc N_m \mb R (X)$, $\codom
g_m = \ms Q$, $d_0 \ms Q = \dom g_{m-1}$, $d_0 \codom g_i = \dom
g_{i-1}$ for $i=1,{\ldots} ,m-1$.

The faces are
\begin{multline*}
d_i(\ms Q, \la g_{m},{\ldots} ,g_0 \ra) \\ = \left\{
\begin{array}{ll}
(\codom g_{m-1}, \la g_{m-1}d_0 g_m, g_{m-2},{\ldots} , g_0 \ra) & \text{ for } i=0, \\
(d_i \ms Q, {\la}d_{i}g_{m},{\ldots} ,d_1g_{m-i+1},\us{m-i-1}
{(g_{m-i-1}\circ d_0g_{m-i})},g_{m-i-2},{\ldots} , g_0 {\ra}) &\text{ for } i=0,{\ldots} ,m-1, \\
(d_m Q, {\la}d_{m}g_{m},{\ldots} , d_1g_1{\ra}) &\text{ for } i=m.
\end{array} \right.
\end{multline*}
The degeneracies are
$$
s_i(\ms Q, \la g_{m},{\ldots} ,g_0 \ra)= (s_i \ms Q, \la s_i g_m, s_{i-1}g_{m-1},{\ldots}
,s_0g_{m-i},e_{m-i},g_{m-i-1},{\ldots} ,g_0\ra).
$$
Define a morphism of simplicial sets
$$\mc W \mr{Prism} (X) \xar{u_{\mr{Prism} (X)}}  \ol{\mc W} \mr{Prism} (X) $$
by the formula $(\ms Q, \la g_m,{\ldots} ,g_0 \ra) \mapsto (\ms Q,
[g_{m-1},{\ldots} ,g_0])$.

The following simple fact is important for us.
\begin{prop} \label{pr7}
We have the following commutative square of maps of simplicial sets:
\begin{equation*}
\bfig
\Square[\mc W \mr{Prism} (X)`\mc W \PL (X)`\ol{\mc W} \mr{Prism} (X)`\ol{\mc W} \PL (X),;\Pi'
`u_{\mr{Prism} (X)}`u_{\PL(X)}`\Pi]
\efig
\end{equation*}
where the horizontal arrows forget the combinatorics of
objects of the groupoid:
$$(\ms Q, \la g_m,{\ldots} ,g_0 \ra) \overset{\Pi'}{\mapsto} \la g_m,{\ldots},g_0 \ra,$$
$$(\ms Q, [ g_{m-1},{\ldots} ,g_0]) \overset{\Pi}{\mapsto} [ g_{m-1},{\ldots} ,g_0 ].$$
\end{prop}
Formulas~(\vref{w_const1}) and (\vref{w_const_2}), which define the $\mc
W$-construction, express the faces of  $\ol{\mc W}
G$ via the faces of  $G$ and the degeneracies of  $\ol{\mc W}G$ via
the degeneracies of $G$. We need to mention the existence of
the inverse expressions. Let $w = (\ms Q(w),[g_{m-1}(w),{\ldots} ,g_0(w)])$ be
an $m$-simplex of
 $\mc W \mr{Prism}(X)$.
We can deduce the following expressions.
\begin{prop}\label{10}
$$\begin{array}{rll}
   g_{i}(w)  &  = g_i(d_0 w)   & \text{ for } i=1,{\ldots} ,m-2, \\
  d_j g_{m-1}(w) & = g_{m-2}(d_j w) & \text{ for } j=1,{\ldots} ,m-1, \\
  d_0 g_{m-1}(w) & = (g_{m-2}(d_0 w))^{-1}\circ g_{m-2}(d_1 w), &
\end{array}$$
where $d_* $ on the left is a face in $G$ and  $d_*$ on the right is a
face in $\ol{\mc W} \mr{Prism}(X)$.
\end{prop}
\subsection{The embedding $\mc N \wt{\bg R} (X) \xar{\Psi} \ol{\mc W} \mr{Prism} (X)$
and the projection  $\ol{\mc W} \mr{Prism} (X) \xar{\Psi^*} \mc N \wt{\bg R} (X)$}\mbox{}
We will build an
embedding $\mc N \wt{\bg R} (X) \xar{\Psi} \ol{\mc W} \mr{Prism}
(X)$ by induction on the skeletons of the simplicial set $\mc N \wt{\bg R}
(X)$, using a sequential choice of prismatic homeomorphisms. The skeletons
of simplicial sets are discussed in detail, for example, in
\cite[Chap.~V.~1]{GJ}.
The skeleton $\mr{sk}_m X \subset X$  is the
simplicial subset generated by all nondegenerate simplices of
dimension at most $m$.
\subsubsection{}
First, for each $m$-simplex
$$\mf{Q} = (\ms Q_0 \sar{\mu_1}{\ldots} \sar{\mu_m}\ms Q_m)
\in (\mc N \wt{\bg R} (X))_m$$
we prepare its linearization
$$\ms L \mf{Q} = (\ms L_0 \mf{Q} \trianglelefteq{\ldots}  \trianglelefteq \ms L_m \mf{Q}) \in (\mc N \mb R (X))_m. $$
By the definition of an abstract assembly (Sec.~\vref{2.7}), any abstract
assembly $\ms Q_{i-1}\sar{\mu_i}\ms Q_i$ of
the chain $\mf{Q}$
is representable by some $\bc R (X)$-morphism $\ms Q_{i-1}\xar{f_i}\ms
Q_i$. Let us fix these representatives for all $i$. We obtain a chain
$$Q=(\ms Q_0 \xar{f_1}{\ldots} \xar{f_m}\ms Q_m) \in (\mc N \bc R (X))_m$$
such that $\bc A (Q) =  \mf{Q}$. Applying the decomposition~(\vref{h-r}),
we obtain the following commutative diagram in $\bc R (X)$:
\begin{equation}
\xymatrix{\ms{\mf{Q}}_0=\ms L_0 \mf{Q}\ar[dr]_{f_1} \ar[r]^{\trianglelefteq}&\ms L_1  \mf{Q}\ar[r]^{\trianglelefteq}\ar[d]& \ms L_2  \mf{Q}\ar[r]^{\trianglelefteq}\ar[d]&\cdots\ar[r]^{\trianglelefteq}& \ms L_m  \mf{Q}\ar[d]\\
           & \ms \mf{Q}_1=\ms L_0 d_0  \mf{Q}\ar[dr]_{f_2}\ar[r]^{\trianglelefteq}&\ms L_1 d_0  \mf{Q}\ar[r]^{\trianglelefteq}\ar[d]&\cdots\ar[r]^{\trianglelefteq}& \ms L_{m-1} d_0  \mf{Q}\ar[d]\\
           & &\ms \mf{Q}_2 = \ms L_0d_0^2  \mf{Q} \ar[dr]_{f_2}\ar[r]^{\trianglelefteq} &\cdots\ar[r]^{\trianglelefteq}&\ms L_{m-2} d_0^2  \mf{Q}\ar[d]\\
           & & &\cdots \ar[dr]_{f_m}&\vdots\ar[d]\\
           & & & & \ms \mf{Q}_m=\ms L_0d_0^m  \mf{Q},}
\end{equation}
where the rows are chains of geometric assemblies starting from
$\ms Q_0,{\ldots} ,\ms Q_m$. The vertical morphisms are cellular
homeomorphisms. The rows are independent from the choice of morphisms $f_i$
representing $\mu_i$, they depend only on $\ms Q_i$. The upper row
is the chain $\ms L \mf{Q}$ canonically associated to
$\mf{Q}$. Note that if $i>0$, then $\ms L d_i \mf{Q} = d_i \ms L
\mf{Q}$.

\subsubsection{} \label{ss342}
Pick $\mf{Q}\in (\mc N \wt{\bg R} (X))_m$. We will look for
$\Psi(\mf{Q})\in (\mc W\mr{Prism}(X))$ in the form $(\ms L
\mf{Q},[g_{m-1}(\mf{Q}),{\ldots} , g_0(\mf{Q})])$, where the
homeomorphisms $g_i(\mf{Q})$ are constructed by induction on the
skeletons $\mr{sk}_m \mc N \wt{\bg R} (X) \xar{\Psi_m}\mr{sk}_m
\mr{Prism}(X)$. A nondegenerate simplex in $\mc N \wt{\bg R} (X)$ is a
chain of abstract assemblies that has no identity assemblies. The
simplices of $\mr{sk}_m \mc N \wt{\bg R} (X)$ are all chains of
abstract assemblies that have at most $m$ nonidentities. Thus
$\mr{sk}_0 \mc N \wt{\bg R} (X)$ is the simplicial set whose typical
$k$-simplex is a chain
$$\us{k}{\underbrace{\ms Q \sar{\id} \ms Q\sar{id} {\ldots}  \sar{\id}}}\ms Q$$ consisting
of $k$ identical
abstract assemblies of a ball complex $\ms Q \in
\mb R (X)$. To this simplex we assign the degenerate simplex
$$(\ms Q, [e_{m-1},{\ldots} ,{e}_0]) \in  \ol{\mc W} \mr{Prism}(X).$$
Assume that $\Psi_{m-1}$ is already constructed. Let us construct a
map $\Psi_m$ extending $\Psi_{m-1}$. Let $\mf{Q}\in (\mr{sk}_m \mc N
\wt{\bg R} (X))_m$ be a nondegenerate simplex. Then $d_i
\mf{Q} \in \mr{sk}_{m-1} \mc N \wt{\bg R} (X)$  and $\Psi d_i
\mf{Q}$ are already constructed. We search for $\Psi(\mf{Q})$
in the form $(\ms L \mf{Q},[g_{m-1}(\mf{Q}),{\ldots} ,g_0(\mf{Q})])$. The
equalities $g_i(\mf{Q})=g_i(d_0 \mf{Q})$, $i=m-2,{\ldots} ,0$ are
already satisfied. Therefore $g_i$ are already constructed for
$i=0,{\ldots} ,m-2$. We must define  $g_{m-1} (\mf{Q})$. According to
Proposition~\vref{10}, the maps $d_i g_{m-1}(\mf{Q})$, $i=1,{\ldots} ,m-1$,
are already constructed and satisfy the requirements of
Proposition~\vref{Alex}. Therefore there exists a prismatic
homeomorphism $g_{m-1}(\mf{Q})$ that is an extension of $d_i
g_{m-1}(\mf{Q})$, $i=1,{\ldots} ,m-1$. For the degenerate simplex $\mf{Q}$,
we can construct $(\ms L \mf{Q},[g_{m-1}(\mf{Q}),{\ldots} ,g_0(\mf{Q})])$
according to formulas~(\vref{w_prism2}) with the help of the unique
reduction to a nondegenerate
 simplex: $\mf{Q}=s_{i_1}s_{i_2}{\ldots}  \mf{Q}'$, where $\mf{Q}'$ is a nondegenerate simplex
 of  dimension at most $m$.
\subsubsection{}
Let us build a projection $\ol{\mc W} \mr{Prism} (X) \xar{\Psi^*}
\mc N \wt{\bg R} (X)$. For a fiberwise homeomorphism $f \in
\PL_m(X)$, we denote by $ f\lfloor_i$ its fiber over the $i$th vertex
of the base. To a simplex $w=(\ms Q, [g_{m-1},{\ldots} ,g_0])\in
(\ol{\mc W} \mr{Prism} (X))_m$ we associate the sequence of maps $\wh
g_i(w)= g_i\lfloor_i\in \PL_0(X)$, $i=0,{\ldots} ,m-1$. Consider the
sequence of ball complexes
$$\ms Q_0 = \wh{\ms Q}_0(w),\;\wh {\ms Q_i} (w)=\wh g_{m-i}\circ \wh g_{m-i+1}\circ{\ldots} \circ
\wh g_{m-1}(w)(\ms Q_i),\quad i=1,{\ldots} ,m.$$
Due to the prismaticity of the homeomorphisms $g_i$, we obtain a chain
$\wh{\Psi^*}(w) \in (\mc N\mb R (X))_m$:
$$ \wh{\Psi^*}(w) =  \wh{\ms Q}_0(w) \xar{\wh g_{m-1}(w)}\wh{\ms Q}_1(w)\xar{\wh g_{m-2}(w)}
{\ldots} \xar{\wh g_0(w)} \wh{\ms Q}_m(w). $$
Set the value
$\Psi^*(w) \in (\mc N \wt{\bg R}(X))_m $ to be the
image of the chain $\wh{\Psi^*}(w)$ under the functor
$\bc{A}: \bc{R}(X) \xar{} \wt{\bg R}(X)$. One can verify that the
definition is correct, $\Psi^*(w)$ is a map of simplicial sets,
$\Psi^*\Psi=\id$, and the following proposition holds.
\begin{prop}\label{psi}
The map $\Psi\Psi^*$ is simplically homotopic to the identity.
\end{prop}
\begin{proof}
A simplicial homotopy is constructed by induction on the skeletons using
Lemma~\vref{Alex}.
\end{proof}
\subsection{The maps $\Phi, \Phi^*$}
We define an embedding $\mc N \mb R (X) \xar{\Phi} \mc W \mr{Prism}(X)$ by
a correspondence on simplices: to an $m$-chain  $\ms Q=(\ms Q_0
\trianglelefteq{\ldots} \trianglelefteq \ms Q_m)$ we associate the $m$-simplex $\Phi(\ms Q
)=(\ms Q,[e_m,{\ldots} ,e_0])$. We define a map $\mc W \mr{Prism}(X)
\xar{\Phi^*}  \mc N \mb R (X)$ by a correspondence on simplices:
to an $m$-simplex $w=(\ms Q,[g_m,{\ldots} ,g_0])$ we associate the $m$-chain
$\Phi^*(w)=g_m^{-1}\ms Q=( g_m^{-1}\lfloor_0\ms Q_0
 \trianglelefteq{\ldots} \trianglelefteq  g_m^{-1}\lfloor_m\ms Q_m)$.
\begin{prop} \label{pr10}
The composition $\Phi^*\Phi$ is the identity, the map
$\Phi\Phi^*$ is simplically homotopic to the identity.
\end{prop}
\begin{proof}
The first assertion follows from the construction, the second one can be
proved by induction on the skeletons using Lemma~\vref{Alex}.
\end{proof}
As a result of the constructions of this section, we obtain the following
commutative diagram:
$$
\bfig
\Square<300>[ \mc N \mb R (X)` \mc W \mr{Prism} (X)`\mc N \wt{\bg R} (X)`\ol{\mc  W }\mr{Prism} (X),;
\Phi `\mc N \mb A `u_{\mr{Prism} (X)}`\Psi]
\efig
$$
where $\Phi, \Psi$ are
simplicial homotopy equivalences.

\subsection{Reduction of  Theorem~A
to the lemma on prismatic trivialization} According to  Proposition~\vref{pr7},
we have the commutative square
\begin{equation*}
\bfig
\Square[\mc W \mr{Prism} (X)`\mc W \PL (X)`\ol{\mc W} \mr{Prism} (X)`\ol{\mc W} \PL (X).;\Pi'
`u_{\mr{Prism} (X)}`u_{\PL(X)}`\Pi]
\efig
\end{equation*}
Consider the principal $\PL$ bundle
$$\wt{\mc W} \mr{Prism} (X) \xar{} \ol{\mc W} \mr{Prism}(X)$$
 induced by the map $\Pi$.
We obtain the commutative diagram
\begin{equation}\label{dia}
\xymatrix{ &\mc N \mb R (X) \ar[dl]_{\mc N \mb P}  \ar[r]^{\Phi} \ar[d]^{\mc N \mb A} & \mc W \mr{Prism} \ar[dr] (X) \ar@{^{(}->}[r]^i &  \wt{\mc W} \mr{Prism} (X)\ar[d] \ar[r]^{\wt \Pi} & \mc W \PL(X) \ar[d]\\
\mc N\bg R &\ar[l]_{\mc N\mb F}\mc N \wt{\bg R} (X)\ar[rr]^{\Psi} &  &\ol{\mc  W }\mr{Prism} (X) \ar[r]^{\Pi}& \ol{\mc W} \PL(X),
}
\end{equation}
where $\wt \Pi$ is the
pullback of $\Pi$, $\Pi'=\wt \Pi \circ i$, and the
left triangle is the nerve of the triangle~(\vref{tria}).

\begin{lemma}[\textbf{On prismatic trivialization}] \label{CombFrag}
The embedding $$|\mc W \mr{Prism}(X)|\hra{|i|} |\wt{\mc W}
\mr{Prism} (X)|$$ is a homotopy equivalence.
\end{lemma}
In Sec.~\ref{cm} Lemma~\ref{CombFrag} will be deduced from
Lemma~\ref{CombFrag1}. And all that follows will be mainly devoted
to the proof of  Lemma~\ref{CombFrag1}.

\subsection{Proof of Theorem A} \label{Aproof}
The homotopy equivalence $B\bg R(X)\approx B\PL(X)$ will be a
geometric realization of the bottom chain of simplicial maps
 in the diagram~(\ref{dia}). By definition, $B\bg R (X)\approx |\mc N \bg R (X)|$ at the
 left end of the chain; by the Eilenberg--MacLane theorem,
 $B\PL(X)\approx |\ol{\mc W}\PL(X)|$ at the right end of the chain. We must check that all intermediate elements of
 the chain are homotopy equivalences.

1. The following assertion is a standard fact of the theory
of simplicial principal bundles.

\medskip\noindent
{\em Let $X\xar{f}BG$ be a simplicial map, and assume that $EG_f\xar{{u_G}_f}X$ is the
induced principal bundle, i.e., the following  square is
Cartesian:
$$
\bfig
\Square[EG_f`EG`X`BG.;\tilde f`{u_G}_f`u_G`f]
\efig
$$
In this situation, if the space $|EG_f|$ is contractible, then in the
square
$$
\bfig
\Square[|EG_f|`|EG|`|X|`|BG|;|\tilde f|`|{u_G}_f|`|u_G|`|f|]
\efig
$$
the maps $|f|, |\tilde f| $ are homotopy equivalences.}

\medskip
This fact
follows from the equivariance of $\tilde f$, the possibility to
extend the square to a morphism of
the exact homotopy sequences of fibrations, the Whitehead theorem, and the 5-lemma.

\smallskip
2. By  Proposition~\vref{pr10}, the map $\Phi$ is a homotopy
equivalence. Therefore,  by  Lemma~\vref{lem2},  the space $|\mc W
\mr{Prism}(X)|$ is contractible. By  Lemma~\ref{CombFrag}, the space $\wt{\mc
W}\mr{Prism}(X)$ is contractible. Thus
we find ourselves in the situation
of Step~1 of the proof. Therefore $|\Pi|$ is a homotopy equivalence.
By Proposition~\vref{psi}, the map $|\Psi|$ is a homotopy equivalence;
and $|\mc N \mb F|$ is a homotopy equivalence because $\mb F$ is an
equivalence of categories by definition. This completes the proof
of Theorem~A. $\Box$

\section{Reduction of the lemma on prismatic trivialization to the lemma
on a common $\mb R(X)$-triangulation of fiberwise homeomorphisms}
\label{cm} It this section, we will deduce  Lemma~ \ref{CombFrag} ``on
prismatic trivialization'' from  Lemma~\ref{CombFrag1} ``on a common
$\mb R(X)$-triangulation of fiberwise homeomorphisms,'' which we will
formulate now.

Let $N \geq 0$. We define a simplicial set $\underline{\mr{Prism}}^N
(X)$. For $N \geq 1$, the $k$-simplices of
$\underline{\mr{Prism}}^N (X)$ are the data sets $  \la \ms Q,
\{g_i\}_{i=1}^N \ra $, where $\ms Q \in \mc N_k \mb R (X)$
and for every $i=1,{\ldots} ,N$
$$g_i \in \Mor_k \mr{Prism}(X),\quad \mr{codom} g_i = \ms Q. $$
For $N=0$, we assume that $\underline{\mr{Prism}}^0 (X)$
coincides with $\mc N \mb R (X)$. For $N \geq 1$, by forgetting
the prismaticity of $g_N$, we obtain an embedding $
\underline{\mr{Prism}}^N (X) \hra{j} \underline{\mr{Prism}}^{N-1}(X)
\times \PL (X)$.
\begin{lemma}[On a common $\mb{R}(X)$-triangulation
of fiberwise homeomorphisms] \label{CombFrag1}
The embedding $ |\underline{\mr{Prism}}^N (X)| \hra{|j|}
|\underline{\mr{Prism}}^{N-1}(X) \times \PL (X)|$ is a homotopy
equivalence.
\end{lemma}
We will prove Lemma~\ref{CombFrag1} in the subsequent sections. In this
section, after some preparations, we will prove the
implication Lemma~\ref{CombFrag1} $\Rightarrow $
Lemma~\ref{CombFrag} (in Sec.~\vref{lem6}).
\subsection{Injective simplicial sets
(or ``$\Delta$-sets'')}\label{inj} Let $\mb{iN}$ be the subcategory of
the category $\mb N$ of finite ordinals that has the same objects but
only injective monotone maps as morphisms. The morphisms of $\mb{iN}$
are generated by the cofaces. An \bem{injective simplicial set} is a
functor  $\mb{iN}^{\op} \xar{}\mb{Sets}$. (This is the same as a
``$\Delta$-set''
 of Rourke and Sanderson \cite{RSII}; we just wish to avoid another use of the
overloaded symbol $\Delta$.) The embedding $\mb{iN} \hra{\mb i} \mb{N}$
generates a forgetful functor $$\mb{Sets}^{\mb{N}^{\op}}
\xar{\mb{d}}
 \mb{Sets}^{\mb{iN}^{\op}},$$ which assigns to a simplicial set its injective part.

A simplicial complex $K$ is \bem{locally ordered} if there is a
partial order on the vertices of $K$ such that the vertices of every
simplex of $K$ are totally ordered. Denote by $\bs i (K)$ the natural structure
of an injective simplicial set on the simplices of $K$.
\subsection{Colorings of polyhedra by simplicial sets
and Brown's theorem on representable functors}\label{color} A
triangulation $|K|\xar{t}N$ of a polyhedron $N$ is \bem{locally
ordered} if the complex $K$ is locally ordered. Let $Y$ be a
simplicial set and $N$ be a polyhedron. A \bem{coloring of $N$ by $Y$}
(or an \bem{$Y$-coloring of $N$})
is a pair $\la t, f \ra$, where $K\xar{t}N$ is a
locally ordered triangulation and $\bs i (K)\xar{f}\mb d (Y)$ is a
morphism of injective simplicial sets. Two $Y$-colorings $\la t_0,
f_0 \ra$ and $\la t_1, f_1 \ra$ of a polyhedron $N$ are
\bem{concordant} if there is an $Y$-coloring of $N\times[0,1]$ that
coincides with $\la t_j, f_j \ra$ on the face $N\times\{j\}$. Let
$$\mb{PL}^\op \xar{\mb{Col} Y} \mb{Sets}$$ be the functor that
associates to a polyhedron $N$ the set of all concordance classes of
$Y$-colorings of $N$. Let $\mb{PL}^\op \xar{\mb{Ho}(-,|Y|)} \mb{Sets}$
be the functor that associates to a polyhedron $N$ the set $\mb
{Ho}(N,|Y|)$ of all homotopy classes of maps from $N$ to $|Y|$.

The relative Zeeman theorem on simplicial approximation \cite{Ze}
guarantees that there is a natural isomorphism of functors $\mb{Col}
Y\approx \mb{Ho}(-,|Y|)$. If $|Y|$ is connected and has countable
homotopy groups, then Brown's theorem \cite[p.~469]{Br} on
representability of homotopy functors on polyhedra guarantees that
$Y$ is determined by $\mb{Ho}(-,|Y|)$ (and hence also by $\mb{Col} Y$)
in a homotopy unique way.

Colorings of polyhedra by simplicial sets is a natural tool for
proving the homotopy triviality of pairs for non-Kan simplicial
sets, since simplicial homotopy works poorly in this case.

\subsection{The ``$\pmb{ \frown}$-lemma'' and ``$\pmb{\sqcap}$-lemma''}
Let $(X,Y)$ be a pair of simplicial sets, i.e., $Y\subseteq X$. We
always suppose that $X$ is connected. Consider the
embedding $|Y|\hra{f}
|X|$. We need to fix two different ways to say that $f$ is a
homotopy equivalence in terms of $(X,Y)$-colorings of
pairs of polyhedra. The difference between these two ways will be pure
combinatorial. The corresponding assertions will be called the ``$\pmb{
\frown}$-lemma'' and the ``$\pmb{\sqcap}$-lemma.''

Colorings of a polyhedron by a simplicial set were discussed in
Sec.\ref{color}~4.2. A \bem{coloring of a pair of polyhedra $(P,Q)$ by a
pair of
 simplicial sets $(X,Y)$}
is a locally ordered triangulation $p$ of $P$  inducing a
triangulation $q$ of $Q$ together with an $X$-coloring of $p$
inducing an $Y$-coloring of $q$.

\subsubsection{} \label{frown}
Let $D^{k+1}$ be a $(k+1)$-ball and  $S^k=\partial D^{k+1}$ be its
boundary sphere. Let $S_e^{k-1}\subset S^{k}$ be an equator
specified in $S^{k}$, and let $D_+^{k}, D_-^k \subset S^{k} $,
$D_+^{k}\cap D_-^k = S_e^{k-1}$, be the two closed hemispheres
determined by this equator. Consider the embedding of pairs
$(D_+^k,S_e^{k-1})\hra{}(D^{k+1},D_-^k)$, which is illustrated in
the following figure:
$$\input{ball_pairs.pstex_t} $$
By the \bem{$\pmb \frown  $-lemma for $(X,Y)$} we mean the following assertion.
\begin{a_lemma}
For any $k$ any $(X,Y)$-coloring of the pair $(D_+^k,S_e^{k-1})$ can
be extended to an $(X,Y)$-coloring of the pair $(D^{k+1},D_-^k)$.
\end{a_lemma}

\begin{prop}\label{44}
The embedding $|Y| \hra{}|X|$ is a homotopy equivalence  if and only if
the $\pmb \frown$-lemma is valid for  $(X,Y)$.
\end{prop}
\begin{proof}
By the Whitehead theorem, a necessary and sufficient condition for the
homotopy equivalence is the triviality of the relative homotopy classes of
the spheres. By the Zeeman theorem on relative simplicial approximation, the
triviality of the relative homotopy classes of the spheres has the form of
the $\pmb \frown$-lemma.
\end{proof}

\subsubsection{}  \label{pilemma}
Let $A\xar{f}B$ be a map of $\mb{iN}$-sets. The map $f$ is a
\bem{Kan fibration-equivalence} if the following sequence of axioms
is satisfied:
\begin{itemize} \label{kf-e}
\item[KF-E(0):] the map $f_0:A_0\xar{}B_0$ is an epimorphism; \\
..............................
\item[KF-E(k):] for any
$b \in (b)_k$, $a_0,{\ldots} ,a_k \in A_{k-1}$ such that $f(a_i)=d_ib$
and $d_i(a_j)=d_{j-1}a_i$ for $i<j$, there exists
 $a \in A_{k}$ such that
$d_i a=a_i$ and $f(a) = b$. \\
..............................
\end{itemize}
Let $(A^0,A^1)$ and $(B^0,B^1)$ be pairs of $\mb iN$-sets. A pair
of maps $$(A^0,A^1)\xar{f^0,f^1}(B^0,B^1)$$ is a Kan
fibration-equivalence of pairs if both maps $f^0$ and $f^1$ are Kan
fibration-equivalences.
\subsubsection{}
Consider the following injective simplicial sets:
\begin{itemize}
\item the $\mb{iN}$-set $\ol{X}$ whose $k$-simplices are locally ordered,
$X$-colored
triangulations of $\Delta^k$ that subdivide the  ball complex $[\Delta^k]$,
\item the $\mb{iN}$-set $\ol{X}_{\times [0,1]}$ whose $k$-simplices are
locally ordered, $X$-colored triangulations of  $\Delta^k \times
\Delta^1$ that subdivide the ball complex $[\Delta^k] \times
[\Delta^1]$.
\end{itemize}

There are two maps $$\ol{X}_{\times [0,1]} \xar{h_0,h_1} \ol{X},$$
corresponding to the two maps
$$\Delta^k \xar{\id\times d^{0},\, \id\times d^{1}} \Delta^k \times \Delta^1.$$
There is a tautological embedding of pairs $(\mb d X, \mb d
Y)\subset (\ol{X}, \ol{ Y}) $. Here a $k$-simplex of $X$ is
regarded as $\Delta^k$ colored by $x$. Consider the pair of
$\mb{iN}$-subsets $(\wt{X}_{\times [0,1]}, \wt{Y}_{\times [0,1]})
\subset (\ol{X}_{\times [0,1]}, \ol{Y}_{\times [0,1]})$ defined as
follows: $\wt x \in (\wt{X}_{\times[0,1]})_k$ if and only if
\begin{itemize}
\item [(i)] $h_0(\wt x) \in \mb d X_k $,
\item [(ii)] $h_0(\wt x) \in \mb d Y_k
 \Rightarrow \wt x \in (\wt{Y}_{\times[0,1]})_k $,
\item [(iii)] $h_1 \wt x \in \ol{Y}_k$.
\end{itemize}
The \bem{$\pmb \sqcap$-lemma for $(X,Y)$} is the following assertion.
\begin{s_lemma}
The projection $$(\wt X_{\times [0,1]},\wt Y_{\times
[0,1]})\xar{h_0,h_0} (\mb d X, \mb d Y)$$ is a Kan
fibration-equivalence.
\end{s_lemma}
\subsubsection{} \label{intscap}

Geometrically, one can imagine the assertion
of the $\pmb \sqcap$-lemma
for  $(X,Y)$ as follows. Consider the
simplicial bucket $\Delta^k
\cup \partial \Delta^k \times [0,1]$
\begin{figure} 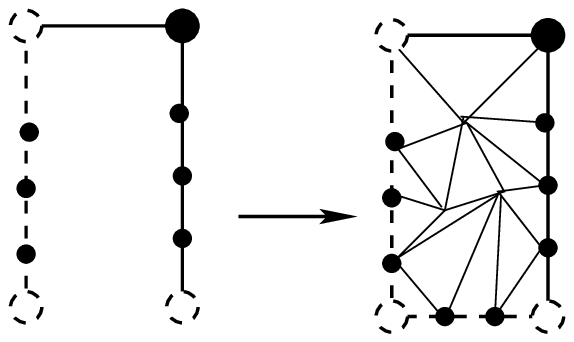 \caption{\label{bin}}\end{figure}
(see Fig.~\ref{bin}, on the left) with a locally ordered triangulation $K$
subdividing the ball complex $[\Delta^k] \cup \partial [\Delta^k] \times
[0,1]$. On the bottom of the bucket, the triangulation $K$ coincides with the standard
ordered triangulation
$[\Delta^k]$. The triangulation $K$ is colored
by $X$ in such a way that the rim $K|_{\Delta^k \times \{1\}}$ of the bucket
is colored by $Y$. If the face $d_i K|_{\Delta^k}$ is colored by $Y$,
then the whole wall $K|_{d_i \Delta^k \times [0,1] }$ is also colored
by $Y$. The $\pmb \sqcap$-lemma states that every such
$(X,Y)$-coloring of the bucket  can be extended to an $(X,Y)$-coloring
of the pair $(\Delta^k \times [0,1], \Delta^k \times \{1\})$.
\begin{prop}\label{45}
The embedding $|X|\hra{}|Y|$ is a homotopy equivalence if and only if the  $\pmb
\sqcap$-lemma
 is valid for
 the pair $(X,Y)$. \end{prop}
\begin{proof}
The $\pmb \sqcap$-lemma follows from the $\pmb \frown$-lemma, because, in
terms of the
 interpretation of Sec.~\ref{intscap}, the $\pmb \frown$-lemma allows us to fill the bucket
 in a required way.

Conversely, the $\pmb \frown$-lemma follows from the $\pmb \sqcap$-lemma.
By induction on
the skeletons of a triangulation $(K,L)$ of $(D^{k+1}, S^k)$, we
construct and fill $(X,Y)$-colored buckets.
\end{proof}

\subsection{Derivation of Lemma~\ref{CombFrag} from Lemma~\ref{CombFrag1}}
\label{lem6} By Propositions~\ref{44} and~\ref{45}, in order to
derive \ref{CombFrag} from \ref{CombFrag1}, it
suffices to prove the following assertion.
\begin{prop}
The $\pmb \frown$-lemma  for the pair $(\wt{\mc W}\mr{Prism}(X), \mc W
\mr{Prism}(X))$ follows from the $\pmb \sqcap$-lemma for all pairs
$$(\underline{\mr{Prism}}^{N-1}(X) \times \PL (X),
 \underline{\mr{Prism}}^N (X))$$
for all finite $N\geq 1$.
\end{prop}

\begin{proof}
Assume that we have a pair $(D^{k+1}, S^{k})$ endowed with a pair of
triangulations $(K,L)$ colored by
$$(\wt{\mc W}\mr{Prism}(X), \mc W \mr{Prism}(X)),$$
i.e., a fixed morphism $$(\bs i K, \bs i L)\xar{(\alpha,\beta)} (\mb d
\wt{\mc W}\mr{Prism}(X),
 \mb d {\mc W}\mr{Prism}(X)).$$
Consider the pair of polyhedra
$$ (D^k \times [0,1], S^{k-1}\times [0,1] \cup D^k \times \{1\})$$
and the pair of ball complexes on these polyhedra
$$(K\times [0,1], L \times [0,1] \cup K \times \{1\}).$$

We will seqrch for an extension of $(K,L)$ to a locally ordered
triangulation
$$(K',L') \trianglelefteq  (K\times [0,1], L \times [0,1] \cup K \times \{1\}),$$
and an extension of $(\alpha,\beta)$ to a coloring
$$(\bs i K', \bs i L')\xar{(\alpha',\beta')} (\mb d \wt{\mc W}\mr{Prism}(X),
 \mb d {\mc W}\mr{Prism}(X)).$$
In order to find an extension of the coloring, we use induction on
the filtration
$$(\Xi_l,\Upsilon_l) = (K \cup (K_{l-1} \times [0,1]), L
\cup (L_{l-1} \times [0,1]) \cup (K_{l-1} \times \{1\})),$$
$$(\Xi_0,\Upsilon_0)=(K,L)\hra{}{\ldots} \hra{}(\Xi_k,\Upsilon_k)=
 (K\times [0,1], L \times [0,1] \cup K \times \{1\}).$$
Let $(\wt{\mc W}\mr{Prism}(X), \mc W \mr{Prism}(X))$ be a coloring
of
 $(\Xi_l,\Upsilon_l)$ extending $(\alpha,\beta)$.
Then the condition for the existence of its extension onto
$(\Xi_{l+1},\Upsilon_{l+1})$ is equivalent to the $\pmb \sqcap$-lemma
for $$(\underline{\mr{Prism}}^{k-l}(X) \times \PL (X),
\underline{\mr{Prism}}^{k-l+1} (X)).$$
\end{proof}

\section{A plan of the proof of the lemma on a common
$\mb R(X)$-triangulation of fiberwise homeomorphisms}

Our project of proving Lemma~\ref{CombFrag1} is to prove the
$\pmb \frown$-lemma (see Sec.~\ref{frown}) for the pair $
(|\underline{\mr{Prism}}^{m-1}(X) \times \PL (X)|,
|\underline{\mr{Prism}}^m (X))|) $. This means that we wish to prove
the following. Let $K$ be a triangulation of the disk $D^k_+$. Let
$\mc Q$ be a coloring of $K$ by the poset $\mb R(X)$. Let
$G_1,{\ldots} ,G_m,U \in \PL_{D^k_+}$ be a family of fiberwise
homeomorphisms such that $G_1,{\ldots} ,G_m$ are $\mc Q$-prismatic and
$U\lfloor_{S_e^{k-1}}$ is prismatic with respect to the coloring $\mc
Q\lfloor_{S_e^{k-1}}$ of the triangulation $K\lfloor_{S_e^{k-1}}$ of
the sphere $S_e^{k-1}=\partial D^k_+$. The claim is that in this
situation $K$ can be extended to a triangulation $\wt K$
of $D^{k+1}$, the coloring $\mc Q$ can be extended to an $\mb
R(X)$-coloring $\wt{\mc Q}$ of $\wt{K}$, the
homeomorphisms $G_1,{\ldots} ,G_m, U$ can be extended to homeomorphisms $\wt
G_1,{\ldots} ,\wt G_m, \wt U \in \PL_{D^{k+1}}(X)$ in such a way that
$\wt G_1,{\ldots} ,\wt G_m$ will remain $\wt{\mc
Q}$-prismatic and $\wt U\lfloor_{D^k_-}$ will become $\wt{\mc
Q}\lfloor_{D^k_-}$-prismatic.

The homotopy
of a family of homeomorphisms to a ``more prismatic''
form is available with the help of a multi-dimensional generalization of
Hudson's construction of
fragmentation of a $\PL$ isotopy. To
describe the generalized Hudson fragmentation and its influence on the
prismaticity,
 we introduce a universal notion of prismaticity, the {\em prismaticity
 of a homeomorphism with respect to an Alexandroff presheaf}.

We start with Sec.~\ref{alex} devoted to Alexandroff topologies on
polyhedra. The set $\mb{Al}(X)$ of all Alexandroff topologies on  a
polyhedron $X$ is partially ordered by strengthening. Any poset is a
$T^0$ Alexandroff space. Therefore, for any Alexandroff topology $\mc
T$ on some other polyhedron $P$, we can consider ``Alexandroff
presheaves'' with values in $\mb{Al}(X)$, i.e., continuous maps of
Alexandroff spaces $(B, \mc T)\xar{} \mb{Al}(X)$. We learn
elementary surgery and homotopy of Alexandroff presheaves.

In Sec.~\ref{gen_prizm} we introduce the notion of  prismaticity of a
homeomorphism $G\in \PL_B(X)$ with respect to an Alexandroff presheaf
$(B, \mc T)\xar{} \mb{Al}(X)$. This notion generalizes the
prismaticity of $G$ with respect to  $\mb R(X)$-colorings of
triangulations of $B$. We study how to deform an Alexandroff
presheaf preserving the prismaticity of a homeomorphism.

Our constructions distinguish
the case of Alexandroff presheaves with
values in $\mb D_c^\infty(K) \subseteq \mb{Al}(X)$. The topologies
from $\mb D_c^\infty(K)$ are discrete everywhere except a disjoint
set of open (in the standard topology) balls that are conic with respect to
some triangulation $K$ of the manifold $X$. The point is that, on
the one hand, homeomorphisms that are prismatic with respect to
presheaves with values in $\mb D_c^\infty(K)$ are the output of
the generalized fragmentation, and, on the other hand, they possess an
$\mb R(X)$-triangulation.

In Sec.~\ref{soglas}, we describe the process of $\mb
R(X)$-triangulating  $\mb D_c^\infty(K)$-prismatic homeomorphisms.

Then we approach the multistage construction of the generalized Hudson
fragmentation of a fiberwise homeomorphism.

In Sec.~\ref{some_scheme}, we introduce a general
nonfiberwise
deformation of a
 fiberwise homeomorphism $G\in \PL_B(X)$ in the class of fiberwise maps.
  The deformation is controlled
by a map $X\times C \xar{F} X\times B$ that is fiberwise
with respect to the projection to $X$. The output of the deformation is a
new map $X\times C \xar{F\rtimes G} X\times C$ that is fiberwise
with respect to the projection to $C$. We prove the important
Proposition~\ref{p52}, which states that if $F\rtimes G$ happens to
be  a fiberwise homeomorphism, then, under some conditions on $F$, the
homeomorphism $F\rtimes G$ ``preserves'' the prismaticity of $G$.

In Sec.~\ref{syst}, we define  ``graph systems'' in the space $X\times
I^k$ (where $I^k$ is a $k$-dimensional cube). With a graph system $\Gamma$
we associate a remarkable map $X\times T^k \xar{F_\Gamma} X \times
I^k$, which is fiberwise with respect to the projection to $X$. The
polyhedron $T^k$ is a $(k+1)$-trapezoid with cubic base $I^k$.

A ``graph system'' $\Gamma$ is a collection of subpolyhedra in
$X\times I^k$ that are the graphs of functions $X\times
I^{k-1}\xar{}I$ embedded along various coordinates. We need some consistency
conditions ({\bf gf1, gf2}) on graphs. The meaning of
these conditions is that all the graphs can be simultaneously blown
up to wide strips without spoiling the overall  picture.
Into $X\times T^k$ the polyhedra are embedded that are the traces
of blowing up the graphs from $\Gamma$ with constant speed. The
map $X\times T^k \xar{F_\Gamma} X \times I^k$ sends a point to its
preimage under the blowing up of $\Gamma$.

In Sec.~\ref{hom_syst}, we study the deformations
$F_\Gamma\rtimes G$. In this sutuation, the general deformations
$\rtimes$ are turned into homotopies between the initial
homeomorphism $G=(F_\Gamma\rtimes G)^0$ and a fiberwise map
$(F_\Gamma\rtimes G)^1$ in the fiber bundle over the upper base
$(T^k)^1$ of the trapezoid $T^k$. Lemma~\ref{l7} states that if
$G^{-1}(\Gamma)$ is also a graph system, then $F_\Gamma\rtimes G$ is
a fiberwise homeomorphism and hence $F_\Gamma\rtimes G$ is a
homotopy of $G$.

Then we observe (see Lemma~\ref{lem8}) that, under an additional condition
{\bf gf4} on  $\Gamma$, the homeomorphism $(F_\Gamma\rtimes G)^1$
becomes prismatic with respect to a remarkable Alexandroff presheaf  $\mb
J_\Gamma$ on $(T^k)^1$.

In Sec.~\ref{hud_gr}, to any triangulation $L$ of a
manifold $X$ and any compact family of homeomorphisms $\mc G \subset
\PL_{I^k}(X)$ we associate a ``Hudson graph system'' in  $\mr H (L, \mc G)$. The
Hudson graph system is consistent
with any homeomorphism  $G\in \mc
G$ and has the property that the image of the Alexandroff presheaf  $\mb
J_{\mr H(L,\mc G)}$ belongs to $\bs S^k(L) \subset \mb{Al}(X)$. The
topologies in $\bs S^k(L)$ are discrete everywhere except the union of
at most
$k$ open stars of $L$. Lemma~\ref{prop48} guarantees that
Hudson graph systems do exist for every $L$ and $\mc G$.

In Sec.~\ref{balls_sep}, we exploit the freedom in the choice of
a triangulation when constructing a Hudson graph system. In Proposition~\ref{fuballsep}
we show that for any triangulation $K$ of a
manifold $X$, we can choose a sufficiently large number $n$ such that
the Alexandroff presheaf  $\mb J_{\mr H(\sd_{n}K,\mc G)}$ can be weakened
to a presheaf $\mb J'$ with  values in $\mb D^\infty_c(K)$. The
weakening preserves the prismaticity of homeomorphisms; therefore for any
$G \in \mc G$ the homeomorphism $(F_{\mr H(\sd_{n}K,\mc G)}\rtimes
G)^1$ is prismatic with respect to the Alexandroff presheaf $\mb J'$ with
values in $\mb D^\infty_c(K)$.

In Sec.~\ref{fragm} we assemble the material accumulated in
Secs.~\ref{some_scheme}--\ref{balls_sep} into Lemma~\ref{lem_frag_on_qube}
``on prismatic fragmentation of fiberwise
homeomorphisms over the cube.'' This lemma guarantees that a finite
family of homeomorphisms in $\PL_{I^k}(X)$ can be deformed to a $\mb
D^\infty_c(K)$-prismatic form by a simultaneous homotopy.

In Sec.~\ref{proof}, Lemmas~\ref{loc_lem} and~\ref{lem13}
perform the final assembly of the $\mb
D^\infty_c(K)$-frag\-menta\-tion of homeomorphisms and the
$\mb R (X)$-triangulation of $\mb D^\infty_c(K)$-prismatic homeomorphisms into
the $\pmb \frown$-lemma for the pair $$
(|\underline{\mr{Prism}}^{m-1}(X) \times \PL (X)|,
|\underline{\mr{Prism}}^m (X))|). $$

\section{Alexandroff spaces (or ``preordered sets'')} \label{alex}
In this section, we present some useful constructions involving
Alexandroff spaces and maps of Alexandroff spaces. We begin by
recalling the definition and standard properties of Alexandroff
spaces \cite{PA} (see also \cite{MC1, Ar}).
\subsection{}\label{preoder}
An \bem{Alexandroff topology} $\mc T$ on a set $Y$ is a topology in
which every point $y\in Y$ has a minimal open neighborhood $o(y)$.
The pair $( Y,\mc T )$ is called an Alexandroff space. An
Alexandroff topology
on $Y$ is equivalent to the structure of a preorder on $Y$.

A \bem{preorder} on a set $Y$ is a transitive and reflective
relation $\preceq$ on $Y$. If, additionally, $\preceq$ is
antisymmetric, i.e., $(a\preceq b)\wedge(b\preceq a) \Rightarrow a=b$,
then it is just a partial order.
Sets endowed with preorders are called \bem{preordered sets}.

An Alexandroff topology $\mc T$ gives rise to the following
preorder on $Y$: $y_0 \us{\mc T}{\leq} y_0 \Leftrightarrow o(y_0)\subseteq o(y_1)$.
Conversely, the lower
ideals of a
preorder on $Y$ form an Alexandroff topology on $Y$. This
correspondence is an isomorphism between the category  all
Alexandroff spaces and continuous maps and the category of
preordered sets and monotone maps. We identify both categories by
this isomorphism and regard an Alexandroff space $(
Y, \mc T )$ also as a preordered set. The category of all Alexandroff
spaces and continuous maps is denoted by $\mb{Al}$. Posets are
the same  as $T^0$ Alexandroff spaces. Thus the category
$\mb{Posets}$ is a full subcategory in $\mb{Al}$.
\subsubsection{ } \label{72}
An Alexandroff topology $\mc T$ on $Y$ has a unique \bem{minimal
base}, which is formed by the minimal neighborhoods of points. We denote the
minimal base of
 $\mc T$ by $\mf{B}( Y, \mc T)$ and regard it as a subposet of the poset  $2^Y$.

In terms of the preorder $\us{\mc T}{\leq}$ on $Y$, the minimal
base is formed by all principal ideals.

\subsubsection{}\label{thA}
The Alexandroff theorem describes the minimal base as a cover.
\begin{thm*}[P.~S.~Alexandroff \cite{Al}]
A cover $\mc U$ of a set $Y$ is the minimal base of an Alexandroff
topology if and only if
\begin{itemize}
\item[{\rm(i)}]  for every pair $U_0,U_1 \in \mc U$, the equality
$U_0 \cap U_1 = \cup_{W\in \mc
W} W$ holds
for some $\mc W \subseteq \mc U$;
\item[{\rm(ii)}] if for $\mc W \subseteq \mc U$ and $U \in \mc U$ the equality
$\cup_{W \in \mc W} W = U$ holds,
 then $U \in \mc W$.
\end{itemize}
\end{thm*}
\noindent The first condition of the Alexandroff theorem says that
$\mc U$ is a base of some topology $\mc T$; the second one guarantees
that the sets $U \in \mc U$  are indeed the minimal neighborhoods in
$\mc T$.

\subsubsection{} \label{sss514}
The correspondence $A \mapsto \mf{B}(A)$ on the objects of  $\mb
{Al}$ has a natural extension to a functor $\mb{Al} \xar{\mf
B} \mb{Posets}$. The map $y\mapsto o(y)$ is a morphism of preordered
sets $(Y,\mc T)\xar{o}\mf{B}(Y,T)$. In terms of $\mb{Al}$, the map
$o$ is a universal map from an Alexandroff space to a $T^0$
Alexandroff space, in the sense that for any map $A\xar{b}P$ in
$\mb{Al}$ such that $P$ is a $T^0$-space, there exists a unique map
$b'$ such that the following
 triangle is commutative:
\begin{equation}\label{zzt}
\bfig \Vtriangle/>`>`<--/[A`P`\mf{B}(A).;b`o`b']\efig
\end{equation}
Thus the maps  $A\xar{o}\mf{B}(A)$ form a natural transformation of
the identity functor into $\mf B$.
\subsubsection{} \label{ss514}
The
intersection of any family of open sets in an Alexandroff
topology is open. Therefore the set $\mc T^c$ of all closed subsets
in $(Y,\mc T)$ is also an Alexandroff topology. The minimal
neighborhoods in $\mc T^c$ are the minimal closed neighborhoods of
points in $(Y,\mc T)$. We denote by  $c(y)$ the minimal closed
neighborhood of $y\in Y$. The minimal base of $\mc T^c$ is formed by
all upper principal ideals of the preordered set $(Y, \mc T)$: $(Y, \mc
T^c)=(Y,\mc T)^\op$. Thus passing to the dual topology is an involution
in $\mb {Al}$. The equivalence classes $\us{\mc T}\sim$  of the
preorder $\us{\mc T}{\leq}$ are the sets $s(y)= c(y)\cap o(y)$. They
form a partition $\Sigma (Y,\mc T)$ of $Y$. There is a one-to-one
map $\mf{B}(Y, \mc T)\xar{\mr{ns}} \Sigma(Y,\mc T)$, which associates to a
principal ideal the equivalence class of its maximal elements: for every
$y\in Y$,
$$\mr{ns}(o(y))=s(y).$$
\subsubsection{\bem{Weakening of Alexandroff topologies}}
Let $\mc R$ and $\mc T$ be two Alexandroff topologies on $Y$. The
identity map of $Y$ induces an $\mb{Al}$-morphism (continuous map)
$\mr (Y, \mc R)\xar{w}\mr (Y, \mc T)$ if and only if $\mc T$ is weaker than
$\mc R$. Therefore we obtain a partial order on  $Y$. We denote this
partial order by the same symbol $\trianglelefteq$ as
the order by subdivision.
We regard the poset of all topologies on $Y$ ordered
by weakening as a subcategory
$\mb{Al}(X)$ of $\mb{Al}$. Simultaneously, $\mb{Al}(X)$ is a poset,
so that it is represented by some object in
$\mb{Posets}\hra{}\mb{Al}$. The poset $\mb{Al}(X)$ has a  maximal
element, the \bem{trivial topology} $X^\mr{triv}$, which has the unique open set $X$.
The poset $\mb{Al}(X)$ has a minimal element, the
\bem{discrete topology} $X^\delta $, in which all the points of $X$ are
simultaneously open and closed.

\subsection{Example} \label{exq}
Let $\ms Q \in \mb R(X)$ be a ball complex on $X$. To $\ms Q$ we can associate
the Alexandroff topology $\mc A (\ms Q)$ on $X$ whose minimal base consists
of the open
(in the standard topology) stars of balls. In the dual
topology $\mc A^c (\ms Q)$, the minimal base consists of the closed
(in the standard topology) balls  of $\ms Q$. The elements of $\Sigma(Y,
\mc A^c (\ms Q))$ are the relative interiors (in the standard topology) of
the balls from $\ms Q$. Consider two ball complexes $\ms Q_0, \ms Q_1
\in \mb R(X)$.

\begin{prop}
$(\ms Q_0 \us{\mb R(X)}{\leq} \ms Q_1) \Leftrightarrow (\mc A(\ms
Q_0)\us{\mb {Al}(X)}{\leq} \mc A (\ms Q_1)) \Leftrightarrow (\mc
A^c(\ms Q_0)\us{\mb {Al}(X)}{\leq} \mc A^c (\ms Q_1)) $.
\end{prop}
\begin{proof} Follows tautologically from definitions.
\end{proof}
Thus the correspondence $\ms Q \mapsto \mc A^c \ms Q$ is an
embedding $\mb R (X) \hra{\mc A^c} \mb {Al}(X)$. Note that the
following diagram is commutative:
$$\bfig \Vtriangle[\mb R(X)`\mb{Al}(X)`\mb{Posets}.;\mc A^c`\mb P` \mf B] \efig $$
\subsection{\bf Alexandroff topologies and covers} \label{cover}
We will associate an Alexandroff topology $\mc A(\mc U)$ to a finite
cover $\mc U$ of a set $Y$. Define the minimal open neighborhood of
$y \in Y$ in $\mc A(\mc U)$ by the equality
$$o(y)=\bigcap_{\{U\,|\, U\in \mc U,\, y \in  U\}}  U.$$
Then the minimal closed neighborhood looks as
$$c(y)=\bigcap_{\{ U\,|\, U\in \mc U,\, y \in  \ol U\}} \ol U,$$
where $\ol U = Y \setminus U$.
\begin{prop}
The topology $\mc A(\mc U)$ on $Y$ is the weakest topology on $Y$
in which all elements of $\mc U$ are open.
\end{prop}
\begin{proof}
Follows directly from the construction.
\end{proof}

\subsection{\bf Inscribing topologies} \label{ss53}

A cover $\mc W_0$ of $Y$ is {\em inscribed} into a cover $\mc W_1$
(we denote this fact by $W_0 \lessdot W_1$) if for every $W_0 \in \mc W_0$
there exists $W_1 \in \mc W_1$ such that $W_0 \subset W_1$. We say
that an Alexandroff topology $\mc R$ on $Y$ is \bem{inscribed} into
an Alexandroff topology
$\mc T$ if the cover  $\mf B(Y, \mc R)$ is inscribed into the cover
$\mf B(Y, \mc T)$. An Alexandroff  topology $\mc T $ is
\label{dense} \bem{dense} if  $ \mf{B}(X,\mc T)$ is a lower
subsemilattice in $2^Y$ (i.e., $U_0, U_1 \in \mf{B}(X,\mc T)$
implies $U_0\cap U_1 \in \mf{B}(X,\mc T)$). In the case where $\mc T$ is
dense, inscribing can be made functorial.
\begin{prop}\label{pro15}
If an Alexandroff topology  $\mc R$ on  $Y$ is inscribed into an
Alexandroff topology $\mc T$ and $\mc T$ is dense, then there is a
poset morphism $\mf B(Y, \mc R)\xar{\phi} \mf B(Y, \mc T)$ such that
for every $R \in B(Y, \mc R)$ the inclusion $R \subset \phi(R)$ holds.
\end{prop}
\begin{proof}
Set $\phi(R)=\bigcap_{R\subseteq T \in \mf R(Y,\mc T)} T $. The
left-hand side of this equality belongs to $\mf R(Y,\mc T)$ by the
definition of a dense topology, and our definition is obviously functorial.
\end{proof}
\subsection{Limits in $\mb{Al}$}
Here we will observe the existence of some limits in $\mb{Al}$. The
following propositions are obvious after switching to the language
of preordered sets.
\begin{prop}
The category $\mb{Al}$ contains
Cartesian squares. To fix the notation, we will
formulate this fact explicitly: the diagram
$$\bfig \square/-->`-->`>`>/[A\us{i,j}{\times}B`A`B`C;j^*`i^*`i`j] \efig $$
 formed by morphisms $i,j$
can be extended to a Cartesian square by morphisms $i^*,j^*$.
\end{prop}
\begin{prop} \label{p15}
The category $\mb{Al}$ contains some co-Cartesian squares, namely,
``pastings,'' i.e., the diagram
$$\bfig \square/<--`<--`<-^{)}`<-^{)}/[A\us{i,j}{\sqcup}B`A`B`C;j^*`i^*`i`j] \efig $$
formed by embeddings $i,j$ can be extended to a co-Cartesian square.
\end{prop}
Let us mention a special case of Cartesian square. Assume that we have
two Alexandroff topologies $\mc R,\mc T$  on a set $Y$. Consider
the morphisms
\begin{equation}\label{eq255} (Y,\mc R)\xar{} (Y,Y^\mr{triv}) \xleftarrow{} (Y,\mc T)\end{equation}
of weakening the topologies $\mc R, \mc T$ to the trivial topology
$Y^\mr{triv}$. The co-Cartesian square generated by the diagram~(\ref{eq255})
is the following square in $\mb{Al}$:
$$\bfig \Square/-->`-->`>`>/[(Y, \mc R\vee \mc T)`(Y,\mc R)`(Y,\mc T)`(Y,Y^\mr{triv}),;```] \efig $$
where $\mc R \vee \mc T$ is the minimal common strengthening of $\mc
R$ and $\mc T$. We can describe the minimal base for $\mc R \vee \mc
T$.
\begin{prop} \label{p16}
For every  $V\in \mf{B}(Y,\mc R \vee \mc T)$ there is a unique
pair $(R,T)$, $R\in \mf{B}(Y,\mc R)$,  $T\in \mf{B}(Y, \mc T)$, such
that $Y=R\cap T$ and $\mr{ns}(R)\cap \mr{ns}(T) \neq \emptyset$.
Conversely, for every pair
 $(R,T)$, $R\in \mf{B}(Y,\mc R)$,  $T\in \mf{B}(Y, \mc T)$, such that
  $\mr{ns}(R)\cap \mr{ns}(T) \neq \emptyset$ there is a unique
element $R\cap T \in \mf{B}(Y,\mc R \vee \mc T)$.
\end{prop}
\subsubsection{}
Let  $A \xar{\phi} B$ be a map of Alexandroff spaces. Consider all
possible extensions of $\phi$ to diagrams of the
form
\begin{equation}\label{cylup1} \bfig \Atriangle/<-`<-`>/[C`A`B;\alpha`\beta`\phi]\efig
\end{equation}
such that
\begin{equation}\label{cylup2}\alpha \geq \beta\circ\phi.
\end{equation}
\begin{prop}\label{pr15}
There exists a diagram
\begin{equation}\label{dia28} \bfig \Atriangle/<-`<-`>/[\overrightarrow{\mr{cyl}}_\phi`A`B;i_0`i_1`\phi]\efig
\end{equation}
such that $i_0 \geq i_1 \circ \phi$ and for every diagram of the
form~{\rm \ref{cylup1}} such that the condition~{\rm(\ref{cylup2})} is satisfied,
there is a unique map
 $\overrightarrow{\mr{cyl}}_\phi\xar{\overrightarrow{\mr{cyl}}_\phi(\alpha,\beta)} C$
such that $$u \circ i_0 = \alpha, \quad u \circ i_1 = \beta.$$
\end{prop}
\begin{proof}
The construction $\overrightarrow{\mr{cyl}}_\phi$ coincides with the
construction of ``homotopy colimit.'' Let $\us{A}{\leq}$ be the
preorder on  $A$ and $\us{B}{\leq}$ be the preorder on $B$. Then we
define $\overrightarrow{\mr{cyl}}_\phi$ as the disjoint union  $\{(a,0)\}_{a\in A}\cup\{(b,1)\}_{b \in B} $ of  $A$
and $B$ endowed with the
preorder
$$(x,i)\us{\overrightarrow{\mr{cyl}}_\phi}\geq (y,j)
\Leftrightarrow \begin{cases} x\us{A}{\geq}y & \text{ for } i=j=0, \text{ or }\\
  x\us{B}{\geq}y & \text{ for } i=j=1, \text{ or } \\
  \phi(x) \us{B}{\geq} y & \text{ for } i =0,\; j =1.   \end{cases}$$
Embeddings $A \xar{i_0} \overrightarrow{\mr{cyl}}(\phi)$ and
 $B \xar{i_1} \overrightarrow{\mr{cyl}}(\phi)$
are defined as $i_0(a)=(a,0)$ and $i_1(b)=(b,1)$. By
construction, $i_0
(a)\us{\overrightarrow{\mr{cyl}}(\phi)}\geq i_1 \circ \phi(a)$, and
the universality property can be checked tautologically.
\end{proof}
Given a morphism $A \xar{\phi}B$, again consider all diagrams
of the form~(\vref{cylup1}), but such that inequality~(\vref{cylup1}) is reversed, i.e.,
\begin{equation}\label{cyld2}
\alpha \leq \beta\circ\phi. \end{equation}
\begin{prop}\label{pr16}
There is a diagram
$$ \bfig \Atriangle/<-`<-`>/[\underrightarrow{\mr{cyl}}_\phi`A`B;i_0`i_1`\phi]\efig
$$
such that $i_0 \leq i_1 \circ \phi$ and for every diagram of the
form~{\rm(\ref{cylup1})} such that {\rm(\ref{cyld2})} is satisfied, there is a map
 $\underrightarrow{\mr{cyl}}_\phi\xar{\underrightarrow{\mr{cyl}}_\phi(\alpha,\beta)} C$
such that $u \circ i_0 = \alpha$, $u \circ i_1 = \beta$.
\end{prop}

\begin{proof}
We can use the duality and observe that after passing to the dual spaces
and maps we find ourselves in the conditions of Proposition~\ref{pr15}.
Thus we can put
$\underrightarrow{\mr{cyl}}_\phi=(\overrightarrow{\mr{cyl}}_{\phi^\op})^\op$.
Therefore
 $\underrightarrow{\mr{cyl}}_\phi$ is the disjoint union
$\{(a,0)\}_{a\in A}\cup\{(b,1)\}_{b \in B} $ of $A$ and $B$ with the preorder
$$(x,i)\us{\underrightarrow{\mr{cyl}}_\phi}\leq (y,j) \Leftrightarrow
 \begin{cases} x\us{A}{\leq}y & \text{ for } i=j=0, \text{ or }\\
  x\us{B}{\leq}y & \text{ for } i=j=1, \text{ or } \\
  \phi(x) \us{B}{\leq} y & \text{ for } i =0,\; j =1.   \end{cases}$$
Embeddings $A \xar{i_0} \underrightarrow{\mr{cyl}}_\phi$ and
 $B \xar{i_1} \underrightarrow{\mr{cyl}}_\phi$
are defined as $i_0(a)=(a,0)$ and $i_1(b)=(b,0)$.
\end{proof}

\subsection{Piecewise linear Alexandroff topologies}
Now  let $Y$ be a closed compact $\PL$ polyhedron.

An \bem{$\mr{OPL}$ Alexandroff topology} on the polyhedron $Y$
is an Alexandroff topology $\mc T$ on $Y$ such that all subsets closed in $\mc T$
are closed $\PL$ subpolyhedra in $Y$.

A \bem{$\mr{CPL}$ Alexandroff topology} on the polyhedron $Y$ is
an Alexandroff topology $\mc T$ on $Y$ such that all subsets open in $\mc T$
are closed $\PL$ subpolyhedra in $Y$.

The duality $\mc T \leftrightarrow \mc T^c$ sends an $\mr{OPL}$
topology to a $\mr{CPL}$ topology and {\it vice versa}.

\subsubsection{\bf $\PL$ Alexandroff  topologies associated with triangulations}

Let $K$ be a finite simplicial complex. The cover of  $|K|$ by the
closed simplices of  $K$ is the minimal base of the $\mr{CPL}$ Alexandroff
topology $\mc A^c(K)$ (see Sec.~\vref{exq}). The minimal base of the
dual $\mr{OPL}$ topology
$\mc A (K)$ is the set of all open stars of the
simplices from $K$.

The following observation essentially belongs to the $\PL$ category.
\begin{prop} \label{prop17}
Any finite cover $\mc U$ of a compact polyhedron $Y$ by closed
subpolyhedra has a triangulation, i.e., there is a triangulation $K$
of $Y$ such that $K$ refines $Y$ as a closed cover.
\end{prop}
\begin{proof}
One can triangulate any $U \in \mc U$ and find a common
triangulation of all these triangulated subpolyhedra. \end{proof}

An immediate corollary is as follows.
\begin{prop}\mbox{}\label{p22}
\begin{itemize}
\item[{\rm1.}] For every  $\CPL$ Alexandroff topology $\mc T$ on a compact
polyhedron $Y$ there exists a linear triangulation
$K$ of $Y$ such that  $\mc A^c(K)\trianglelefteq\mc T$.
\item[{\rm2.}] For every  $\OPL$ Alexandroff topology $\mc R$ on a compact
polyhedron $Y$ there exists a linear triangulation $K$ of $Y$ such
that $\mc A(K)\trianglelefteq\mc R$.
\end{itemize}
\end{prop}
\begin{proof}
1. From the minimality, compactness, and consideration of the dual $\OPL$
topology it follows that the minimal base of $\mc T$ is finite. By
Proposition~\ref{prop17}, we can choose a triangulation $K_{\mc T}$
of the minimal base of $\mc T$. Then $\mc A^c(K_{\mc T})\trianglelefteq\mc T$.

2. Switch to the dual $\mr{CPL}$ topology $\mc R^c$. Choose $K_{\mc
R^c}$. Then $\mc A(K_{\mc R^c})\trianglelefteq \mc R$.
\end{proof}

In Sec.~\vref{dense} we have introduced the notion of a dense Alexandroff topology.
Not every Alexandroff topology is dense, but for any
simplicial complex $K$ the topologies $\mc A K$ and $\mc A^c K$ are
dense for obvious reasons. Therefore  Proposition~\ref{p22} has the following
corollary.
\begin{prop} \label{p24}
Any OPL[CPL] Alexandroff topology on a compact polyhedron has a
dense OPL[CPL] strengthening.
\end{prop}

\subsection{\bf Some properties of  $\mr{CPL}$ Alexandroff topologies} \label{7.4}
Here we observe that maps of $\mr{CPL}$ Alexandroff spaces to
$T^0$ Alexandroff spaces can be regarded as generalized colorings of
a polyhedron by the nerve of a poset (see Sec.~\vref{color}).
\subsubsection{}
Let $Y$ be a polyhedron, $P$ be a poset, and $f$ be a coloring of $Y$
by the nerve $\mc NP$  of $P$. Without loss of generality this means
that we have a locally ordered simplicial complex $K$ such that
$|K|=Y$ and a morphism of injective simplicial sets (see
Sec.~\vref{inj}) $\bs i K \xar{f} \mb d \mc N P$. Consider the $\mr{CPL}$
Alexandroff topology $\mc A (K)$. The coloring $f$ gives rise to
a map $\mr{Max}f$ of the space $(Y, \mc A^c (K))$ to the
poset $P$. The map $\mr{Max}f$ is defined by the following
commutative triangle:
\begin{equation} \label{e33} \bfig \Vtriangle/>`>`<-/[\mr (Y, \mc A^c (K))`P` \mf{B}(Y, \mc A (K)),;
 \mr{Max} f `o` \max f]  \efig\end{equation}
 where $K=\mf B(Y, \mc A^c (K))\xar{\max f} P$ is the map if posets that sends
 a simplex $\delta \in K$ to the value of $f$ at
 the maximal (in the local order) vertex of the simplex.
\subsubsection{} \label{742}
Assume that we have a $\mr{CPL}$ Alexandroff topology $\mc T$ on $Y$. Let
$P$ be a poset, and let a map  $(Y,\mc T)\xar{\xi} P$ be fixed. A
\bem{triangulation of the map $\xi$} is  a coloring
$\bs i K \xar{f}\mb d \mc N P$ of $Y=|K|$
by the nerve of $P$ such that the topology  $\mc
A^c(K)$ is stronger than $\mc T $ and the diagram
\begin{equation} \bfig \dtriangle(0,0)/<-`<-`>/<600,400>%
[P`(Y, \mc A^c (K))`(Y,\mc T);\mr{Max} f`\xi`] \efig\end{equation}
is commutative. The lower arrow in the diagram is the morphism of
weakening the topology.
\begin{prop}\label{prop19}
Assume that
\begin{itemize}
\item $(Y,Z)$ is a pair of compact polyhedra,
\item $\mc T$ is a $\mr{CPL}$ Alexandroff topology on $Y$,
\item $\mc T\lfloor_Z$ is the induced topology on $Z$,
\item $P$ is a poset,
\item $\mr (X,\mc T)\xar{\xi} P$ is a map,
\item $\xi_Z=\xi\lfloor_{(Z,\mc T_Z)}$ is the restriction of $\xi$.
\end{itemize}
Let  $f_Z$  be a triangulation of  $\xi_Z$. Then there exists a
triangulation $f$ of $\xi$ extending $f_Z$.
\end{prop}

\begin{proof}
This is a standard  $\PL$-theorem on the extension of a
triangulation from a closed subpolyhedron to the entire polyhedron.
\end{proof}
\subsubsection{}
We will need a slightly stronger assertion than Proposition~\ref{prop19}
in the case $Z=\emptyset$. Consider a finite poset $R$ and the
geometric realization $|\mr{Ord}(R)|$ of its order complex. On
$|\mr{Ord}(R)|$ we have the $\mr{CPL}$ Alexandroff topology
$\mc A^c(\mr{Ord}(R))$. The vertices of simplices of
$\mr{Ord}(R)$ are
canonically ordered and indexed by the elements of $R$. Therefore
we have the tautological coloring $ \mr{Ord}(R) \xar{f(R)} \mb d \mc N
R$ and the canonical map $(|\mr{Ord}(P)|, \mc A^c(\mr{Ord}(R)) )
\xar{\varkappa(R) } R$ sending a point $x \in |\mr{Ord}(P)|$ to the
maximal vertex of the simplex that contains $x$ in the relative
interior. In this situation, $\mr{Max} f(R)=\varkappa(R)$ and the
coloring $f(R)$ can be recovered from $\varkappa(P)$ in a unique way.
That is, if $\mr{Max} g=\varkappa(P)$ for some other coloring of $|\mr{Ord}
R |$, then $g=f(R)$.
\begin{prop} \label{prop18}
Let $\mc T$ be a $\mr{CPL}$ Alexandroff topology on $Y$. Consider a map $(Y,\mc
T)\xar{\xi} P$, where  $P$ is a poset. Then there is a
finite poset
 $O$ such that $\mr{Ord} (O)$ triangulates $Y$, $|\mr{Ord} (O)|=Y$,
the topology $\mc T$ is weaker than
 $\mc A^c(\mr{Ord}(O)) $, and the commutative diagram of solid arrows
$$\bfig \Square/..>`<-`<-`>/[O`P`(Y, \mc  A^c(\mr{Ord}(O)))`(Y,\mc T);`\varkappa(O)`\xi`] \efig$$
can be uniquely completed by a dashed arrow.
Here the lower arrow is the
morphism of weakening the topology.
\end{prop}
\begin{proof}
It is sufficient to find a common triangulation $K$ of all the
elements from $\mf{B}(Y,\mc T)$. Then $O$ is the
flag poset of this
triangulation and its order complex is the
first barycentric
subdivision of $K$.
\end{proof}
\subsubsection{}\label{744}
We need to record a standard fact. Let $O$ be a finite poset and $P$ be
an arbitrary poset. Assume that we have two maps $O\xar{f_0,f_1}P$. Let $f_1
\geq f_0$. Then there is a canonical homotopy between $f_0$ and
$f_1$, a poset map $O \times \bs 1 \xar{F} P$, where $\bs 1 = (0 <
1)$ is the total order on two elements. The map $F$ is defined as
follows: $F(x,0)=f_0(x)$, $F(y,1)=f_1(y)$. We will specify the
properties of this homotopy in the language of $\mr{CPL}$
Alexandroff spaces. Put $$( |\mr{Ord} O|,\mc A^c \mr{Ord} O)
\xar{\xi_0 = f_0 \circ \varkappa (O),\, \xi_1 = f_1 \circ \varkappa
(O)} P $$ and
$$(|\mr{Ord}(O \times \ms 1)|, \mc A^c (\mr{Ord}(O \times \ms 1)) %
\xar{\Xi = F \circ \varkappa ((O \times \ms 1))} P. $$
We have the commutative diagram
\begin{equation} \xymatrix{ & P & \\
   ( |\mr{Ord} O|,\mc A^c \mr{Ord} O) \ar[ur]^-{\xi_0} & \ar[u]^\Xi \ar@{<-}[l]_-{h_0} \ar@{<-}[r]^-{h_1}  (|\mr{Ord}(O \times \ms 1)|, \mc A^c (\mr{Ord}(O \times \ms 1)) & %
   \ar[ul]_-{\xi_1} ( |\mr{Ord} O|,\mc A^c \mr{Ord} O),
}
\end{equation}
where $h_0, h_1$ are the natural embeddings into the $0,1$ faces, which are
also embeddings in the Alexandroff topology. There is a
homeomorphism $ |\mr{Ord} O \times \bs 1| \xar{g} |\mr{Ord} O |
\times ([0,1]=|\bs 1|)$. The homeomorphism $g^{-1} $ sends the
product of simplices
 $|s|\times [0,1]$ to $|s \times \ms 1| \subseteq
|\mr{Ord} O \times \bs 1|$ with the help of the standard geometric
triangulation of the product of ordered simplices. We have the
Alexandroff topology $\mc A^c \mr{Ord} O \times \mc A^c \mr{Ord} \bs
1$ on $|\mr{Ord} O | \times [0,1]=|\bs 1|$. The minimal base of
this $\mr{CPL}$ topology
(see \cite{Ar}) is formed by all products of
elements of the minimal bases of the factors. That is,
\begin{multline}\mf B ( |\mr{Ord} O | \times [0,1] ,\mc A^c \mr{Ord} O \times \mc A^c \mr{Ord} \bs 1)
\\ = \{\delta \times \{0\}\}_{\delta \in \mr{Ord} O} \cup \{\delta
\times [0,1]\}_{\delta \in \mr{Ord} O} \cup \{\delta \times
\{1\}\}_{\delta \in \mr{Ord} O}.
\end{multline}
Therefore the $\PL$ homeomorphism $g$ induces the
weakening morphism
for $\mr{CPL}$ Alexandroff topologies:
\begin{equation}(|\mr{Ord} O| \times [0,1], \mc A^c \mr{Ord}(O \times \bs 1)) \xar{\tilde g }%
( |\mr{Ord} O | \times [0,1] ,\mc A^c \mr{Ord} O \times \mc A^c \mr{Ord} \bs 1).\end{equation}
The following diagrams are commutative for $i=0,1$:
\begin{equation}\bfig
\Vtriangle(0,0)/>`<-`<-/<900,400>[(|\mr{Ord} O| \times {[}0,1{]}, \mc A^c \mr{Ord} (O \times \bs 1) )`%
(|\mr{Ord} O| \times {[}0,1{]}, \mc A^c \mr{Ord} O \times \mc A^c \mr{Ord} \bs 1) )`%
 ( |\mr{Ord} O|,\mc A^c \mr{Ord} O).;\tilde g`h_{i}`h_i]
\efig
\end{equation}
Define a morphism
$$\Psi': \mf B ( |\mr{Ord} O | \times [0,1] ,\mc A^c \mr{Ord} O \times \mc A^c \mr{Ord} \bs 1)
 \xar{} P $$
 by the rule
$$\Psi' (\delta \times [0,1]) =   \Psi' (\delta \times\{1\})=\mr{max} f_1 (\delta),$$
$$\Psi' (\delta \times \{0\}) = \mr{max} f_0 (\delta).$$
This is a well-defined map, and the following diagram is
commutative:
\begin{equation} \label{e39}\bfig \Atriangle(0,0)/<-`<-`<-/<900,400>[P`(|\mr{Ord} O| \times {[}0,1{]}, \mc A^c \mr{Ord} O \times \mc A^c \mr{Ord} \bs 1) )`%
(\mr{Ord} (O \times \bs 1), \mc A^c \mr{Ord} (O \times \bs 1) ).;\Psi`\Xi`\tilde g]
\efig
\end{equation}
\subsection{Approximation of an $\OPL$ Alexandroff topology
by   a $\CPL$ Alexandroff topology} \label{761} Assume that we have an
$\OPL$ Alexandroff topology $\mc T $ on a compact polyhedron $Y$.
We say that a $\CPL$ Alexandroff topology $\mc C$ \bem{strictly
approximates}
$T$ if the following is true:
there exists a map $\mf B (Y, \mc C) \xar{\xi} \mf B (Y, \mc T)$ such
that $C \subseteq \xi(C)$ for every $C \in \mf B (Y, \mc C) $.
\begin{prop}\label{reg_neighb}
For any $\mr{OPL}$ Alexandroff topology $\mc T$ on a compact
polyhedron $Y$ there is a $\mr{CPL}$ Alexandroff topology that
strictly approximates $\mc T$.
\end{prop}

\begin{proof}
This follows from the theory of regular neighborhoods. In some sufficiently
fine common triangulation of the closed compliments to the elements
of the minimal base of $\mc T$ we choose open regular neighborhoods to
all these complements. Switching again to the compliments, we obtain
a required approximating $\mr{CPL}$ topology.
\end{proof}
\section{Generalized prismatic homeomorphisms} \label{gen_prizm}

\subsubsection{}
Let $(B,\mc T)$ be an Alexandroff space. Let
 $A \xar{f} B$ be a map.
Denote by $\mc T\lfloor_f$ the weakest topology on $A$ in which $f$ is
continuous. The topology $\mc T\lfloor_f$ is again an
Alexandroff topology. Let $(B,\mc T)\xar{\xi}(D,\mc U)$ be an
$\mb{Al}$-morphism, and let $A \xar{f} B$ be a map of sets. We denote by
$\xi\lfloor_f$ the map $(A,\mc T\lfloor_f)\xar{}(D, \mc U)$ such
that $\xi\lfloor_f=\xi\circ f$.
\subsubsection{Alexandroff presheaves and generalized
prismatic homeomorphisms}\label{ss522} Let   $X\times B
\xar{\pi_2}B$ be the trivial bundle, and let  $G \in \PL_{B}(X)$ be a
fiberwise homeomorphism. Let $\mc T$ be an Alexandroff topology on
$B$. The poset $\mb{Al}(X)$ is a poset,
therefore it is an
Alexandroff space. Thus we are allowed to consider a continuous map
\begin{equation}\label{preAl}(B,\mc T)\xar{\xi} \mb{Al}(X). \end{equation}
A continuous map of the form~(\ref{preAl}) will be called an
\bem{Alexandroff presheaf} on $B$.

We say that {\em the fiberwise homeomorphism $G$ is} \bem{$\xi$-prismatic}
 {\em if for every $U \in \mf{B}(B,\mc T)$, every pair of points $b_1,b_2 \in U$,
every $V \in \mf{B}(X, \xi'(U))$, the following condition holds:
$G^{-1}\lfloor_{b_1}(V)=
G^{-1}\lfloor_{b_2}(V)$.}
(For the definition of $\xi'$, see~(\vref{zzt}).)

Generally, we do not assume any connection between the topologies of
the polyhedra $B,X$ and the Alexandroff
topologies from the definition of a
prismatic homeomorphism. If $\mc T$ is the discrete topology $B^\delta$
on $B$, then any $G$ is $\xi$-prismatic for any $\xi$. Similarly, any $G$
is $\xi$-prismatic if $\xi$ sends the whole $(B,\mc T)$ to the trivial
topology $X^{\mr{triv}}$. Note also the following fact.
\begin{prop}
In the definition of the $\xi$-prismaticity of $G$, all
occurrences of $\mf B$ can be replaced by $\Sigma$.
\end{prop}
Piecewise linear homeomorphisms of $X$ act on Alexandroff
topologies on $X$. If a homeomorphism $G$ is  $\xi$-prismatic, then
the correspondence $b \mapsto G^{-1}\xi(b) $ is a new Alexandroff
presheaf $ (B,\mc T)\xar{G^{-1}\xi} \mb{Al(X)}$. Therefore if some
other homeomorphism is  $G^{-1}\xi$-prismatic, then the composition
$F\circ G$ is $\xi$-prismatic. We obtain a groupoid $\mb{GPrism}_B(X)$
whose objects are the Alexandroff presheaves on $B$  with values in
$\mb{Al}(X)$. A $\mb{GPrism}_B(X)$-morphism $\phi \xar{} \xi$ is
a $\xi$-prismatic homeomorphism such that $G^{-1}\xi=\phi$. The
groupoid $\mb{GPrism}_B(X)$ is the groupoid of \bem{generalized
prismatic homeomorphisms} of the trivial bundle $X\times
B\xar{\pi_2}B$.

Let $A\xar{f}B$ be a $\PL$ map and $\phi\xar{G}\xi  \in
\Mor\mb{GPrism}_B(X)$. Then
$$f^*(G)= (\phi\lfloor_{f}\xar{G\lfloor_f}\xi\lfloor_f) \in \Mor \mb{GPrism}_A(X).$$
The correspondence $G \mapsto G\lfloor_f$ is a morphism of groupoids
$$ \mb{GPrism}_B(X) \xar{f^*} \mb{GPrism}_A(X).$$
Therefore $\mb{GPrism}(X)$ is a contravariant functor
$$\PL\xar{\mb{GPrizm}}\mb{Groupoids}.$$
There are three maps associated with every Alexandroff
presheaf $(B,\mc T)\xar{\xi}\mb{Al}(X)$:
\begin{equation}(B,\mc T)\xar{\mf B \xi} \mb{Posets},\quad \mf{B}(B,\mc T)\xar{\wt \xi} \mb{Posets}
,\quad \mf{B}(B,\mc T)\xar{\wt \xi'} \mb{Al}.\end{equation} The following tetrahedral diagram
of functors is commutative:
\begin{equation}\bfig \square(0,0)/>`<-`<-`>/<600,400>[\mb{Al}`\mb{Posets}
`(B,\mc T)`\mf{B}(B,\mc T);\mf B`\xi`\wt{\xi}`o]
\morphism(600,0)/@{->}_<>(.15){\wt{\xi}'}/<-600,400>[\mf{B}(B,\mc T)`\mb{Al};]
\morphism(0,0)/@{->}|<>(.45)\hole^<>(.15){\mf B \xi}/<600,400>[(B,\mc T)`\mb{Posets};]
\efig\end{equation}
By the definition of prismaticity, the homeomorphism $\phi\xar{G}\xi$
induces some isomorphism of functors $\wt \phi \xar{\wt G} \wt \xi$.

\subsubsection{Example}\label{ss604}
Let us show in what sense the notion of a generalized prismatic homeomorphism
contains the notion of a $\ms Q$-prismatic homeomorphism (see
Sec.~\vref{ss322}) with a simplex as the base. 
It was mentioned in
Sec.~\vref{exq} that identifying a ball complex $\ms Q$ with the
Alexandroff topology
$ A^c(\ms Q)$ yields an embedding $\mb
R(X)\xar{}\mb{Al}(X)$. Let $Q=\ms Q_0,{\ldots} ,\ms Q_k$ be a coloring of
$[\Delta^k]$ by $\mb R(X)$. Let $(\Delta^k, \mc
A^c([\Delta^k]))\xar{\mr{Max}\, Q} \mb R(X)$ be the map of Alexandroff
spaces defined in~(\ref{e33}) ) on page~\pageref{e33}. Comparing the
definitions, we obtain the following proposition.
\begin{prop} A homeomorphism  $G\in \PL_k$ is $Q$-prismatic
in the sense of Sec.~{\rm \ref{ss322}} if and only if
 $G$ is $\mr{Max}\,Q$-prismatic in the generalized sense.
\end{prop}
Thus the $\PL$ groupoid $\mb{GPrizm}$ contains the simplicial
groupoid $\mb{Prism}$.

We will usually omit the adjective ``generalized.''

\subsection{Operations on prismatic homeomorphisms}

\subsubsection{Pasting prismatic homeomorphisms} \label{ss611}
Let $A\hra{i_0}B_0$, $A\hra{i_1}B_1$ be embeddings of polyhedra.
Fix $\la G_0, \alpha \ra \in \mb{GPrizm}_{B_0}(X)$, $\la G_1, \beta
\ra
 \in \mb{GPrism}_{B_1}(X)$ such that
 $\la G_0, \alpha \ra\lfloor_{i_0}=\la G_1,
  \beta \ra\lfloor_{i_1} \in \mb{GPrism}_A(X)$.
Consider the following co-Cartesian square in $\PL$:
$$\bfig \Square/<--`<--`<-^{)}`<-_{)}/[B_0\sqcup_{i_0,i_1}B_1`B_0`B_1`A.;i_0'`i_1'`i_0`i_1]\efig$$
\begin{prop} \label{p17}
There is a unique prismatic homeomorphism
$$F=\la G_0 \sqcup_{i_0,i_1} G_1,\alpha
 \sqcup_{i_0,i_1} \beta \ra \in \mb{GPrism}_{B_0\sqcup_{i_0,i_1}B_1}(X)$$
such that $F\lfloor_{i_0'}=\la G_0,\alpha\ra$, $F\lfloor_{i_1'}=\la
G_1,\beta \ra$.
\end{prop}
\begin{proof} By  Proposition~\vref{p15}, the Alexandroff presheaf
$\alpha \sqcup_{i_0,i_1} \beta$ is uniquely defined. The
homeomorphism   $G_0 \sqcup_{i_0,i_1} G_1$ can be defined pointwise.
Its prismaticity with respect to $\alpha \sqcup_{i_0,i_1} \beta$ can
easily be checked.
\end{proof}
\subsubsection{Manipulations with different prismatic structures on a fixed
fiberwise homeomorphism} A homeomorphism $G\in \PL_B(X)$ can be
prismatic with respect to different Alexandroff presheaves. Here we
will describe some useful manipulations with such presheaves.

Let $(B,\mc R)\xar{\alpha} \mb{Al}(X)$ and $(B,\mc T)\xar{\beta}
\mb{Al}(X) $ be Alexandroff presheaves. We write $\alpha \trianglelefteq
\beta$ if $\mc R \trianglelefteq \mc T$ and $\alpha \trianglerighteq \beta\circ w$, where
$(B,\mc R)\xar{w}(B,\mc T)$ is the morphism of weakening the topology.
The definition of prismaticity implies the following proposition.
\begin{prop} \label{p31}
If a homeomorphism  $G$ is $\beta$-prismatic and $\alpha \trianglelefteq
\beta$, then $G$ is $\alpha$-prismatic.
\end{prop}
Therefore all Alexandroff presheaves $\xi$ such that $G$
is $\xi$-prismatic form a poset whose minimal element is the
constant Alexandroff presheaf sending the discrete topology on $B$
to the trivial topology on $X$.

Consider the coordinate embeddings $B\xar{h_0,h_1}B\times [0,1]$ of
a polyhedron $B$ into the ``faces'' of the polyhedron $B\times
[0,1]$, where $h_i=\id\times d^i$.
\begin{prop} \label{p20}
Let $\mc{R,T}$  be two Alexandroff topologies on $B$ and
$\mc R \trianglelefteq\mc T $.
Let $(B,\mc R)\xar{w}(B,\mc T)$ be the morphism of weakening the
topology. Let $\mb C$ be a subposet in $\mb{Al}(X)$. Let $(B,\mc
R)\xar{\alpha} \mb C$, $(B,\mc T)\xar{\beta} \mb C $ be two
Alexandroff presheaves such that $\alpha \trianglerighteq \beta\circ w$. Let
$G\in \PL_B(X)$ be $\beta$-prismatic.
Then
\begin{itemize}
\item[{\rm A.}] The homeomorphism $G$ is $\alpha$-prismatic.
\item[{\rm B.}] There exist an Alexandroff topology
 $\ol{\mr{Cyl}}_w$ on $B\times [0,1]$ and an Alexandroff presheaf
 $(B\times [0,1]\ol{\mr{Cyl}}_w)\xar{\ol{\mr{Cyl}}_w(\alpha,\beta)} C$
such that the homeomorphism $G\times\id\in \PL_{B\times [0,1]}$ is
 $\ol{\mr{Cyl}}_w(\alpha,\beta)$-prismatic,
$h^*_0 (G\times\id, \ol{\mr{Cyl}}_w(\alpha,\beta))=(G,\alpha)$, and
$h^*_1 (G\times\id, \ol{\mr{Cyl}}_w(\alpha,\beta))=(G,\beta)$.
\end{itemize}
\end{prop}
\begin{proof}
Assertion~A coincides with Proposition~\vref{p31}.
Assertion~B must be proved.

\smallskip
B1. We will construct the minimal base of a topology
 $\ol{\mr{Cyl}}_w$ on $B\times [0,1]$.
Define $\mf{B}(B \times [0,1], \ol{\mr{Cyl}}_w)$ as follows:
\begin{equation}\label{eq19}
\{R \times\{0\} \cup \mf{B}(w)(R)\times (0,1]\}_{R \in \mr \mf{B}(B, \mc R)}  \cup
\{T \times (0,1]\}_{T \in \mr N (B, \mc T)}. \end{equation}
Using the Alexandroff theorem (see Sec.~\vref{thA}), it is easy to check that
$\mf{B}(B \times [0,1], \ol{\mr{Cyl}}_w)$ is indeed a minimal base.
According to the construction,
$$(B \times [0,1], \ol{\mr{Cyl}}_w)\lfloor_{h_0}=(B,\mc R),$$
$$(B \times [0,1], \ol{\mr{Cyl}}_w)\lfloor_{h_1}=(B,\mc T).$$

\smallskip
B2. Consider the morphism $\mf{B}(w)$ defined by
$$\bfig \Square/>`<-`<-`>/[\mf{B}(B,\mc R)`\mf{B}(B,\mc T)`(B,\mc R)`(B,\mc T).;\mf{B}(w)`o`o`w] \efig$$

There is a canonical isomorphism
 $\mf{B}(\ol{\mr{Cyl}}_w)\approx \ol{\mr{cyl}}_{\mf{B}(w)}$.
 We will present a construction of this isomorphism.

According to the construction from Proposition~\vref{pr15}, the
Alexandroff space $\ol{\mr{cyl}}_{\mf{B}(w)}$ is defined as the
disjoint union
$$\{(R,0)\}_{R\in \mf{B}(B,\mc R) }\cup\{(T,1)\}_{T \in \mf{B}(B, \mc T)} $$
 of $\mf{B}(B,\mc R)$ and $\mf{B}(B, \mc T)$ with the preorder
$$(U,i)\us{\ol{\mr{cyl}}_{\mf{B}(w)}}\geq (V,j)
\Leftrightarrow \begin{cases} U \supseteq V & \text{ for } i=j=0, \\
  U \supseteq V & \text{ for } i=j=1, \\
  \mf{B}(w)(U) \supseteq V & \text{ for } i =0,\; j =1.   \end{cases}$$

Consider the correspondence defined by the rules
 $(R,0) \mapsto R\times\{0\}\cup \mf{B}(w)(R)\times (0,1]$ for
  $R \in \mf{B}(B,\mc R)$, and
$(T,1)\mapsto T \times (0,1]$ for $T \in \mf{B}(B, \mc T)$.
One can check that this correspondence yields a required isomorphism, the
diagram
$$\xymatrix{ \mr{N}(B,\mc R) \ar[r]^-{\mf{B}(h_0)} &
\ol{\mr{cyl}}_{\mf{B}(w)}& \ar[l]_-{\mf{B}(h_1)} \mf{B}(B,\mc T) \\
(B,\mc R) \ar[r]^-{h_0} \ar[u]^{o}& (B\times[0,1],\ol{\mr{Cyl}}_w) \ar[u]^{o} & \ar[l]_-{h_1} (B,\mc T) \ar[u]^{o}
}$$
is commutative, and the diagram
$$ \bfig \Atriangle/<-`<-`>/[\ol{\mr{cyl}}_{\mf{B}(w)}`\mf{B}(B,\mc R)`\mf{B}(B,\mc T);\mf{B}(h_0)`\mf{B}(h_1)`\mf{B}(w)]\efig
$$
is exactly the universal diagram~(\vref{dia28}).

\smallskip
B3. By the universality of the ``$o$'' morphisms~(\vref{zzt}), we can
replace $\alpha$ and $\beta$ by  unique poset morphisms
$$\mf{B}(B,\mc R)\xar{\alpha'} \mb C, \quad \mf{B}(B, \mc T)\xar{\beta'} \mb C$$
such that $\alpha=\alpha'\circ o$, $\beta=\beta'\circ o$. We can
verify that $\alpha' \geq \beta' \circ \mf{B}(w)$. Next, applying
Proposition~\vref{pr15},
we construct a morphism
$\ol{\mr{cyl}}_{\mf{B}(w)}\xar{\ol{\mr{cyl}}_{\mf{B}(w)}(\alpha',\beta')}\mb
C$ and define an Alexandroff presheaf by
 $$\ol{\mr{Cyl}}_{w}(\alpha,\beta)=\ol{\mr{cyl}}_{\mf{B}(w)}(\alpha',\beta') \circ o.$$

\smallskip
B4. We must verify the
$\ol{\mr{Cyl}}_{w}(\alpha,\beta)$-prismaticity of the homeomorphism
$G\times \id$. It follows from the description~(36) of
elements of the minimal base of  $\ol{\mr{Cyl}}_{w}$.
Proposition~\vref{pr15} and the universality of  $o$-maps
guarantee that with this definition of
$\ol{\mr{Cyl}}_w(\alpha,\beta)$ we have
$$h^*_0 (G\times\id, \ol{\mr{Cyl}}_w(\alpha,\beta))=(G,\alpha)\quad\text{and}\quad h^*_1 (G\times\id, \ol{\mr{Cyl}}_w(\alpha,\beta))=(G,\beta).$$
\end{proof}

The following proposition can be proved similarly, with
the $\ol{\mr{cyl}}$-construction replaced by the
$\underline{\mr{cyl}}$-construction.
\begin{prop} \label{p21}
Let $\mc{R,T}$ be two Alexandroff topologies on $B$ and $\mc R \trianglelefteq
\mc T $. Let
 $(B,\mc R)\xar{w}(B,\mc T)$ be the morphism of weakening the topology.
 Let $\mb C$ be a subposet of $\mb{Al}(X)$.
Let $(B,\mc R)\xar{\alpha} \mb C$, $(B,\mc T)\xar{\beta} \mb C $ be
two Alexandroff presheaves and $\alpha \trianglelefteq \beta\circ w$. Let
a homeomorphism $G\in \PL_B(X)$ be both $\alpha$-prismatic
and $\beta$-prismatic.

Then there exist an Alexandroff topology  $\underline{\mr{Cyl}}_w$
on $B\times [0,1]$ and an Alexandroff presheaf
 $$(B\times [0,1]\underline{\mr{Cyl}}_w)\xar{\underline{\mr{Cyl}}_w(\alpha,\beta)} C$$
such that the homeomorphism  $G\times\id\in \PL_{B\times [0,1]}$ is
 $\underline{\mr{Cyl}}_w(\alpha,\beta)$-prismatic and
$$h^*_0 (G\times\id, \underline{\mr{Cyl}}_w(\alpha,\beta))=(G,\alpha),$$
$$h^*_1 (G\times\id, \underline{\mr{Cyl}}_w(\alpha,\beta))=(G,\beta).$$
\end{prop}
\begin{prop} \label{p32}
Let $\mc{R,T}$ be two Alexandroff topologies on $B$. Let
$\mf{B}(B,\mc R)\xar{\phi}\mf{B}(B,\mc T)$ be a map such that
for every $R \in \mf{B}(B,\mc R)$ the inclusion $R \subset \phi(R)$ holds.
 Let $\mb C \subseteq \mb{Al}(X)$, and let $(B,\mc T)\xar{\beta} \mb C$
 and $(B,\mc R)\xar{\alpha} \mb C$ be two
 Alexandroff presheaves such that
 $\alpha\geq\beta'\circ\phi\circ o$.
Let a homeomoprhism $G\in \PL_B(X)$ be $\beta$-prismatic.
Then
\begin{itemize}
\item[{\rm A.}] The homeomorphism $G$ is $\alpha$-prismatic.\\
\item[{\rm B.}] There exist an Alexandroff topology $\wt{\mr{Cyl}}_\phi$ on
$B\times[0,1]$ and an Alexandroff presheaf
$$(B\times [0,1], \wt{\mr{Cyl}}_\phi) \xar{\wt{\mr{Cyl}}_\phi(\alpha,\beta)} \mb C$$
such that the homeomorphism $G\times\id\in \PL_{B\times [0,1]}$ is
 $\wt{\mr{Cyl}}_\phi(\alpha,\beta)$-prismatic and
$$\la G\times\id, \wt{\mr{Cyl}}_\phi(\alpha,\beta)\ra\lfloor_{h_0}=\la G,\alpha\ra,$$
$$\la G\times\id, \wt{\mr{Cyl}}_\phi(\alpha,\beta)\ra\lfloor_{h_1}=\la G,\beta\ra.$$
\end{itemize}
\end{prop}

\begin{proof}

A. This is obvious.

\smallskip
B1. Consider the common strengthening $\mc R \vee \mc T$ of the
topologies $\mc R$ and $\mc T$ on $B$. Let $$(B, \mc
R)\xleftarrow{w_0}(B, \mc R \vee \mc T )\xar{w_1}(B,\mc T)$$ be the
morphisms of weakening the topology.

\smallskip
B2.1. Consider the Alexandroff presheaves $(B, \mc R \vee \mc T
)\xar{\alpha \circ w_0} \mb C$ and $(B, \mc R \vee \mc T
)\xar{\gamma=\beta \circ w_1} \mb C$. The homeomorphism $G$ is both
$\alpha$-prismatic and
 $\gamma$-prismatic.  We will show that the data set $\la w_0, \gamma, \alpha, G\ra$
satisfies the conditions of Proposition~\vref{p21}.
Namely, we will show that $\gamma \leq \alpha \circ w_0$.

\smallskip
B2.2. By Proposition~\vref{p16}, an element
 $S\in \mf{B}(B, \mc R \vee \mc T )$ is  identified with a unique pair
 $R_S\in \mf{B}(B,\mc R)$, $T_S\in \mf{B}(B,\mc T)$ such that
  $s(R_S)\cap s(T_S) \neq \emptyset$, $R_S\cap T_S=S$, $w_0(S)=R_S$, $w_1(S)=T_S$.
Note that $\phi(R_S) \supseteq T_S$. Indeed, $\phi (R_S)
\supseteq R_S$ and $s(T_S)\cap R_s \neq \emptyset$, whence
 $s(T_S)\cap \phi (R_S) \neq \emptyset$. By the duality between Alexandroff
 topologies and preorders,
this is possible only if
$\phi(R_S) \supseteq T_S$.

\smallskip
B2.3. Thus $\gamma(S)=\beta(T_S)\leq \beta(\phi(R_S))\leq
\alpha(R_S)=\alpha \circ w_0 (S)$. Therefore the data set $\la  w_0,
\gamma,\alpha, G \ra$ satisfies
the conditions of Proposition~\ref{p21}.

\smallskip
B3. By construction, the data set $\la w_1,\gamma, \beta, G \ra $
satisfies the conditions of
Proposition~\ref{p20}.

\smallskip
B4.  We apply Proposition~\vref{p21} and construct an Alexandroff
presheaf
$$(B\times [0,1], \underline{\mr{Cyl}}_{w_0})
\xar{ \underline{\mr{Cyl}}_{w_0}(\gamma,\alpha)} \mb C.$$
We apply Proposition~\vref{p20} and construct an Alexandroff
presheaf
$$(B\times [0,1], \ol{\mr{Cyl}}_{w_1})
\xar{ \ol{\mr{Cyl}}_{w_0}(\gamma,\beta)} \mb C.$$
Consider the following diagram of embeddings in $\PL$:
$$ B\times [0,1]\xleftarrow{h_0}B\xar{h_0}B\times[0,1]. $$
We obtain a condition on prismatic homeomorphisms:
$$\la G \times \id,\underline{\mr{Cyl}}_{w_0}(\gamma,\alpha)\ra
\lfloor_{h_0}=\la G ,\gamma\ra = \la G\times\id,\ol{\mr{Cyl}}_{w_1}(\gamma,\beta)\ra
\lfloor_{h_0}.$$
By Propositions~\vref{p15} and~\vref{p17}, the
pastings are defined.
We put
$$\wt{\mr{Cyl}}_\phi=\underline{\mr{Cyl}}_{w_0}\sqcup_{h_0,h_0} \ol{\mr{Cyl}}_{w_1},$$
$$\wt{\mr{Cyl}}_\phi(\alpha,\beta)=\underline{\mr{Cyl}}_{w_0}(\gamma,\alpha)\sqcup_{h_0,h_0}\ol{\mr{Cyl}}_{w_1}(\gamma,\beta).$$
The pasted homeomorphism $ G\times [0,1]\sqcup_{h_0,h_0} G\times
[0,1]$ is automatically $\wt{\mr{Cyl}}_\phi(\alpha,\beta)$-prismatic.
\end{proof}

\subsection{Triangulation of a fiberwise homeomorphism}
Let $G\in \PL_B(X)$ and $\mb C \subset \mb{Al}(X)$. Let $f$ be a
$\mb C$-coloring of $B$. We call $f$ a \bem{$\mb
C$-triangulation of $G$} if $G$ is $\mr{Max}\,f $-prismatic (see (\vref{e33})).
In particular, according to the remark from
Sec.~\vref{ss604}, for a $\ms Q$-prismatic homeomorphism $G\in\PL_k(X)$, the chain $\ms Q$
is a special case
of an $\mb R(X)$-triangulation of $G$.
 If $|K|=B$ and $K\xar{\mc Q} \mc N\mb R(X)$ is an
$\mb R (X)$-triangulation of $G$, then for every simplex $s\in K$ the
homeomorphism $G\lfloor_s$ is $\mc Q\lfloor_s$-prismatic in the sense
of the definition from Sec.~\ref{ss322}.
\begin{prop}\label{pr64}
Assume that $(Y,Z)$ is a pair of closed polyhedra, $G\in \PL_Y(X)$,
$\mc T$ is a $\mr{CPL}$ Alexandroff topology on $Y$, $\mc
T\lfloor_Z$ is the induced topology
 on $Z$, $\mb C\subset \mb{Al}(X)$, and $\mr (X,\mc T)\xar{\xi} \mb C$ is an
Alexandroff presheaf such that $G$ is $\xi$-prismatic. Let $f_Z$ be
a triangulation of $\xi\lfloor_Z$ that is a triangulation of
$G\lfloor_Z $. Then there is a triangulation $f$ of $\xi$ that
extends $f_Z$ and triangulates $G$.
\end{prop}
\begin{proof}
This is a corollary of Proposition~\vref{prop19}  and
Proposition~\ref{p20}(A),  on page \pageref{p20}.
\end{proof}
\section{Ball complexes compatible with disjoint configurations of balls
} \label{soglas}
\subsection{}\label{mcd}

We need to introduce some general definitions and notations. Let $A$
be a subset of a set $S$.  With $A$ we associate the Alexandroff
topology $\mc D_A$ on $S$ with the minimal base
$\{A\}\cup\{s\}_{s\in S\setminus A}$. Obviously, if $A,B \subseteq
S$ and $A\subseteq B$, then $\mc D_A \trianglelefteq \mc D_B$. Therefore
\begin{equation} \label{eq57}2^S \hra{\mc D} \mb{Al}(S)\end{equation} is an embedding of posets.

Denote by $\mb D(X)$ the poset of all balls of full dimension on a
manifold $X$. The order is defined as follows:
$$D_0 \us{\mb D(X)}{<}D_1 \Leftrightarrow D_0 \subset \int D_1.$$
Denote by $\mb D^\infty(X)$ the poset of all finite sets of disjoint
closed full-dimensional balls  on $X$. The order is defined as
follows:
$$A_0 \us{\mb D^\infty(X)}{<}A_1 \Leftrightarrow
\cup_{D\in A_0} D  \subset \cup_{D\in A_1} \int D.$$

\label{dc} Let $K$ be a simplicial manifold. The simplicial
structure generates a locally conic structure on $|K|$ (see
\cite{RS}). We say that a closed ball $B\subset |K|$
is \bem{consistent}
with $K$ if and only if
\begin{itemize}
\item[(i)] it is full-dimensional and conic,
\item[(ii)] it has the following property: if  $B\cap s \neq \emptyset$ for
some  $s \in K$,
then $\int B \cap s \neq \emptyset$.
\end{itemize}
\begin{figure}
\includegraphics[scale=.5]{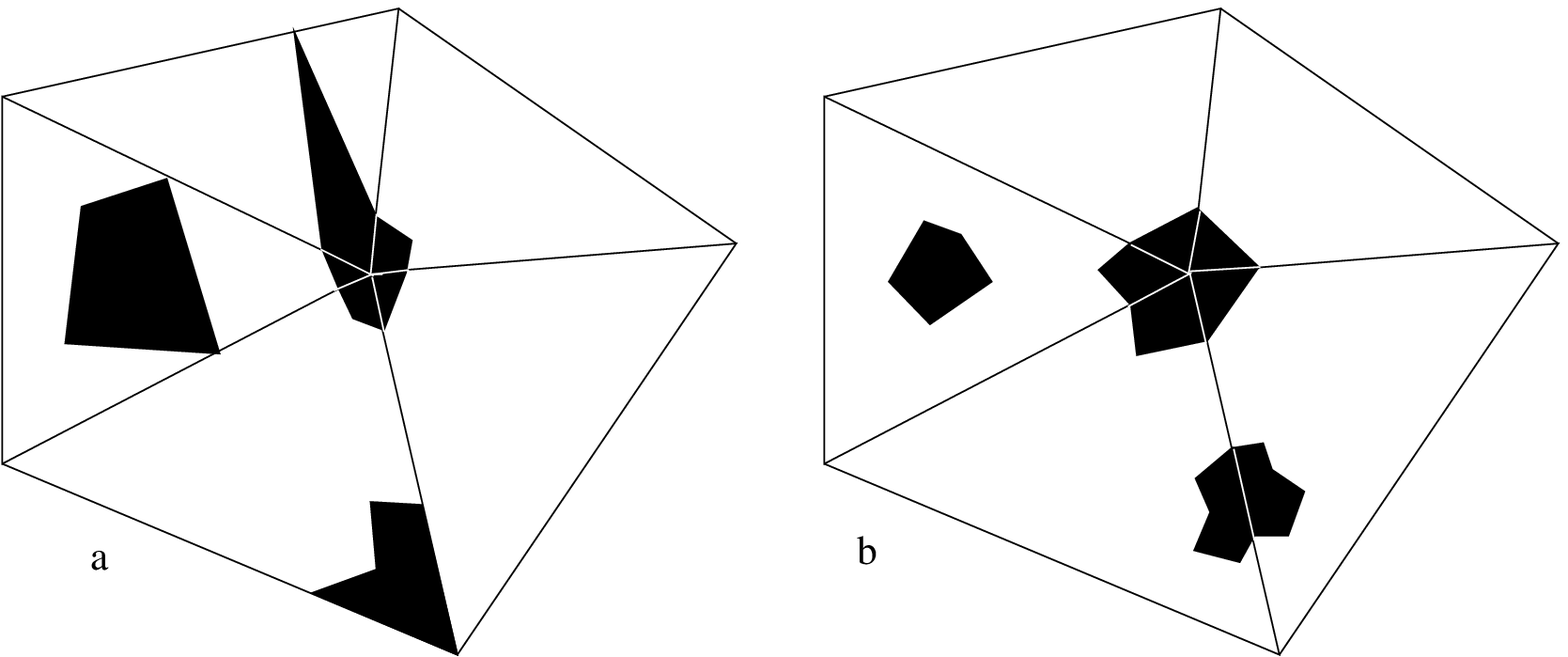}
\caption{\label{bal}}
\end{figure}
The balls shown in Fig.~\ref{bal}\,(a)   are conic with respect
to the triangulation, but condition~(ii) is violated. The balls shown in
Fig.~\ref{bal}\,(b)  satisfy condition~(ii).

\label{dc1} By $\mb D_c^\infty (K)$ we denote the subposet of $\mb
D^\infty(|K|)$ formed by all configurations of
disjoint closed balls consistent with $K$. The functor $\mc D$ identifies the posets $\mb D(X)$,
$\mb D^\infty (X)$, $\mb D^\infty_c (K)$ with subposets of
$\mb{Al}(X)$.
\subsection{}
We will need the following proposition.
\begin{prop} \label{prop55}
Let $P$ be a finite poset, and let $P\xar{\phi} \mb D_c^\infty(K)$
be a poset morphism.
Then there exists a morphism $P\xar{\xi} \mb D_c^\infty(K) $ such
that $\xi(p)>\phi(p)$ for every $p \in P$, i.e., for every $p$, every ball
from $\phi(p)$ belongs to the interior of some ball from $\xi(p)$.
\end{prop}
\begin{proof}
We use the finiteness of  $P$.  By induction on decreasing ranks,
 we can
replace configurations of conic balls with configurations of small
conic neighborhoods of these balls.
\end{proof}

\subsection{}
As usual, let $X$ be an $n$-dimensional compact manifold. Let $\ms Q
\in \mb R(X)$ be a ball complex on $X$. Let $K$ be a triangulation
of $X$, $K\trianglelefteq  \ms Q$, and $\ms D \in \mb D_c^\infty (K)$. We say
that the configuration of balls \bem{$\ms D$ is inscribed into $\ms
Q$} if the interior of any  ball from $\ms D$ is contained in the
interior of an $n$-ball from $\ms Q$. We say that a subdivision  $\ms
Q_0 \trianglelefteq \ms Q$ \bem{does not touch} $\ms D$ if for every ball $D \in
\ms D$ the equality $\ms Q_0|_{\int D}=\ms Q|_{\int D}$ holds.

\begin{prop}\label{p53}
Let $\ms Q \in \mb R(X)$ be a ball complex and
 $K \trianglelefteq \sd_1 \ms Q$. Let $B$  be a
ball consistent with  $K$ with conic vertex $b$. Then there exists a structure of a ball
complex $\ms C(B)$ on $B$ such that the set of all balls of $\ms C(B)$ is the
union of two disjoint sets
$I\ms C(B)$ and $S\ms C(B)$ and the following properties are satisfied:
\begin{itemize}
\item the balls from $S\ms C(B)$ form a ball complex on
 $\partial B$;
\item the balls from $I\ms C(B)$ are the balls touching
$\int B$; and
\item for every ball $ D\in I\ms C(B)$ there exists a ball $D' \in \ms
Q$ such that $D\cap \int B=D'\cap \int B$.
\end{itemize}
\end{prop}
\begin{proof}
1. By the consistency condition, for every ball $Q\in \ms Q$ the intersection
$Q\cap B$ is empty if $b\not \in Q$.

\smallskip
2.  Consider the set  $\wt \st_\ms{Q}b$ of all balls of $\ms Q$
containing the point $b$. In the set $\wt \st_\ms{Q}b$ there is a
unique ball $A$ such that $b \in \relint A$. By the theorem on
invariance of  $\PL$ stars, the following is true: for every $D \in \wt
\st_\ms{Q}b$ the intersection $D\cap B$ is a ball that is a cone over $S(D)=\partial
B \cap D$. Here $S(D)$ is a ball when $D\neq A$, and  a
sphere when $D=A$. For every $Q \in \wt \st_\ms{Q}b$ put $I(Q)=B\cap
Q$, $S(Q)=\partial B\cap Q $.

\smallskip
3. Put $I\ms W (B)=\{I(Q)\}_{Q \in \wt \st_\ms{Q}b}$. Put $S \ms W
(B)=\{S(Q)\}_{Q \in \wt \st_\ms{Q}b}$. Consider the collection of
closed subsets $\ms W(B)= I\ms W(B) \cup S\ms W (B)$.
It has the following properties. The collection $\ms W(B)$ is a cover of $B$
by closed subsets. The interiors of the elements of $\ms W(B)$
form a partition of $B$. The boundary of every element of $\ms
W(B)$ is formed by elements of the cover. All the elements except
$S(A)$ are balls. If we remove  $S(A)$ from $\ms W(B)$ and replace
it by the set $\ms T S(A)$ of simplices of any triangulation of
$S(A)$, then we will obtain a ball complex $\ms C(B)$ on  $B$. By
construction, the decomposition $I\ms C(B)=I\ms W(B)$, $S\ms
C(B)=(S\ms W(B) \setminus S(A))\cup \ms TS(A)$ is what was promised
in the statement of the proposition.
  \end{proof}

\begin{lemma}[Functorially inscribing
disjoint configurations of balls into ball complexes] \label{lem9}
Let $\Lambda$ be a finite poset. Consider a poset morphism
$ \Lambda\xar{\mc Q} \mb R(X)$. Choose a triangulation $K$ of the
common
subdivision of all $\mc Q_\lambda$:
$K \trianglelefteq \mc Q_\lambda$ for every $\lambda \in \Lambda$.
Let  $\Lambda \xar{\xi} \mb D_c^\infty (K)$ be a morphism.
Then there exist a pair of morphisms $\Lambda \to^{\mc Q_0, \mc Q_1}
\mb R(X)$ and a pair of natural transformations
 $\mc Q \toleft^\alpha_\trianglerighteq \mc Q_0 \to^\beta_\trianglelefteq \mc Q_1 $
such that for every $\lambda \in \Lambda $ the subdivisions
$\alpha_\lambda$ do not touch the balls from $\xi(\lambda)$ and the
configuration of balls $\xi(\lambda)$ is inscribed into  $\mc Q_1
(\lambda)$.
\end{lemma}
\begin{proof}
1. First note the following fact.  Let $\ms L$ be a ball
complex on the manifold $X^n$. Let  $\ms B \subseteq \ms L$ be a
cellular embedded
ball complex on an
$n$-ball.  Let us do the
following: delete from $\ms L$ all balls that touch the
interior of $|B|$ and then add the ball $|B|$ itself. This
operation gives a new ball complex on $X^n$. Denote it by $\ms L /
\ms B $. By definition, there is a canonical subdivision $\ms L \trianglelefteq
\ms L / \ms B$. Similarly, if in  $\ms L$ we have a configuration of
 subcomplexes on disjoint  $n$-balls
 $\{\ms B_1,{\ldots} ,\ms B_l\}$, then the  ball complex $\ms L / \{\ms B_1,{\ldots} ,\ms B_l\}$
 is well-defined:
$$\ms L \trianglelefteq  \ms L / \{\ms B_1,{\ldots} ,\ms B_l\} =
 ({\ldots} ((\ms L/\ms B_1)/\ms B_2)/\ms B_3 {\ldots} )/\ms B_l . $$

\smallskip
2. Assume that the poset  $\Lambda$ is of rank $r$. Let
$$\Lambda^r \hra{} \Lambda^{r-1}\hra{}{\ldots} \hra{} \Lambda^0=\Lambda$$
be a filtration of  $\Lambda$ by subposets,  where the subposet
$\Lambda^k$ is formed by all elements of    $\Lambda$ of
rank at least $k$. Denote by $\mc Q^k$ the  restriction $\mc Q
|_{\Lambda^k}$.

We construct natural transformations $\mc Q
\toleft^\alpha_\trianglerighteq \mc Q_0 \to^\beta_\trianglelefteq \mc Q_1 $ inductively on
the inverse rank filtration
$$\xymatrix{
\mc Q^r \ar@{^{(}-}[d] & \ar@{^{(}-}[d] \ar[l]_{\alpha^r}^\trianglerighteq \mc Q^r_0 \ar[r]^{\beta^r}_\trianglelefteq & \ar@{^{(}-}[d] \mc Q^r_1\\
\mc Q^{r-1} \ar@{^{(}-}[d] & \ar@{^{(}-}[d] \ar[l]_{\alpha^{r-1}}^\trianglerighteq \mc Q^{r-1}_0 \ar[r]^{\beta^{r-1}}_\trianglelefteq & \ar@{^{(}-}[d] \mc Q^{r-1}_1\\
    \vdots \ar@{..}[r] \ar@{^{(}-}[d]    & \vdots \ar@{..}[r]             \ar@{^{(}-}[d]                                                 & \ar@{^{(}-}[d] \vdots \\
\mc Q^0 \ar@{=}[d] & \ar@{=}[d] \ar[l]_{\alpha^0}^\trianglerighteq \mc Q^0_0 \ar[r]^{\beta^0}_\trianglelefteq & \ar@{=}[d]\mc Q^0_1 \\
\mc Q  & \ar[l]_{\alpha}^\trianglerighteq \mc Q_0 \ar[r]^{\beta}_\trianglelefteq &\mc Q_1.} $$

\smallskip
3. Let $\lambda \in \Lambda^r$. For  $\mc
Q(\lambda)$ we construct a diagram
$$\mc Q (\lambda) \trianglerighteq \mc Q_0 (\lambda)\trianglelefteq \mc Q_1 (\lambda),$$
which will provide the $\alpha^r(\lambda), \beta^r(\lambda)$-components of
required natural transformations. Pick linear triangulations of
the spherical ball complexes $\partial \ms C_{\mc Q(\lambda)}(D)$, $D
\in \xi(\lambda)$, where $\ms C_{\mc Q(\lambda)}(D)$ are constructed
by Proposition~\ref{p53}. As a result, we obtain some subdivisions of the
ball complexes $\ms T(D) \trianglelefteq \ms C_{\mc Q(\lambda)}(D)$,  $D \in
\xi(\lambda) $, affecting only the boundary spheres. Then we extend the
geometric simplicial complexes $\partial \ms T(D)$, $D \in
\xi(\lambda)$, to a triangulation of $X \setminus \cup_{D \in
\xi(\lambda)} \int D $ subdividing all subpolyhedra $B \setminus
\cup_{D \in \xi(\lambda)} \int D$, $B \in \mc Q(\lambda)$. We obtain a
subdivision $\mc Q \trianglerighteq
\mc Q_0 (\lambda)$ that does not touch the balls of  $ \xi(\lambda)$.
Here $\ms T(D)$, $D \in \xi(\lambda)$, are balls
embedded into $\mc Q_0 (\lambda)$. By Step~1 of the proof, we set
 $$\mc Q_0 (\lambda) \trianglelefteq \mc Q_0 (\lambda)/\{\ms T(D)\}_{D\in \xi(\lambda)} =
 \mc Q_1 (\lambda).$$
By construction, the balls of $\xi(\lambda)$ are inscribed into $\mc Q_1
(\lambda)$. Constructing
 $$\mc Q (\lambda) \trianglerighteq \mc Q_0 (\lambda)\trianglelefteq \mc Q_1 (\lambda)$$
 for all $\lambda \in \Lambda^r$, we
obtain
 $\mc Q^r \toleft^{\alpha^r}_\trianglerighteq \mc Q^r_0 \to^{\beta^r}_\trianglelefteq \mc Q^r_1 $.

\smallskip
4. {\em Inductive step.} Assume that $\mc Q^{k+1}
\toleft^{\alpha^{k+1}}_\trianglerighteq \mc Q^{k+1}_0 \to^{\beta^{k+1}}_\trianglelefteq \mc
Q^{k+1}_1 $ is constructed. Let us extend the construction to
$\Lambda^k$. Let $\lambda \in \Lambda$,
 $\mr{rank}\, \lambda =k $. Consider the upper ideal
$\lambda^\leq = \{\mu \in \Lambda\,|\, \lambda \leq \mu\}$.

Choose some linear triangulations of the spherical simplicial complexes
$\partial \ms C_{\mc Q(\lambda)}(D)$, $D \in \xi(\lambda)$. Thus we obtain
some  subdivisions $\ms T(D) \trianglelefteq \ms C_{\mc Q(\lambda)}(D)$,  $D \in
\xi(\lambda) $, affecting only the boundary spheres. Extend the
geometric simplicial complexes
 $\partial \ms T(D)$, $D \in \xi(\lambda)$,
to a triangulation $X \setminus \cup_{D \in \xi(\lambda)} \int D $
subdividing all subpolyhedra $B \setminus \cup_{D \in \xi(\lambda)}
\int D$, where
 $B \in \mc Q_0(\mu)$, $\mu \in \lambda^\leq $.
We obtain a subdivision $\mc
Q(\lambda) \trianglerighteq \mc Q_0 (\lambda)$
that does not touch the balls of $ \xi(\lambda)$. Here $\ms T(D)$, $D \in
\xi(\lambda)$, are balls embedded into $\mc Q_0 (\lambda)$. By Step~1
of the proof, we set
 $$\mc Q_0 (\lambda) \trianglelefteq \mc Q_0 (\lambda)/\{\ms T(D)\}_{D\in \xi(\lambda)} =
 \mc Q_1 (\lambda).$$
By construction, the balls of  $\xi(\lambda)$ are inscribed into
$\mc Q_1 (\lambda)$. We need to verify that for every $\lambda \leq
\mu $, the following commutative diagram of subdivisions of
ball complexes holds:
$$\xymatrix{
\mc Q(\mu)  & \ar[l]^\trianglerighteq \mc Q_0(\mu) \ar[r]_\trianglelefteq &  \mc Q_1(\mu)\\
\mc Q(\lambda) \ar[u]^\trianglelefteq & \ar[l]^\trianglerighteq  \ar[u]^\trianglelefteq \mc Q_0(\lambda) \ar[r]_\trianglelefteq & \ar[u]^\trianglelefteq \mc Q_1(\lambda).
}$$
It suffices to check that
 $\mc Q_0(\lambda)$ indeed subdivides $\mc Q_0(\mu)$ and
 $\mc Q_1(\lambda)$ indeed subdivides $\mc Q_1(\mu)$.
Here it helps that  $\xi(\mu){\geq}
\xi(\lambda)$ in the sense of the order on $\mb D_c^\infty (K)$. Thus
$\mc Q_0(\mu)|_{\int B}=\mc Q(\mu)|_{\int B}$ for $B \in \xi
(\lambda)$, because $\alpha(\mu)$ does not touch balls from
$\xi(\mu)$, and hence it does not touch balls from $\xi(\lambda)$.
Therefore $\mc Q_0(\lambda)|_B \trianglelefteq \mc Q_0(\mu)|_B $ for $B \in
\xi (\lambda)$. By construction, on the complements of balls from
$\int \xi(\lambda)$ we also have a subdivision. Therefore the whole
$\mc Q_0(\lambda)$ indeed subdivides $\mc Q_0(\mu)$. Similar
arguments show that  $\mc Q_1(\lambda)$ subdivides $\mc Q_1(\mu)$.
\end{proof}
We can develop a slightly more refined version of Proposition~\ref{p53}
for the case of triangulations rather than general ball complexes.
\begin{prop}\label{p54}
Let $\ms Q \in \mb T(X)$, $K \trianglelefteq \sd_1 \ms Q$. Let $B$ be
a closed $K$-consistent $n$-ball with conic vertex $b$. Let
$B_0, B_1 $ be two other closed $K$-consistent balls
with conic vertex $b$ such that $B \subset \int B_{1}$, $B_1
\subset \int B_{0}$.

Under these conditions,  there exist two $\PL$ triangulations $\ms E_0
B_0$, $\ms E_1 B_0$ of  $B_0$
such that
\begin{itemize}
\item[\rm(a)] $\ms E_0 B_0\trianglelefteq \ms C_{\ms Q}B_0$, where $\ms C_{\ms Q}B_0$ is
constructed by Proposition~{\ref{p53}},
\item[\rm(b)] $\ms E_0 B_0\lfloor_{\int B}=\ms Q\lfloor_{\int B}$,
\item[\rm(c)] $\ms E_0 B_0 \trianglelefteq \ms E_1 B_0$,
 $\partial \ms E_0 B_0 = \partial \ms E_1 B_0$,
\item[\rm(d)] $\ms E_1 B_0\lfloor_{B_1}\us{\PL}{\approx} [\Delta^n]$.
\end{itemize}
\end{prop}
\begin {proof}
See Fig.~\ref{triball}.
\begin{figure}
\input{triball.pstex_t}
\caption{\label{triball}}
\end{figure}

1. First, by Proposition~\ref{p53}, we construct a ball complex $\ms
C_{\ms Q}B_0$ and linearly triangulate its boundary. This will be
the linear triangulation $\ms T_0=\partial \ms E_0 B_0 = \partial
\ms E_1 B_0 $. Next we map  $ [\Delta^n]$ onto
$B_1$ piecewise linearly. Denote by  $\ms T_1$ the resulting triangulation of $B^1$
by the
simplex. Construct a $\PL$ triangulation of the annulus
$B^0\setminus \int B^1$ extending $T_0$ and $\partial T_1$ on the
borders. Thus we obtain $E_1 B_0$
satisfying~(d).

\smallskip
2. Now we construct  $\ms E_0 B_0$ satisfying (a), (b), (c). We
need one more $\mc K$-consistent ball $B_2$  with
conic vertex $b$ such that $B \subset \int B_2$, $B_2 \subset\int B_1$.
Using the
pseudoradial projection  with center at $b$, we
can map the closed star of $b$ in $\ms Q$ both onto  $B$ and onto $B_2$. Thus
we obtain two triangulations  $\ms T_2 B_2$ and $\ms T B$ such
that $\ms T\lfloor_{\int B}=\ms Q|_{\int B}$ and $\ms
T_2\lfloor_{\int B_2}=\ms Q\lfloor_{\int B}$. Then we take an arbitrary $\PL$
triangulation of $ B_0 \setminus \int \ms B_2$
that coincides with  $T_0$ on $\partial B_0$ and
subdivides $\ms E_1$ and $D\cap B_0 \setminus \int \ms B_2$ for all
$D\in \ms Q$. We obtain a triangulation $\ms T_{0,2}$ on $B_0 \setminus
\int \ms B_2$. Further, we can extend $T_{0,2}$ to a triangulation
$\ms T_2'$ of
 $B_2$ such that $\partial \ms T_2'=\ms T_{0,2}\lfloor_{\partial B_2}$,
 $\ms T_2'|_{B_0}=\ms T $, and $\ms T_2' \trianglelefteq \ms T^2$.
 Pasting $\ms T_{0,2}$ and $\ms T_2$ by the common part,
 we obtain $\ms E_0 B_0$.
\end{proof}
\begin{lemma}[Functorially inscribing
disjoint configurations
of balls into combinatorial manifolds]
 \label{le10} Let $\Lambda$ be a finite poset.
Consider a poset map $ \Lambda \xar{\mc K} \mb T(X)$. Choose a
triangulation $K$ of the common subdivision of all $\mc K_\lambda$:
$K \trianglelefteq \mc K_\lambda$ for all
$\lambda \in \Lambda$.
Let $\Lambda \xar{\xi} \mb D_c^\infty (K)$ be
a morphism. Then there exist a pair of morphisms $\Lambda
\to^{\mc K_0, \mc K_1} \mb T(X)$ and a pair of natural
transformations
 $\mc K \toleft^\alpha_\trianglerighteq \mc K_0 \to^\beta_\trianglelefteq \mc K_1 $
such that for every $\lambda \in \Lambda $ the subdivisions
$\alpha_\lambda$ do not touch the balls from $\xi(\lambda)$ and the
configuration of balls $\xi(\lambda)$ is inscribed into $\mc K_1
(\lambda)$.
\end{lemma}
\begin{proof}
The proof repeats the proof of Lemma~\ref{lem9}, with
the construction from Proposition~\ref{p53} replaced
by the construction from Proposition~\ref{p54}.
\end{proof}
\subsection{} Now we are able to formulate and prove one of our key lemmas on surgery
of prismatic homeomorphisms.
\begin{lemma} \label{lemma7}
Let $\ms Y$ be a finite simplicial complex, $Y=|\ms Y|$.
Let  $\ms Z\subseteq \ms Y$ be a subcomplex of $\ms Y$, $Z=|\ms Z|$.

Let $G_1,{\ldots} ,G_m, U\in \PL_Y(X)$. Let $\mc Q$ be an
$\mb R(X)$-coloring of $\ms Y$, which induces a coloring
 $\mc Q_{\ms Z}$ of $\ms Z$.

Let $K$ be a common triangulation of all complexes
$\{\ms Q\}_{\ms Q \in  \mr{Max}\mc Q(y),\, y \in Y}$.

Let $\mc T$ be an  $\mr{OPL}$ Alexandroff topology on $Y$. Let
$(Y, \mc T)\xar{\xi} \mb D^\infty_c(K)$ be an Alexandroff presheaf.

Assume that the coloring  $\mc Q$ is an $\mb R(X)$-triangulation of
$G_1,{\ldots} ,G_m$ and the coloring $\mc Q_Z$ is an $\mb
R(X)$-triangulation
of $U\lfloor_Z$.

Assume that the homeomorphisms  $G_1,{\ldots} ,G_m,U$ are also $\xi$-prismatic.

Consider the ball complex $\ms Y \times [I]$ and its subcomplex
$\ms \Pi = \ms Z \times [I] \cup \ms Y \times \{1\}$.

Under these conditions, there exist a triangulation
 $\ms U\trianglelefteq \ms Y \times \ms I$
and an $\mb R(X)$-coloring $\wt{\mc Q}$ of this triangulation
such that
\begin{itemize}
\item$\wt{ \mc Q}|_{\ms Y \times \{0\}} \trianglelefteq \mc Q$,
\item $\wt{ \mc Q}$ triangulates $\wt G_1,{\ldots} ,\wt G_m$,
where $\wt G_i = G_i \times \id \in \PL_{Y \times I}(X)$,
\item $\wt{\mc Q}|_{\ms \Pi}$ triangulates $\wt U\lfloor_{\ms \Pi}$, where
$\wt G_i = G\times \id \in \PL_{Y \times I}(X)$, $\wt U = U\times
\id \in \PL_{Y \times I}(X)$.
\end{itemize}
\end{lemma}
\begin{proof}
1. First we mention a general fact: if a $\mr{CPL}$ Alexandroff
topology $\mc T'$ on $Y$ strictly approximates (see Sec.~\vref{761}) an
$\mr{OPL}$ Alexandroff topology $\mc T$ and $(Y,\mc T')\xar{\xi'}
\mb D_c^\infty(K)$ is the induced morphism, then any $\xi$-prismatic
homeomorphism
 $f\in \PL_Y(X)$
is $\xi'$-prismatic. By Proposition~\vref{reg_neighb}, there always
exists a strict approximation $\mc T'$
for $\mc T$. Let us fix such an approximation.

2. Now pick the common strengthening $\mc A^c (\ms Y) \vee \mc T' $
of the $\CPL$ topologies
$\mc A^c (\ms Y)$ and $\mc T'$.
The canonic weakening morphisms are
$$(Y, \mc A^c(\ms Y))\xleftarrow{\varphi_{\mc A^c (\ms Y)}}
(Y, \mc A^c(\ms Y) \vee \mc T')\xar{\varphi_{\mc T}} (Y, \mc T').$$
Consider the map
$$\mr (Y, \mc A(\ms Y) \vee \mc T')
\xar{\gamma=(\mr{Max}\mc Q \circ
 \varphi_{\mc A^c (\ms Y)})
\times (\xi'\circ\varphi_{\mc T'})} \mb R(X)\times \mb D^\infty_c ( K).$$
In this situation, by Proposition~\vref{prop18}, there is a finite
poset $O$ such that $\mr{Ord} (O)$ linearly triangulates $Y$,
$|\mr{Ord} (O)|=Y$, the topology $\mc A^c (\ms Y) \vee \mc T'$ is
weaker than $\mc A^c(\mr{Ord} (O))$, and the commutative diagram of
solid arrows
$$\bfig \Square/>`>`>`..>/[(Y, \mc A^c(\mr{Ord} (O))`
 (Y, \mc A^c(\ms Y)
 \vee \mc T)`O` \mb R(X)\times \mb D^\infty_c(K);
 `\varkappa(O)`\gamma`P_R \times P_D] \efig$$
can be uniquely completed by a dashed arrow.
(Here the top arrow
$$(Y, \mc A^c(\mr{Ord} (O) )\to (Y, \mc A^c(\ms Y) \vee \mc T)$$
is the morphism of weakening the topology.)

3. \label{comb} Now, by Lemma~\vref{lem9}, we can build a pair of
functors $O\xar{ P_0, P_1} \mb R(X)$ and a pair of natural
transformations
$ P_R \toleft^\alpha_\trianglerighteq  P_0\to^\beta_\trianglelefteq  P_1 $
such that for every
$o \in O$, the subdivision
$\alpha_o$ does not touch the balls from  $P_D (o)$ and the
configuration  $P_D(o)$ is inscribed into $P_1(o)$.

4. Now for $\alpha$ and $\beta$ we can use
the arguments from Sec.~\vref{744}.
By pasting together the  two colorings $[\mr{ord} O \times \bs
1]\xar{\Xi(\alpha)} \mb R(X)$ and $[\mr{ord} O \times \bs
1]\xar{\Xi(\beta)}\mb R(X)$ on $ Y\times [0,\frac{1}{2},1]$, we obtain
the required   $\ms U$ and $\wt{\mc Q}$.
\end{proof}

\section{Nonfiberwise deformations of  fiberwise homeomorphisms}
\label{some_scheme} Our key tool is a set of lemmas on fragmentation
of fiberwise homeomorphisms of the trivial fiber bundle over the
cube.
\subsection{The operation $\rtimes$.}\label{ss821}
Let $H\in \PL_{B}(X)$  be a fiberwise homeomorphism of the trivial
bundle with  fiber $X$ and base $B$: $H(x,b)=(H\lfloor_{b}(x),b)$.
Consider a map $X\times C \xar{F} X\times B$ that is fiberwise
with respect to the projection to $X$: $F(x,c)=(x,F\lfloor_x(c))$.
With the pair
$$\bfig \Vtriangle(0,0)<300,300>[X\times C`X\times B,`X;F`\pi_1`\pi_1] \efig\quad
 \bfig \Vtriangle(0,0)<300,300>[X\times B`X\times B`B;H`\pi_1`\pi_1] \efig  $$
we associate a map
$$ \bfig \Vtriangle(0,0)<300,300>[X\times C`X\times C`C;F\rtimes H`\pi_1`\pi_1] \efig$$
that is fiberwise with respect to the projection to $C$. The map
$F\rtimes H$ is defined by the following correspondence:
\begin{equation}\label{eqt0}(x,c)\overset{F\rtimes H}{\mapsto} (H\lfloor_{(F\lfloor_x(c))}(x),c).\end{equation}
Generally, $F\rtimes H$ is only a fiberwise map, but not a homeomorphism.
\subsection{}\label{ss822} Consider a special case of the construction $\rtimes$.
Consider the following embedding:
$$X\times B \xar{ i^0=\id\times\id\times d^0} X\times B \times I.$$
Let $X\times B\times I \xar{F} X\times B$ be a map that is
fiberwise with respect to the projection to $X$ and such that $F \circ
i^0 = \id$. Let $H\in \PL_B(X)$. Let $h_t$ be the embedding
$B\hra{h_t}B\times I$ defined by $h_t(b)=(b,t)$. The map
 $(F\rtimes H) \lfloor_{h_t}: X\times B \xar{} X\times B$ induced by $h_t$
is fiberwise with respect to the projection to $B$. Generally, this is
not a homeomorphism when $t\neq 0$. But if $t=0$, it follows from the
construction that $ (F\rtimes H) \lfloor_{h_0} \equiv H$. In
this situation, $F\rtimes H$ is a one-parameter deformation of $H$
in the class of fiberwise maps.
\subsection{Prismaticity properties of the operation $\rtimes$}\label{ss823}
Let $(B, B^\mr{triv})\xar{\delta_{\mc L}}\mb{Al}(X)$ be the
Alexandroff presheaf that sends the whole $B$ to the Alexandroff
topology $\mc L$ on $X$. Let $H$ be a homeomorphism that is
prismatic with respect to $\delta_{\mc L}$. By the definition from
Sec.~\vref{ss522}, we have the topology $\mc L'=H^{-1}\mc L$  and
the morphism $(B, B^\mr{triv})\xar{H^{-1}\delta_{\mc L}'}\mb{Al}(X)$
that sends the whole
  $B$ to
$\mc L'$. In this situation, for every
 $ b \in B$ the map
$(X, \mc L)\xar{H^{-1}\lfloor_{b}}(X, \mc L')$ is a homeomorphism of
Alexandroff spaces that induces a constant (independent of $b$)
 isomorphism $\mf{B}(X,\mc L)\xar{\wt H^{-1}} \mf{B}(X,\mc L')$.
\begin{prop}\label{p50}
Let $\mc L $ be an Alexandroff topology on $X$. Let $H$ be prismatic
with respect to the constant Alexandroff presheaf $(B,
B^\mr{triv})\xar{\delta_{\mc L}}\mb{Al}(X)$ that sends the whole
$B$ to $\mc  L$.

Then  $F\rtimes H^{-1}$ has the following property: for every $L \in
\mf{B}(X, \mc L)$, for every $c \in C $, the inclusion $(F\rtimes H^{-1})\lfloor_{c}(L)
\subseteq \wt H^{-1}(L)$ holds.
\end{prop}
\begin{proof}
By the definition~(\vref{eqt0}),
\begin{equation}\label{eqt} (F\rtimes H^{-1})\lfloor_{c}(x)=H^{-1}\lfloor_{(F\lfloor_x(c))}(x).\end{equation}
The right-hand side of~(\vref{eqt}) is contained in $\wt H^{-1}(L)$ for
$x\in L$.
\end{proof}
Therefore, under the conditions of Proposition~\ref{p50}, if
$F\rtimes H^{-1}$ is  a fiberwise homeomorphism, then  $(F\rtimes
H^{-1})^{-1}$ is a prismatic homeomorphism with respect to
$(C,C^\mr{triv})\xar{\delta_{\mc L}}\mb{Al}(X)$.
\begin{prop}\label{p51}
Let $\mc L $ be an Alexandroff topology on $X$. Let  $H$ be
prismatic with respect to
 $(B, B^\mr{triv})\xar{\delta_{\mc L}}\mb{Al}(X)$. Let
$F\rtimes H^{-1}$ be a fiberwise homeomorphism.
Then the inverse fiberwise homeomorphism is prismatic with respect to the
Alexandroff presheaf $(C,C^\mr{triv})\xar{\delta_{\mc L}}\mb{Al}(X)$
and $(F\rtimes H^{-1})\mc L \equiv H^{-1}\mc L$, $\wt{F\rtimes
H}^{-1}\equiv \wt{H}^{-1}$.
\end{prop}
\begin{proof} The assertion of Proposition~\ref{p50} is equivalent
to the same assertion with $\mf{B}(X,\mc L)$ replaced by the partition
$\Sigma(X,\mc L)$ (see Sec.~\vref{sss514}). Under our conditions, the maps
 $F\rtimes H^{-1}\lfloor_c$ are one-to-one, so the elements of the partitions
are mapped ``onto.'' Therefore the elements of the minimal bases are
also mapped ``onto.'' This provides the required assertion for
$(F\rtimes H^{-1})^{-1}$.
\end{proof}

Now we will slightly generalize the previous arguments.
\begin{prop}\label{p52}
Let $\mc T$ be an Alexandroff topology on $C$ and $R$ be an
Alexandroff topology on $B$. Let $\mf{B}(C,\mc
T)\xar{\phi}\mf{B}(B,\mc R)$ be a map of Alexandroff spaces.
Let $(B,\mc R)\xar{\xi}\mb{Al}(X)$ be an Alexandroff presheaf. Let
 $H\in \PL_B(X)$ be a $\xi$-prismatic homeomorphism.
Assume that the map $F$ defined by the commutative triangle
$$\bfig \Vtriangle(0,0)<300,300>[X\times C`X\times B`X;F`\pi_1`\pi_1] \efig$$
has the following ``block'' property: for every $T \in \mf{B}(C,\mc T)$, the
equality $F(X\times T)=X\times \phi (T)$ holds.
Let $F\rtimes H^{-1}$ be a fiberwise homeomorphism.

Then the homeomorphism $(F\rtimes H^{-1})^{-1}$ is
 $(\xi' \circ \phi)$-prismatic.
\end{prop}
\begin{proof}
On each separate block of $X\times T$, we are under the conditions of
Proposition~\ref{p51}. Pasting the conclusion of
Proposition~\ref{p51} over all $T\in \mf{B}(C,\mc T)$ provides the
conclusion of the current proposition.
\end{proof}

\section{Graph systems}\label{syst}
\subsection{}
Let $\mb k$ be a finite set of indices. The cube $I^{\mb k}$ is a
coordinate cube in $\R^\mb k$ with coordinate functions $t_k$,
$k\in \mb k$. For a subset $\mb s \subseteq \mb k$ and $\epsilon
\in \{0,1\}^{\mb k \setminus \mb s}$, denote by $I^{{\mb s},\epsilon}$
the face of $I^{{\mb k}}$ defined by the equality $t_i=\epsilon_i$
for coordinates with indices from $\mb k \setminus \mb s$. The
``coordinate'' faces $I^{\mb{s},\bs 0}$ are denoted by $I^\mb{s}$. Let
$\mb s \subseteq \mb k$ be fixed. Let $\Gamma =\{\Gamma_i\}_{i \in
\mb s} $ be a collection of subpolyhedra in $X\times I^{\mb{k}}$.
Consider the following condition on $\Gamma$:

\medskip\noindent
{\bf gf:} \label{gf} For every subset $\mb l \subseteq \mb s$ there
exists a map $X \times I^{\mb {k}\setminus \mb{l}}\xar{f_\mb{l}}
I^\mb{l}$ such that $\Gamma_{\mb l}=\bigcap_{i\in \mb{l}} \Gamma_i=
\Gamma f_\mb{l}$. Here $\Gamma f_\mb{l}$ is the graph of the map $f_\mb{l}$ in $X\times
I^\mb{k}$.

\medskip\noindent
If the condition {\bf gf} is satisfied, then the maps
 $\{f_{\mb l}\}_{\mb l \subseteq \mb s}$
can be recovered from $\Gamma$ in a unique way. Let $N$ be a positive integer.
Denote by $\bcc{N}$ the set $\{0,{\ldots} ,N\}$.
Define a set $\mb{GF}(X,\mb k, N)$ as follows. An element $\Gamma
\in \mb{GF}(X,\mb k, N)$ is a collection
$\Gamma=\{\Gamma_i^j\}_{i\in \mb k,\, j\in \bcc{N}}$ of subpolyhedra
in $X\times I^\mb{k}$ for which the following conditions {\bf gf1},
{\bf gf2} are satisfied.

\medskip\noindent
{\bf gf1:} For every set $\mb s \subseteq \mb k$ and every map $\mb s
\xar{b} \bcc{N}$, the condition {\bf gf} is satisfied for the collection
 $\{\Gamma_i^{b(i)}\}_{i\in \mb s}$.

\medskip
This means that there exists a map $X\times I^{\mb k \setminus \mb
s}\xar{f^b}I^\mb s$ such that $\Gamma^b=\bigcap_{i\in \mb s}
\Gamma^{b_i}_i = \Gamma f^b$, where
 $\Gamma f^b$ is the graph of the map
$f^b$ in $X\times I^\mb{k}$.

The maps  $f^b$ are uniquely determined by $\Gamma$. If  {\bf gf1} is
satisfied, then
 $\Gamma^j_i=\Gamma f^{i\mapsto j}$
for a unique map $X\times I^{\mb k \setminus
\{i\}}\xar{f^{i\mapsto j}} I^{\{i\}}$. These special maps
$f^{i\mapsto j}$ will be denoted by $f_i^j$. Further, we assume that for
$\{f^j_i\}$
the following condition is satisfied:

\medskip\noindent
{\bf gf2:} $0\equiv f^0_i\leq f^1_i \leq {\ldots} \leq f^{N-1}_i \leq
f^N_i\equiv 1$ for all $ i \in \mb k$.

\subsection{}\label{ss831} We need to use indices that are maps of
finite sets into $\N$. Let us introduce some notations and
agreements.

Let $\mb{r,s}$ be subsets of $\mb k$. Let $\mb s \xar{c} \N$, $\mb r
\xar{d} \N$ be two maps such that $c|_{\mb s \cap \mb r}\equiv
d|_{\mb s \cap \mb r}$. Denote by $c\sqcup d$ the map $\mb s \cup \mb r
\xar{c\sqcup d} \N$ such that $c\sqcup d|_{\mb s}\equiv c$, $c\sqcup
d|_{\mb r}\equiv d$. For a numerical function $\N\xar{p}\N$ and $a\in
\N^\mb{k}$, we denote by $p(a)\in \N^\mb{r}$ the result of the
componentwise application of $p$ to $a$.

\subsection{}\label{ss832}

To every subset $\mb s \subset \mb k $ and a map $\mb s \xar{b}\bcc{N}$ we
can associate a map
$$\mb{GF}(X,\mb k, N)\xar{\delta_b} \mb{GF}(X,\mb k\setminus \mb s, N) $$
as follows. For $i\in \mb k \setminus \mb s $, $j\in \bcc{N}$, we put
\begin{equation}\label{equ49}(\delta^b \Gamma )_i^j=\pi_{x,t_{\mb k \setminus s}}
(\Gamma^{b\sqcup (i\mapsto j)})=\Gamma f_i^{b\sqcup (i\mapsto j)},\end{equation}
where $\pi_{x,t_{\mb k \setminus \mb s}}:X\times I^\mb{k}\xar{}
X\times I^{\mb k \setminus \mb s}$ is the coordinate projection.
\subsection{}
Consider a map $\mb k \xar{a} \bcc{2N-2}$. Define a polyhedron
$\Omega^0_a(\Gamma) \subset X\times I^\mb{k}$ as follows. Put
\begin{multline*} \shoveright{(\Omega_a^0)^o =  \{(x,t_\mb{k})\,|\, \forall i \mbox{ with } a_i=2l_i-1 :
 t_i=f_i^{l_i}(x,t_{\mb k \setminus \{i\}})\} = \bigcap_{a_i=2l_i-1} \Gamma_i^{l_i},}\end{multline*}
\begin{multline*}(\Omega_a^0)^e =  \{(x,t_\mb{k})\,|\, \forall i \mbox{ with } a_i=2l_i :
 f_i^{l_i}(x,t_{\mb k \setminus \{i\}})\leq t_i\leq
f^{l_i+1}_i(x,t_{\mb k \setminus \{i\}})\} \\ = \bigcap_{a_i=2l_i}
 ((\Gamma^{l_i+1}_i)^\leq \cap (\Gamma^{l_i}_i)^\geq). \end{multline*}
(Here $(\Gamma_i^j)^\geq$ and $(\Gamma_i^j)^\leq$ are
the closed supergraph and the subgraph of $\Gamma_i^j$, respectively).
Put
\begin{equation}\label{equ50} \Omega^0_a = (\Omega_a^0)^o\cap (\Omega_a^0)^e \end{equation}
and $$\Omega^0(\Gamma)=\{\Omega^0_a\}_{a\in \bcc{2N-2}}.$$
\subsection{}
Fix the following data: $\mb s \subset \mb k$, $\mb s \xar{b}
\bcc{N}$, $\Gamma \in \mb{GF}(X,\mb k, N)$. According to
Sec.~\ref{ss832},
 $\delta^b \Gamma \in \mb{GF}(X,\mb k\setminus \mb s, N) $ is defined.
Let us study the relations between $\Omega^0(\Gamma)$ and
$\Omega^0(\delta^b \Gamma)$. Let $\mb k \setminus \mb s \xar{c}
\bcc{2N-1}$. According to  (\ref{equ49}), there are canonical
mutually inverse  homeomorphisms
\begin{equation}\label{equ51}\Omega^0_c(\delta^b \Gamma)\two/<-`->/^{\pi_{x,t_{\mb k \setminus \mb s}}}_{f^b}
\Omega^0_{c\sqcup 2b-1}(\Gamma).\end{equation}

\subsection{}
Construct a trapezoid $T^\mb{k}(N) \subset \R^\mb{k}\times
\R^{\{v\}}$, where $\R^{\{v\}}$ is a copy of $\R$ with
coordinate function $v$.
By $e_v$ we denote the corresponding basis vector.

\begin{figure}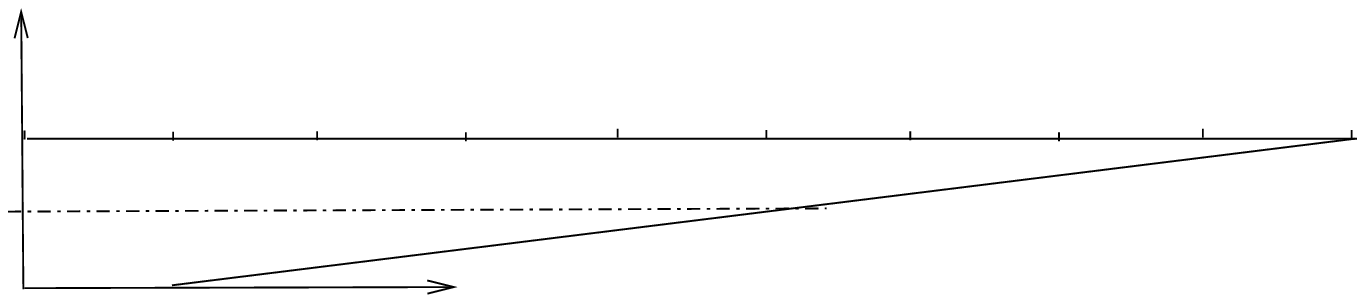\caption{}\label{hudson2}\end{figure}

Denote by $ \br{2N-1} $ the interval  $[0, 2N-1] \subset \R$. Denote
by $T^\mb{k}$ the trapezoid in $\R^\mb{k}\times \R^{\{v\}}$ that is
the convex hull of two embedded cubes: $T^\mb{k} = \conv (I^\mb{k}, e_v
+ \br{2N-1}^{\mb k})$. Figure \ref{hudson2} presents the trapezoid $T^1$.

The trapezoid $T^\mb{k}$ is the product of  $\mb k$
copies of $T^1$ fibered by  $v$. We have the embeddings
$$I^\mb{k} \xar{h_0} T^\mb{k} \xleftarrow{h_1} \br{2N-1}^\mb{k}, $$
where $h_0(t_\mb{k})=(t_\mb{k}, 0_{v})$,
$h_1(t_\mb{k})=(t_\mb{k},1_v)$. We denote the images of $h_0$ and
$h_1$ by $H_0^k$ and $H_1^\mb{k}$, respectively. For each face
$I^{\mb{k}\setminus \mb s, \epsilon}\subset I^\mb{k}$ of the cube,
we have a map $T^\mb{k\setminus \mb
s}\hra{d_{s,\epsilon}}T^\mb{k}$, which is defined as follows:
$$(t_\mb{s},v)\mapsto (t_\mb{s},g^\epsilon_{\mb k \setminus \mb s}(t_\mb{s},v)),$$
where
$$
g^\epsilon_{i}(t_\mb{s},v)=\begin{cases} 1+2(N-1)v & \text{ for } \epsilon_i=1, \\
0 & \text{ for } \epsilon_i=0.
\end{cases}
$$
Denote the image of $d_{s,\epsilon}$ by
 $T^{\mb{s},\epsilon} \subseteq T^\mb{k}$.
Denote the union
$\cup_{i\in \mb k,\epsilon_i=0,1} T^{\mb{k}\setminus\{i\},\epsilon_i}$
of all polyhedra
by $W^\mb{k}$. Denote the faces of $W^\mb{k}$ that
correspond to $v=0$ and $v=1$  by
$h^0(W^\mb{k})$ and $h^1(W^\mb{k})$, respectively.
\begin{prop}
There exists a noncanonical cellular  $\PL$ homeomorphism
$[I^\mb{k}\times I^{v}]\xar{} [T^\mb{k}]$ that sends
$I^\mb{k}\times{0}$ to $H^\mb{k}_0$ and  $I^\mb{k}\times{1}$ to
$H^\mb{k}_1$. Under this homeomorphism,
the face $I^{\mb{s}, \epsilon}\times I^{\{v\}}$ goes
to $T^{\mb{s},\epsilon}$, and $W^\mb{k}$ goes to $\partial I^\mb{k}\times
I^{\{v\}}$.
\end{prop}
\begin{proof}
This is a standard assertion for $\PL$ topology. The required homeomorphism
can be built using Alexander's trick.
\end{proof}

\subsection{}
Now we will construct a family $\Omega(\Gamma) = \{\Omega_a\}_{a \in
\bcc{2N-2}^{\mb k}}$ of closed subpolyhedra in $X\times T^\mb{k}$.
Let $x^* \in X$. Put $\Omega^0_a(x^*)=\{(t_\mb{k})\,|\,(x^*,t_\mb{k})\in
\Omega_a^0\}\subset I^\mb{k}$. Define
\begin{equation}\label{equ52}\Omega_a(x,v)=
\Omega^0_a(x)\oplus 2vI^{o(a)}+2[\frac{a}{2}]v \subset \R^\mb{k},
\end{equation}
where $[\cdot]$ is the integral part. By the agreement in Sec.~\ref{ss831},
we regard $2[\frac{a}{2}]$ as an integer vector. The symbol
$\oplus$ means the Minkowski sum, and
$$o(a)=\{i\in \mb k\,|\, a_i \text{ is odd}\}.$$
Define a map
\begin{align}\Omega_a(x,v)\xar{F_a\lfloor_{x,v}} \Omega^0_a(x) \\
\intertext{ as follows: }
F_a\lfloor_{x,v}(t_\mb{k})=(t_{e(a)},
f^{\frac{O(a)-1}{2}}(x, t_{e(a)})),\end{align}
where
$$e(a)=\{i\in \mb k\,|\, a_i \text{ is even}\}$$ and
$O(a)=a|_{o(a)}:o(a)\xar{}\N$. Define
\begin{equation}\label{equ53}\Omega_a=\{(x,t_\mb{k},v)\,|\, t_\mb{k}\in \Omega_a(x,v)\}.\end{equation}
Put
\begin{equation}\label{equ54}\Omega_a \xar{F_a} \Omega_a^0:
F_a(x,t_\mb{k},v)=F_a\lfloor_{x,v}(t_{\mb{k}}).\end{equation}
(According to  (\ref{equ51}), the image of $F_a$ belongs to
$\Omega_a^0$.)
\subsection{}\label{ss837}
Consider the sets $\Omega^1_a(x)=\Omega_a(x,1) \subset
\br{2N-1}^\mb{k}$. For $a\in \bcc{2N-2}$, define an interval $\Lambda_a$ in
$\br{2N-1}$ by
$$\Lambda_a=\begin{cases} [2l,2l+1] & \text{ for } a=2l, \\
[2l-1,2l+1] & \text{ for } a=2l-1. \end{cases}$$
For $a\in \bcc{2N-2}^\mb{k}$, define a parallelepiped $\Lambda_a$ by
the formula $\Lambda_a=\prod_{i\in \mb k}\Lambda_{a_i}$. The set
$\{\Lambda_a\}_{a\in \bcc{2N-2}^\mb{k}}$  of closed
parallelepipeds is the base
of some  $\CPL$ Alexandroff topology $\Lambda (N,\mb k)$ on
$\br{2N-1}^\mb{k}$.
\begin{prop} \label{pro53}
For any $a \in \bcc{2N-2}^\mb{k}$, $x\in X$, the inclusion $\Omega^1_a(x) \subset
\Lambda_a$ holds.
\end{prop}
\begin{proof} This follows immediately from (\ref{equ52}).
\end{proof}

\subsection{}
Here we will proof the following proposition.
\begin{prop}\label{pro54}\mbox{}\\
{\rm1.} The set $\Omega(\Gamma)$ of polyhedra  is a cover of $X\times T^\mb{k}$.\\
{\rm2.} Setting $F_\Gamma|_{\Omega_a}=F_a$, we obtain a
well-defined map $X\times T^\mb{k}\xar{F_\Gamma} X\times I^\mb{k}$ that
is fiberwise with respect to the projection to $X$.
\end{prop}
\begin{proof}
The proof proceeds by induction on the number of elements in $\mb{k}$.

\smallskip
1. Let $\# \mb{k}=1$. In this case, we are in the classical situation of
Hudson's proof \cite[Theorem~6.2, p.~130]{Hu}. For $x^*\in X$, the sets
$\Omega_j(x^*)$ form a cover of  $T^1$ by parallelograms and
triangles (see Fig.~\ref{e57}).
\begin{figure}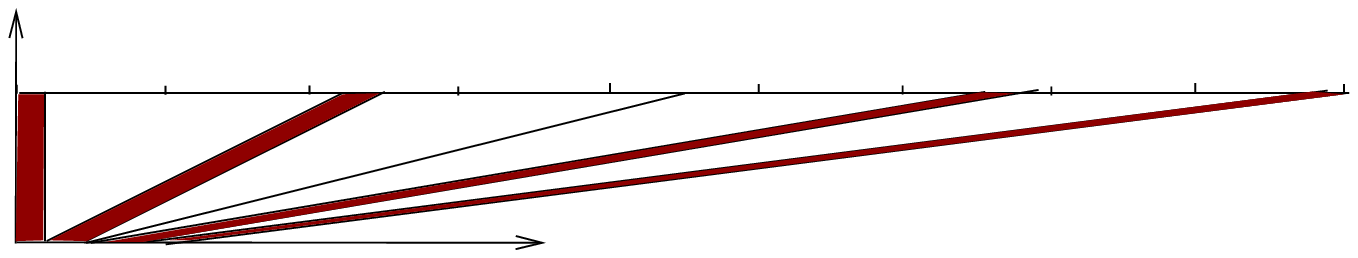\caption{}\label{e57}\end{figure}
Sets of the type $\Omega_{2l}(x^*)$ are parallelograms, sets
of the type  $\Omega_{2l-1}(x^*)$ are triangles. Under the maps $F_a$,
a parallelogram is projected to its base along the side edge, and a
triangle is projected into one of its vertices. All these maps are
automatically pasted into a global map $F_\Gamma$.

\smallskip
2. Inductive step.

\smallskip
2.1 We can build $2N$ new functions $T^{\mb{k}\setminus
\{i\}}\xar{\ol f_i^j} \R^{\{i\}}$ with indices from $\bcc{2N-1}$ as follows:
\begin{equation}\label{equ55}\begin{array}{llll}  \ol f_i^{2l}(\Gamma) &= & f_i^l\circ F_{\delta_i^l \Gamma}(x,t_{\mb{k}\setminus \{i\}},v)+2lv & \text{ for } l=0,{\ldots} ,N-1,\\
 \ol f_i^{2l-1}(\Gamma) &= & f_i^l\circ F_{\delta_i^l \Gamma}(x,t_{\mb{k}\setminus \{i\}},v)+(2l-2)v & \text{ for } l=1,{\ldots} ,N. \end{array}\end{equation}
By the inductive assumption,  these functions are well-defined.
Put $$\ol \Gamma_i^j = \{(x,t_\mb{k},v)\,|\, t_i= \ol f_i^{j}(x,t_{\mb
k\setminus \{i\}},v)\}.$$ By construction,
\begin{equation}\label{equ56}\ol f_i^{2l}\equiv\ol f_i^{2l-1}+2v.\end{equation}
The collections  $\ol f_i^j$ are monotone for every $i\in \mb{k}$:
$$0\equiv\ol f^0_i\leq \ol f_i^1 \leq {\ldots}  \leq \ol f^{2N-1}_i\equiv 1+2(N-1)v. $$
For $\mb s \xar{c} \bcc{2N-1}$, define a map $X\times T^{\mb k
\setminus \mb s}\xar{\ol f^c} \R^{\mb s}$ by the following rule:
\begin{equation}\label{equ61}\ol f^c(x,t_{\mb k \setminus \mb s},v)=f^{[\frac{c+1}{2}]}\circ
F_{\delta^{[\frac{c+1}{2}]}\Gamma}(x,t_{\mb k \setminus \mb s},v)+2[\frac{c}{2}]v.\end{equation}
We obtain
\begin{equation}\label{equ62}\Gamma \ol f^c = \cap_{i\in \mb s} \ol \Gamma^{c(i)}_i.\end{equation}

\smallskip
2.2. We must verify that the definition~(\ref{equ53}) has the
following property: for every $a\in \bcc{2N-2}$,
\begin{equation}\label{equ57}\Omega_a=\{(x,t_\mb{k},v)\,|\, \ol f_i^{a_i}
(x,t_{\mb k\setminus \{i\}},v)
 \leq t_i \leq \ol f_i^{a_i+1}(x,t_{\mb k\setminus \{i\}},v)  \}.\end{equation}
Thus $\Omega_a$ is a cover of $X\times T^\mb{k}$.

\smallskip
2.3. Let us construct a map
$$X\times T^\mb{k}\supset  \bigcup_{i\in \mb k,\, j \in \bcc{2N-1}}
\ol\Gamma_i^j \xar{\wt F}
\bigcup_{i\in \mb k,\, j \in \bcc{N}} \Gamma_i^j \subset {X\times I^\mb{k}}. $$
We define  $\wt F|_{\ol \Gamma_i^j}$  by the following rule:
\begin{equation}\label{equ58}(x,t_\mb{k},v)\mapsto f^{[\frac{j+1}{2}]} \circ F_{\delta_i^{[\frac{j+1}{2}]} \Gamma}(x,t_{\mb{k}\setminus \{i\}},v).\end{equation}
We must verify that the map
 $\wt F|_{\ol \Gamma_i^j}$ is well-defined. To this end, we must
check that the rule~(\vref{equ58}) has the following property: if
$(x^*,t^*_\mb{k},v^*) \in \ol \Gamma_{i_1}^{j_1}\cap
\ol\Gamma_{i_2}^{j_2}$, then $\wt F|_{\ol
\Gamma_{i_1}^{j_1}}(x^*,t^*_\mb{k},v^*) =\wt F|_{\ol
\Gamma_{i_2}^{j_2}}(x^*,t^*_\mb{k},v^*)$. There are two cases:
$i_1\neq i_2$ and $i_1=i_2$. If $i_1\neq i_2$, then we use
(\ref{equ62}), (\ref{equ61}). If $i_1=i_2$, then we should directly
unfold (\ref{equ58}).

\smallskip
2.3. From (\ref{equ57}) it follows that for any polyhedron $\Omega_a$ its
boundary $\partial \Omega_a$
is contained in
 $\bigcup_{i\in \mb k,\, j \in \bcc{2N-1}}\ol\Gamma_i^j$,
and for any boundary point $(x,t_\mb{k},v)\in \partial \Omega_a$ we have
$F_a(x,t_\mb{k},v)=\wt F(x,t_\mb{k},v)$. This observation completes
the inductive step.
\end{proof}
\subsection{}
We will separately mention
some facts about the map $F_\Gamma$. Let
$(x^*,t^*_\mb{k},v^*) \in T^\mb{k}$, $v^*>0$. Fix $i\in \mb k$,
$\varepsilon>0$ and consider the set of points
$$V_\epsilon^i=\{(x^*,t^*_{\mb k \setminus \{i\}},t_i,v^*)|
\mod{(t_i-t^*_i)}<\varepsilon\}. $$
\begin{prop} \label{pro55}
If any of the conditions
\begin{itemize}
\item[{\rm(a)}] $(x^*,t^*_\mb{k},v^*)\in \int \Omega_a$ and $a_i=2l-1$; or
\item[{\rm(b)}] $(x^*,t^*_\mb{k},v^*)\in \Omega_a$, $a_i=2l$, and
$\ol f^{a_i+1}(x^*,t^*_\mb{k},v^*)=\ol f^{a_i}(x^*,t^*_\mb{k},v^*)$
\end{itemize}
is satisfied, then there exists $\varepsilon$ such that
$F_\Gamma|_{V_\varepsilon^i}\equiv\mr{const}$.
\end{prop}
\begin{proof}
1. Assume that $a\in\bcc{2N-2}$ and $a_i$ is odd. Consider the line
$L_i(x^*,t^*_\mb{k},v^*)$ passing through $(x^*,t^*_\mb{k},v^*)$ and
parallel to the $t_i$-th coordinate axis. By the definition of
$\Omega_a$, we have
$$L_i(x^*,t^*_\mb{k},v^*)\cap \Omega_a=\{(x^*,t^*_{\mb k \setminus \{i\}},t_i,v^*)\,|\,
\ol f^a_i(x^*,t^*_{\mb k \setminus \{i\}},v^*) \leq t_i
\leq \ol f^{a_i+1}(x^*,t^*_{\mb k \setminus \{i\}},v^*)\}.$$
By the definition of $F_a$, we have
$F_a|_{L_i(x^*,t^*_\mb{k},v^*)\cap \Omega_a}\equiv \mr{const}$.

\smallskip
2. According to (\vref{equ56}), the length of
$L_i(x^*,t^*_\mb{k},v^*)\cap \Omega_a$ is equal to $2v^*$.\\

\smallskip
3.1. If (a) is true, then the conclusion of the proposition is true
by Steps~1 and~2 of the proof.

\smallskip\noindent
3.2. If (b) is true, then, by (\ref{equ56}), the point
$(x^*,t^*_\mb{k},v^*)$ lies in $\Omega_b\cap\Omega_c$, where
$$b_j=\begin{cases}a_i-1 & \text{ for } j=i,\\
                   a_i    & \text{ for } j\neq i
\end{cases}\quad\mbox{and}\quad c_j=\begin{cases}a_i+1 & \text{ for } j=i,\\
                   a_i    & \text{ for } j\neq i.
\end{cases}
$$
According to Steps~1 and~2 of the proof, the map $F_\Gamma$ is constant on
both  intervals
 $L_i(x^*,t^*_\mb{k},v^*)$ of length $2v^*$
  starting at the point $(x^*,t^*_\mb{k},v^*)$.
Therefore the conclusion of the proposition is true.
\end{proof}
\subsection{}\label{ss8311}
Here we will consider the case where the map $F_\Gamma$ has a block
structure (see Sec.~\vref{ss823}). We introduce the following
assumption on graph systems
$\Gamma \in \mb{GF}(X,\mb{k},N)$:

\medskip\noindent
{\bf gf3:} Let $M$ be a positive integer. We assume that $N$ is
a multiple of $M$:  $N=MN'$. Further, we assume that the
graphs $\Gamma_i^{lN'}$ in the system  $\Gamma$ are the graphs of constant functions
$f_i^{lN'}\equiv\frac{l}{M}$ for
$i\in \mb{k}$, $l=0,{\ldots} ,M$.

\medskip
Denote by $\mb{GF}(X,\mb{k},N,M)$ the set of all graph systems
satisfying {\bf gf3}. On the unit interval $I$, consider the
structure of the ball complex $\frac{1}{M}[I]=(I, \ol\Psi)$.
Recall that
we identify a ball complex with the $\CPL$ Alexandroff space
whose minimal base is the set of balls. We parameterize the set of balls  $\mf
B(\ol\Psi)$ by the set $\wt{\bcc{2M-1}}=\{-1,{\ldots} ,2M-1\}$.
Put
$$\begin{array}{lll}
\ol\Psi_{2l}&= & [\frac{l}{M},\frac{l+1}{M}]\subset I,\quad l=0,{\ldots} ,M-1, \\
\ol\Psi_{2l-1}&= & \{\frac{l}{M}\}\in I,\quad l=0,{\ldots} ,M.
\end{array}$$
Now extend  $\frac{1}{M}[I]$ to a $\CPL$ Alexandroff topology
 $ \ol\Theta$
on $T^1$. The elements of the minimal base $\mf B ( \ol\Theta)$
 are indexed by the same set $\wt{\bcc{2M-1}}$.
Put
 \begin{equation}\label{e63}\begin{array}{lll}
 \ol\Theta_{2l}&=&\begin{cases}
 \conv((\frac{l}{M},0),(\frac{l+1}{M},0)(2lN',1),(2(l+1)N',1))  \text{ for } l=0,{\ldots} ,M-2,\\
 \conv((\frac{M-1}{M},0),(1,0),(2(M-1)N',1)(2N-1)) \text{ for } l=M-1;\end{cases}\\
 \ol\Theta_{2l-1}&=&\begin{cases}\conv((\frac{l}{M},0),(2lN',1))\text{ for } l=0,{\ldots} ,M-1, \\
                              \conv((0,1),(2N-1)) \text{ for } l=M.\end{cases}
 \end{array}\end{equation}
Here is an illustration of these ugly formulas:
\begin{equation}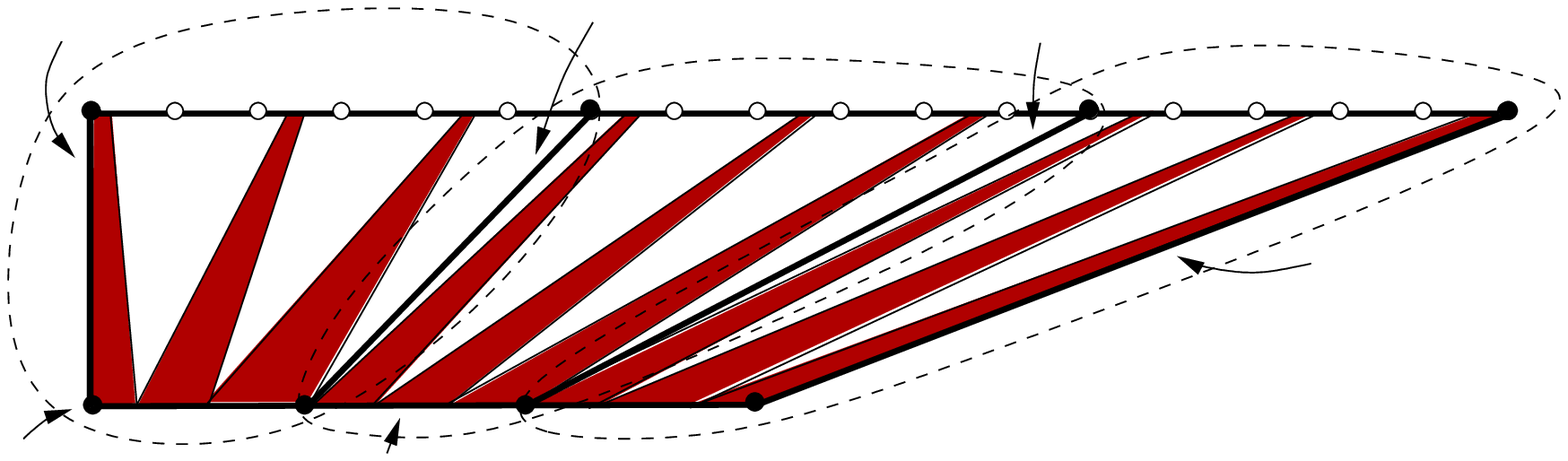\end{equation}
Consider the poset morphism $\mf B(\ol\Theta) \xar{\ol\zeta}\mf
B(\ol\Psi)$ that sends
 $\ol\Theta_a$ to $\ol\Psi_a$, where $a\in \wt{\bcc{2M-1}}$.
\begin{prop}\label{p58}
Let $\Gamma\in \mb{GF}(X,\mb{1},N,M)$. Then $F_{\Gamma}(X\times
\ol\Theta_a)=X\times \ol\zeta(\ol\Theta_a)=X\times\ol\Psi_a$.
\end{prop}
\begin{proof} Consider the triangles
$\Omega_{2jN'-1}(x) \subset T^1$, $j=1,{\ldots} ,M-1$. According to {\bf
gf3}, we have $f^{jN'}\equiv \frac{j}{M}$. From this fact and
(\vref{equ57}) it follows that for every $x\in X$ \begin{multline*}
\pi_{v,t}(\Omega_{2jN'-1}(x))\equiv \Omega_{2jN'-1}(x_2)\\=
\left\{(t,v)\,|\,\frac{j}{M}\leq t+2(jN'-1)v \leq \frac{j}{M}+2jN'v\right\}\supset
\ol\Theta_{2j-1}. \end{multline*} Comparing this equation with the definition
of  $F_{\Gamma}$, we obtain the required assertion.
\end{proof}
In the case of a general $\mb k$, we consider the complex
$\frac{1}{M}[I^\mb{k}]=(\frac{1}{M}[I])^\mb{k}= (I^k,\ol\Psi)$, where
$\mf B(\ol\Psi)=\{\ol\Psi_{a}\}_{a\in \wt{\bcc{2M-1}^\mb{k}}}$,
$\ol\Psi_a=\ol\Psi_{a_1}\times{\ldots} \times \ol\Psi_{a_k}$. Define an
Alexandroff space $(T^\mb{k}, \ol\Theta)$ as a
product fibered by $v$:
$\mf B(\ol\Theta)=\{\ol\Theta_{a}\}_{a\in
\wt{\bcc{2M-1}^\mb{k}}}$,
$$\ol\Theta_a=\{(t_\mb{k},v)\,|\,(t_i,v)\in \ol\Theta_{a_i}, i=1,{\ldots} ,k\}.$$
Define a morphism $\mf B(\ol\Theta)\xar{\ol\zeta} \mf B(\ol\Psi)$ by
the rule  $\ol\zeta \ol\Theta_a=\ol\Psi_a$. Applying Proposition~\ref{p58}
to the components of the fibered product, we obtain the following proposition.
\begin{prop}\label{pro58}
Let $\Gamma\in \mb{GF}(X,\mb{k},N,M)$. Then for every
$a \in \wt{\bcc{2M-1}}^\mb{k}$,
$$F_{\Gamma}(X\times \ol\Theta_a)=X\times \ol\zeta(\ol\Theta_a)=X\times\ol\Psi_a.$$
\end{prop}
The elements of the minimal bases $\mf B (\Theta)=\{\Theta_a\}$ and $\mf
B (\Psi)=\{\Psi_a\}$ of the dual $\OPL$ Alexandroff topologies
 $\Theta =\ol\Theta^c$ on $T^\mb k$ and $\Psi = \ol\Psi^c$ on $I^\mb k$
are indexed by the same set $\wt{\bcc{2M-1}}$. Consider the poset
morphism $\mf B(\Theta) \xar{\zeta}\mf B(\Psi)$ that sends
$\Theta_a$ to $\Psi_a$. By the duality and Proposition~\ref{pro58},
we conclude that the following proposition holds.
\begin{prop}\label{pro59}
Let $\Gamma\in \mb{GF}(X,\mb{k},N,M)$. Then for every
$a \in \wt{\bcc{2M-1}}^\mb{k}$,
$$F_{\Gamma}(X\times \Theta_a)=X\times \zeta (\Theta_a)=X\times\Psi_a.$$
\end{prop}
\section{Homotopies of fiberwise homeomorphisms associated
 with graph systems} \label{hom_syst}
Let  $G \in \PL_{I^\mb{k}}{X}$ be a fiberwise homeomorphism. Let
$\Gamma \in \mb{GF}(X,\mb{k},N)$ be a graph system. We say that  $G$
and  $\Gamma$ are  \bem{consistent} if $G^{-1}(\Gamma)\in
\mb{GF}(X,\mb{k},N)$.
\begin{lemma}  \label{l7}
If $G$ and $\Gamma$ are consistent, then $(F_\Gamma\rtimes
G^{-1})^{-1}$ is a fiberwise homeomorphism and $(F_\Gamma\rtimes
G^{-1})^{-1}=F_{G^{-1}(\Gamma)}\rtimes G$.
\end{lemma}
\begin{proof} The equality
 $(F_\Gamma\rtimes G^{-1})^{-1}=F_{G^{-1}(\Gamma)}\rtimes G$ is a direct
 corollary
of the construction of  $F_{\Gamma} $ (see Proposition~\vref{pro54}).
\end{proof}
\subsection{}
We impose an additional assumption on graph systems
from $\Gamma\in \mb{GF}(X,\mb{k},N)$.

\medskip\noindent
{\bf gf4:} The functions $f_i^j(x,t_{\mb{k}\setminus \{i\}})$ are
independent both from $t_{\mb{k}\setminus \{i\}}$ and the index $i$.
That is, for every $j\in \bcc{N}$ there exists a function $X\xar{f^j} I$
such that for every
$i \in \mb k$ the equality $f_i^j(x,t_{\mb{k}\setminus \{i\}})=f^j(x)$ holds.

\medskip
Denote by $\mb{HGF}(X,\mb{k},N)$ the set of all $\Gamma\in
\mb{GF}(X,\mb{k},N)$ satisfying {\bf gf4}. Denote by
$\mb{HGF}(X,\mb{k},N,M)$ the set of all $\Gamma\in
\mb{GF}(X,\mb{k},N,M)$ satisfying {\bf gf4}.
\subsection{Flat points of a fiberwise homeomorphism}\label{ss81}
We will need some general facts about ``flat points'' of a fiberwise
homeomorphism. Let $G\in \PL_B(X)$ and $(x,b)\in X\times B$.

\medskip\noindent
{\em We say that $(x,b)$ is a} \bem{flat point of $G$} {\em if
for some open neighborhood $V\subset B$ of the point $b$ the following is true:
for any  $b' \in V$ the equality $G^{-1}\lfloor_{b'}(x)=G^{-1}\lfloor_{b}(x)$ holds.}

\medskip
We can rephrase the definition of a flat point in prismatic terms.
Consider the Alexandroff topology $\delta_x$ on $X$
with minimal
base  $\mf{B}(X,
\delta_x)=\{\{X\},\{x\}\}$ consisting of two elements. Consider the constant Alexandroff presheaf
$\zeta_{V,x}$ that sends $(V, V^{\mr{triv}})$ to $(X,\delta_x)$.
\begin{prop}
The point $(x,b)$ is a flat point of $G$ if and only if for some neighborhood
$V \subset B$ of the point $b$ the homeomorphism $G^{-1}\lfloor_{V}$ is
$\zeta_{v,x}$-prismatic.
\end{prop}\label{pr48}
Let $\mb{k}=\{1,{\ldots} ,k\}$. Pick $G \in \PL_{\R^\mb{k}}(X)$,
$(x,b)\in X\times \R^\mb{k}$, and $i \in \mb{k}$. Let $e_1,{\ldots} ,e_k$
be a basis in
 $\R^\mb{k}$. Consider the line
 $l_i(b)=\{b+e_it\,|\, t\in \R\}$. Thus $l_i(b)$ is the line
 passing through the point  $b$ parallel to $e_i$.

\medskip\noindent\label{def1}
{\em We say that} \bem{$G$ is horizontal in the direction $i$ at
the point  $(x,b)$} {\em if there exists a neighborhood $V
\subset \R^\mb{k}$ of the point $b$ such that for every $b'\in V$
the homeomorphism $G\lfloor_{l_i(b')}$
is flat at $(x,b')$.}

\begin{prop} \label{pr49}
If $G\in \PL_{\R^\mb{k}}(X)$ and $(x,b)\in X\times \R^\mb{k}$ are
such that $G$ is horizontal for all directions $i\in\mb{k}$, then
$(x,b)$ is a flat point of $G$.
\end{prop}

\begin{proof}
We may assume that there is a small cube $Q \subset \R^\mb{k}$
with edges parallel to the basis vectors and barycenter $b$
such that for every $b'\in Q$ the homeomorphism $\PL\lfloor_{l_i(b')}$ is
horizontal at $(x,b')$ in any direction $i\in \mb{k}$. We can
approach any point $b'\in Q$ by a finite polygonal line
lying inside  $Q$ and formed by
intervals  parallel to the coordinate axes. From this fact and
the conditions of the proposition it follows
that $G^{-1}\lfloor_{b'}(x)=G^{-1}\lfloor_b(x)$.
\end{proof}

\subsection{}

Denote by $(T^\mb{k})^{v^*}$ the set $\{(t_\mb{k},v)\in
T^\mb{k}\,|\,v=v^*\}$.

\begin{prop} \label{ppro59}
Let $\Gamma\in \mb{HGF}(X,\mb{k},N)$ be a graph system. Let $G\in
\PL_{I^\mb{k}}$ be a fiberwise
homeomorphism consistent with  $\Gamma$. Assume that $(x^*,t^*_\mb{k},v^*)\in X\times
T^\mb{k}$, $v^*>0$. Then the following is true.

If the fiberwise homeomorphism
  $F_\Gamma \rtimes G^{-1} \lfloor_{(T^\mb{k})^{v^*}}$
\bem{ is not horizontal} in the direction $i$ at the point
$(x^*,t^*_\mb{k},v^*)$, then there exists an index $a\in \bcc{2N-1}$
such that $a_i=2l$, $(x^*,t^*_\mb{k},v^*)\in \Omega_a$, and $x^*\in
\supp(f_i^{l+1}-f_i^l)$.
\end{prop}
\begin{proof} We consider two cases.

\smallskip
1. Assume that for every $a$ such that $(x^*,t^*_\mb{k},v^*)\in
\Omega_a$ it is true that  $a_i$ is odd. This means that
$(x^*,t^*_\mb{k},v^*)$ belongs to the complement of the closed set
$\cup_{a\in \bcc{2N-2},\,i\in e(a)}\Omega_a$. By the proof of
Proposition~\vref{pro54}, we can conclude that $a$ belongs to the
interior of some $\Omega_a$ such that $i\in o(a)$.  Then from the
definition of
  $F_a$, the definition of
$F_\Gamma\rtimes G^{-1} $, and Proposition~\vref{pro55}\,(a) it follows
that $F_\Gamma\rtimes G^{-1} $ is horizontal
at the point $(x^*,t^*_\mb{k},v^*)$ in the direction $i$.

\smallskip
2. Let $a \in \bcc{2N-1}$ be such that $a_i=2l$ and $x^* \not\in
\supp(f_i^{l+1}-f_i^l)$. In this case, by the proof of Proposition~\vref{pro54}, the point
 $t^*_\mb{k}$ belongs to
$\Omega_i^{2l-1}(x^*,v^*) \cap \Omega_i^{2l+1} $, where
$$\Omega_i^{j}(x^*,v^*)=\{t_\mb{k}\,|\,\ol f_i^j(x^*,t_{\mb k \setminus \{i\}},v^*)\leq t_i \leq \ol f_i^{j+1}(x^*,t_{\mb k \setminus \{i\}},v^*)\}.$$
The special feature  of the condition {\bf gf4} is that in this case
$$\ol f_i^j=f_i^{[\frac{j+1}{2}]}(x^*)+2[\frac{j}{2}]v^* \equiv \mr{const}.$$
Therefore in some neighborhood of the point  $(x^*,t^*_\mb{k},v^*)$
in $T^\mb{k}(x^*,v^*)$, all the points satisfy the conditions of
Proposition~\vref{pro55}. This completes the proof.
\end{proof}
\subsection{}\label{843}
Let $G$ be consistent with  $\Gamma$. We will describe the structure of
the homeomorphism $(F_\Gamma\rtimes G^{-1})^1 \in
\PL_{\br{2N-1}^\mb{k}}(X)$.

Let us introduce an $\mr{OPL}$
Alexandroff topology
 $\mc E(N,k)$ on the cube $\br{2N-1}^\mb{k}$.
First, consider an Alexandroff topology $\mc E(N,1)$ on
$\br{2N-1}^1=[0,2N-1]$. The elements of the minimal base
$\mf{B}(\br{2N-1}^1, \mc E(N,1) )$ are indexed by the elements of
$\bcc{2N-2}^1 =\{0,1,{\ldots} ,2N-2\}$. Put
$$E_i =
\begin{cases} (i-1,i+2) & \text{ if } 2 \leq i \leq 2N-4 \text { and $i$ is even},\\
              (i,i+1)   & \text{ if } 1 \leq i \leq 2N-3  \text { and $i$ is odd},\\
              (0,2) & \text{ if } i=0,\\
              (2N-3,2N-1) & \text{ if } i =2N-2.
 \end{cases}
$$
Consider the partial order $\us{e}{\leq}$ on  $\bcc{2N-2}^1$ generated
by the relation
 ``an even number is larger than both neighboring odd numbers.''
 We can see that
 $E_i\subseteq E_j $ if and only if $i\us{e}{\leq}j$.
Thus we obtain a canonical isomorphism of posets
\begin{equation}\label{equ67}\mf{B}(\br{2N-1}^1,
 \mc E(N,1) )\approx \bcc{2N-2}^1_{\us{e}{\leq}}.\end{equation}
Define
\begin{equation}\label{equ68}(\br{2N-1}^\mb{k}, \mc E(N,\mb k))=
 (\br{2N-1}^1, \mc E(N,1) )^\mb{k}.\end{equation}
Taking a power of the isomorphism (\ref{equ67}), we obtain a
canonical isomorphism
\begin{equation}\label{equ69}\mf{B}(\br{2N-1}^1, \mc E(N,\mb k) )
\approx \bcc{2N-2}^\mb{k}_{\us{e}{\leq}}\end{equation}
(here we denote the power $(\us{e}{\leq})^\mb{k}$
 of the partial order $\us{e}{\leq}$
by the same symbol $\us{e}{\leq}$). Thus the elements of the minimal
base
 $\mf{B}(\br{2N-1}^\mb{k}, \mc E(N,k))$
are indexed by the elements of $\bcc{2N-2}^\mb{k}$, and
$E_{a_\mb{k}}=\prod_{i\in \mb k} E_{a_i}$, $E_a \subset E_b
\Leftrightarrow a \us{e}{\leq} b$.

Let $\Gamma \in \mb{GF}(X,\mb{k},N)$ be a graph system and assume that the axiom
{\bf gf4} is satisfied. Pick $a\in \bcc{2N-2}$. Consider the subset
\begin{equation}\wt J_\Gamma(a)\subset X, \quad\wt J_\Gamma(a)=\bigcup_{i\in e(a),\,a_i=2l_i}
\supp(f_i^{l_i+1}-f_i^{l_i}).\end{equation}
Put $J_\Gamma(a)=\mc D_{\wt J_\Gamma(a)}$ (the functor $\mc D_\star$
is defined in Sec.~\vref{mcd}). We obtain a poset morphism
$$\bcc{2N-2}^{\mb k}_{\us{e^{\mb k}}{\leq}}\xar{ J_\Gamma} \mb{Al}(X).$$
Using the isomorphism~(\ref{equ69}), we define a morphism
$$(\br{2N-1}^\mb{k}, \mc E(N,\mb k))\xar{\mb J_\Gamma} \mb{Al}(X)$$
as the composite
$$ (\br{2N-1}^\mb{k}, \mc E(N,\mb k))\xar{o}\mf B (\br{2N-1}^\mb{k}, \mc E(N,\mb k))
\approx \bcc{2N-2}^\mb{k}_{\us{e}{\leq}} \xar{J_\Gamma} \mb{Al}(X).$$
\subsection{The prismaticity of $F_\Gamma\rtimes G^{-1}$}
Here we will assemble the prismaticity properties of
$F_\Gamma\rtimes G^{-1}$ into a final lemma.

\label{844}
\begin{prop} \label{pro60}
Let $\Gamma \in \mb{GF}(X,\mb{k},N)$ be a graph system
satisfying {\bf gf4}.
If $G$ is consistent with $\Gamma$, then the fiberwise homeomorphism
$(F_\Gamma\rtimes G^{-1})^1$ is $\mb J_\Gamma$-prismatic.
\end{prop}
\begin{proof}
Note that
$$\mc E(N,\mb k)=\Lambda(N,\mb k)^c$$
(see Secs.~\vref{ss837} and~\vref{ss514}). It follows from
Proposition~\vref{pro53} that for $(x^*,t^*_\mb{k})\in X\times
\br{2N-1}^\mb{k}$ the following is true. If $t^*_\mb{k}\in E_a$ and
$(x^*,t^*_\mb{k})\in \Omega_c^1(\Gamma)$, then $c\us{e}{\leq}a$.
Comparing this fact with Propositions~\vref{pro53} and~\vref{pr49},
we can see that the following assertion is true: if $t^*_\mb{k}\in E_a$ and $(F_\Gamma\rtimes
G^{-1})^1$ is not horizontal at the point
  $(x^*,t^*_\mb{k})$,
then $x^*$ is contained in $\wt J_\Gamma(a)$. This means exactly the
required $\mb J_\Gamma$-prismaticity of $(F_\Gamma\rtimes
G^{-1})^1$.
\end{proof}

Thus the following objects are associated
with a graph system  $\Gamma \in \mb{HGF}(X,\mb k,N,M)$:
\begin{itemize}
\item an $\OPL$ Alexandroff topology $\Theta=\Theta(M)$ on
$T^\mb{k}(N)$ and a morphism
$$\mf B(\Theta)\xar{\zeta} \mf B(\Psi)=\mf B(\frac{1}{M}[I]^\mb{k})^\op$$
(see Sec.~\vref{ss8311});
\item an $\OPL$ Alexandroff topology $\mc E = \mc E(N)$ on
$T^\mb{k}_1=T^\mb{k}_1(N)\approx \br{2N-1}^\mb{k}$ and an
Alexandroff presheaf
 $(T^\mb{k}_1, \mc E)\xar{\mb J_\Gamma} \mb{Al(X)}$
(see Sec.~\vref{843}).
\end{itemize}

The topology $\mc E$ on $T^\mb{k}_1$ is inscribed into the topology
$\Theta^\op\lfloor_{T^{\mb k}_1}$ on $H^\mb{k}_1$. The following lemma holds.
\begin{lemma}\label{lem8}
If a fiberwise homeomorphism  $G$ is consistent with $\Gamma \in
\mb{HGF}(X,\mb{k},N,M)$,
then
\begin{itemize}
\item[{\rm(a)}] $F_\Gamma \rtimes G^{-1}\lfloor_{H_1^{\mb k}}$ is
$\mb J_\Gamma$-prismatic;
\item[{\rm(b)}] if $\mc L\in \mb{Al}(X)$ and
 $G$ is $\delta_{\mc L}$-prismatic, then
 $F_\Gamma \rtimes G^{-1}$ is  $\delta_{\mc L}$-prismatic;
\item[{\rm(c)}] if $G\lfloor_{\partial I^\mb{k}}$ is prismatic with respect to an
Alexandroff presheaf
 $( \partial I^\mb{k}, \Psi\lfloor_{\partial I^\mb{k}} )
\xar{\xi} \mb{Al}(X)$, then $F_\Gamma \rtimes
G^{-1}\lfloor_{T_1^\mb{k}}$ is prismatic with respect to the Alexandroff
presheaf
\begin{equation} \label{eq78}(T^{\mb{k}_1},\Theta\lfloor_{T_1^\mb{k}})\xar{\xi\circ\zeta|_{\mf B(\Theta^\op\lfloor_{T^\mb{k}})}\circ o}\mb{Al}(X).\end{equation}
\end{itemize}
\end{lemma}
\begin{proof}
Assertion~(a) is a special case of Proposition~\vref{pro60}.
Assertion~(b) follows from Proposition~\vref{p51}.
Assertion~(c) follows from Propositions~\vref{ppro59} and~\vref{p52}.
\end{proof}

\section{Hudson functions and their graphs} \label{hud_gr}
Here we use a construction from \cite[Theorem~6.2, p.~130, and
its corollaries]{Hu}. Hudson suggests especially simple functions
needed for fragmentation.

Let $K$ be a simplicial complex triangulating  $X$ and $|K|\xar{t}
X$ be a triangulation. Let $f:|K|\xar{} I$ be a function that is
linear on all simplices of $K$ (i.e., $f$ is determined by its values at
the vertices of $K$). We call such functions \bem{Hudson functions
on the simplicial complex $K$}. A function $X\xar{g}I$ that
decomposes as
$$\bfig
\Vtriangle/>`>`<-/<300,300>[X`I`|K|,;g`t^{-1}`f]
\efig$$
where $f$ is a Hudson function on $K$ will be called a \bem{Hudson function
on the triangulation $t$ of $X$} .

Denote by $\mr{Hud}(t)$ the set of all Hudson functions on $t$.
\bem{The diameter $\mr{diam} f$ of a Hudson function} $f$ is defined
as
$$\max_{\{x,y \in X\}} |f(x)-f(y)|.$$
The diameter measures the difference between  $f$ and a
constant
function. The set $\mr{Hud}_\delta (t)$ of all Hudson functions
of diameter at most $\delta$ is a compact finite-dimensional polyhedron (obviously, it is just a cube).

Let $\mc G$ be a  \bem{compact family of homeomorphisms in
$\PL_{I^k}(X)$}, i.e., $\mc G = \{G_b\}_{b\in B} \subset
\PL_{I^\mb{k}}$, where $B$ is a compact polyhedron. The map $B\times
X \times I^\mb{k}\xar{\mc G} X \times I^\mb{k}$ defined by the
correspondence $(b,x, t_\mb{k}) \mapsto G_b(x,t_\mb{k})$ is
supposed to be piecewise linear.

Pick $\mb s \subseteq \mb k$. Consider a collection of functions
$f_{\mb s}, f_i \in \mr{Hud}(t)$ indexed by $\mb{s}$.  Consider the
collection of subpolyhedra $\Gamma_{\mb{s}}$, where $\Gamma_i =
{\{(x,t_{\mb k})\,|\, t_i=f_i(x)}\}$.

A homeomorphism $G\in \PL_{I^\mb{k}}(X)$ is \bem{consistent} with
$f_{\mb s}$ if both  $\Gamma_s$ and $G^{-1}(\Gamma_\mb{s})$
satisfy the condition {\bf gf} from Sec.~\vref{gf}.

\begin{prop} \label{pr45} Under our conditions, there exists
$\delta (\mc G) >0$ such that if $\mr{diam} f_i \leq \delta$  for all
 $i \in \mb s$,
then $f_{\mb{s}}$ is consistent with all homeomorphisms from $\mc
G$.
\end{prop}
\begin{proof}
We will apply Hudson's arguments inductively by the number of
functions in
$f_{\mb s}$.

\smallskip
1. Let $\# \{\mb s\} =1$ and $\mb s =\{i\}$. We must prove that
there exists $ \delta>0$ such that if $\mr{diam} f_{i} \leq \delta$,
then every leaf
$$L_{(x^*,t^*_{\mb k \setminus \{i\}})}(G)=
\{(x,t_\mb{k})\,|\,x=G(x^*)\lfloor_{t_\mb{k}},
t_{\mb k \setminus \{i\}}=t^*_{\mb k \setminus \{i\}}\}$$ of every
$G \in \mc G$ intersects $\Gamma_i$ in a single point. Consider the
projection $X\times I^{\mb{k}}\xar{\pi_{\{i\}}} I^{\{i\}}$. The
homeomorphisms from
 $\mc G$
are fiberwise with respect to $\pi_{i}$, therefore Hudson's arguments
can be directly applied.

\smallskip
2. Assume that we can prove the assertion for $\#\{\mb s\}\leq
l-1$. Let $\#\{\mb s\}=l$.

Pick $i \in \mb s$. By the previous step, we can find $\delta >0$ such
that if $\mr{diam}(f_i)\leq \delta$, then $f_i$ is consistent with
$\mc G$. Given $f_i$ and $G$, there is a homeomorphism
$W(f_i,G)\in \PL_{I^{\mb k  \setminus \{i\}}}(X)$
defined as
follows. Denote by $X\times I^{\mb k\setminus \{i\}}\xar{G^{-1}f_i}
I^{\{i\}}$ the function with graph $G^{-1} (\Gamma_i)$. Set
$$W(f_i,G)\lfloor_{t_{\mb k \setminus \{i\}}}(x)=
G\lfloor_{(t_{\mb k \setminus \{i\}}(x),
t_i=G^{-1}f_i(x,t_{\mb k \setminus \{i\}}))}(x).$$
The set of homeomorphisms
$$\mc W=\{{W(f_i, G)}\}_{f_i\in \mr{Hud}_\delta,\, G\in \mc G}
 \subset \PL_{I^{\mb k \setminus \{i\}}}(X)$$
is a compact set of homeomorphisms from $\PL_{I^{\mb k \setminus
\{i\}}}(X)$. By the inductive assumption, for this set $\delta$ exists for
the set of Hudson functions $f_{\mb s \setminus \{i\}}$. It
remains to observe that the set of Hudson functions
 $f_{\mb s \setminus \{i\}}$, which is consistent with
$\mc W$, together with the Hudson function $f_{\{i\}}$, which is consistent with $\mc
G$, forms a set $f_\mb{s}$ consistent with $\mc G$.
\end{proof}

Let $|K|\xar{t}X$ be a triangulation. Let $\dot{star}_v(t)$ be the
cover of $X$ by the open stars of $t$. Let $\mc G \in
\PL_{I^\mb{k}}(X)$ be a compact family of homeomorphisms. Let $N,M$
be positive integers. A \bem {Hudson system} for the data
$\la\mc G, N, M, t\ra$ is a graph system $\Gamma\in
\mb{HGF}(X,\mb{k},N,M)$ such that $\Gamma $ is simultaneously
consistent  with all homeomorphisms from $\mc G$ and
$\supp (f_i^{j+1} - f_i^j) \in \dot{star}_v (t)$ for all $i,j$.

As a result of the constructions of this section, we obtain
the following assertion.
\begin{lemma} \label{prop48}
For any compact family of homeomorphisms $\mc G \in
\PL_{I^\mb{k}}(X)$, any triangulation $|K|\xar{t}X$, and any positive
integer $M$ there exist a positive integer $N(\mc G)$ and a Hudson
system $$\mr H(\mc G, t,  M)\in \mb{HGF}(X,\mb{k},N(\mc G),M).$$
\end{lemma}
\begin{proof}
First, order the vertices of $K$ in an arbitrary way: ${v_1,{\ldots} ,v_{V}} $.
Choose a special sequence $\{ h_{j}\}_{j=0,{\ldots} ,V}$ of Hudson
functions on $K$. Put $h_0 \equiv 0$; for $j=1,{\ldots} ,V$, define functions
on the vertices by the following rule:
$$ h_{j} (v_i)=\begin{cases} 1 & \text{ for } i < j, \\
                                              0 & \text{ for } i \geq j.
\end{cases}  $$
It is clear that $0=h_0\leq h_1\leq{\ldots} \leq h_V$ and $\supp (h_{j+1}
- h_j)= \dot \st (v_j)$.

Further, choose  $\delta(\mc G)$ by Proposition~\vref{pr45}
and choose a positive integer $M$ such that $\frac{1}{M} < \delta$.

Now we can  define the system of functions $\{f_i^j\}_{i \in \mb k,\,
j \in \bcc{MV}} $ that determines $\mr H(\mc G, t,  M)$. Put $f^j_i
= \frac{l}{M}+\frac{1}{M}h_{j'}$ for $j = lV + j'$, $0 < j' < V$.
\end{proof}
\subsection{}\label{851}
Comparing Lemma~\ref{prop48} with Lemma~\ref{lem8}, we can see that
for a compact family of homeomorphisms $\mc G \subset \PL_{I^k}(X) $
and a Hudson system $\mr H(\mc G, t,  M)$, the image of the
Alexandroff presheaf
 $\mb J_{\mr H(\mc G, t,  M)}$ lies in the subposet
 $\dot{star}^k (t)\subset \mb{Al}(X)$ that is generated by the unions of
at most $k$  open stars of $t$.
\section{Small balls, separation of small balls} \label{balls_sep}
We need a sequence of propositions on  ``separation'' of
configurations of balls on a simplicial manifold $K$.  We say that a
set of balls $\{A_i\}_i \in \mb D^\infty(K)$ \bem{separates} a (not
necessarily disjoint) set of balls
 $\{B_j\}_j$, $B_j \subset
|K|$, if  $\cup_i B_j \underset{X}{\subseteq} \cup_i A_i $.

Consider a compact simplicial geometric $n$-dimensional
$\PL$ manifold $K$. Denote by $ star (K)$ the cover of $|K|$ by
the closed stars of $K$. Denote by $\sd_i K$ the $i$th barycentric
subdivision of $K$.

\begin{prop} \label{balls0} If $i\geq j\geq 0$, then
$ star (\sd_i K)$ is inscribed into  $ star (\sd_j K)$.
\end{prop}
\begin{proof}
This is  a general property of subdivisions of simplicial complexes:
$$K_0 \trianglelefteq K_1 \Rightarrow star (K_0) \text{ is inscribed into } star (K_1).$$
\end{proof}
\begin{prop}\label{balls1}
For any two balls $a,b \in star(\sd_2 K)$ such that
 $a \cap b \neq \emptyset $
there exists $c \in star( K)$ such that $a\cup b \subseteq c$.
\end{prop}
\begin{proof}
Recall that the vertices of $\sd_1 K$ are indexed by the simplices of
$K$, the edges of $\sd_1 K$ are indexed by the 2-flags  of simplices of
$K$, {\ldots} , the $k$-simplices of $\sd_1 K$ are indexed
by the $k$-flags of $K$.
Thus let
$a,b \in star \sd_2 K$ and $a \cap b \neq \emptyset$. Let $t$ be a
vertex of $\sd_2 K$, $t \in a \cap b$. The complex $\sd_2 K $ is
covered by $\sd_1 K$, therefore the vertex $t$ belongs to some
complete flag $s_0,{\ldots} ,s_n$. We claim that the closed star  $|\st_K
(s_0)|$ contains  $|a\cup b|$ (see Fig.~\ref{stell}).
\begin{figure}
\includegraphics[scale=0.6]{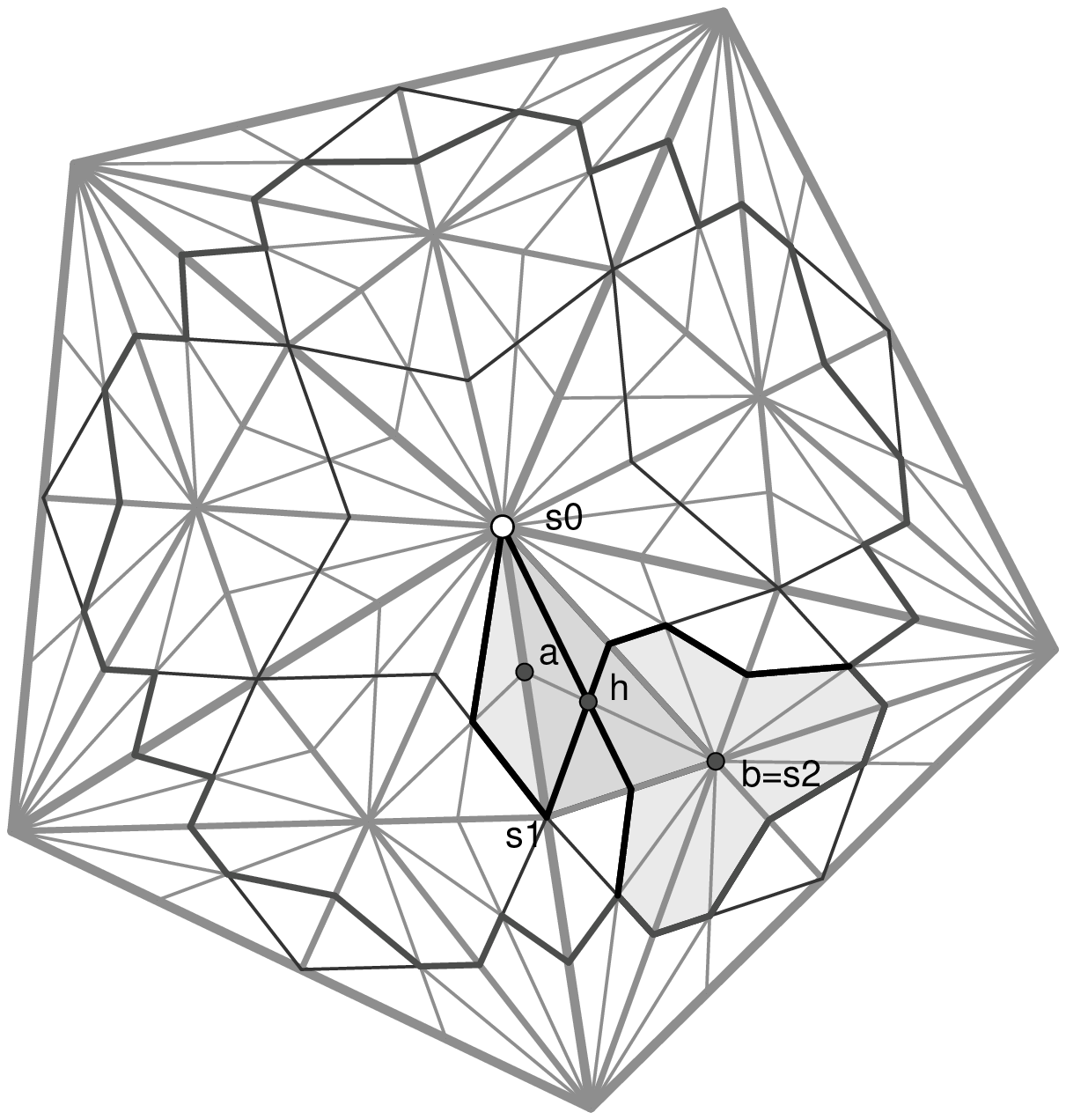} \caption{\label{stell}}
\end{figure}
 This follows from two observations:

\smallskip
1. The simplicial ball $\sd_2 \st_K (s_0) $
is a full subcomplex of $\sd_2 K$. This means
that if all vertices of some simplex from $\sd_2 K$ belong to
$\sd_2 \st_K (s_0)$, then the simplex itself belongs to  $\sd_2 \st_K
(s_0)$.

\smallskip
2. The minimum length of a chain of edges connecting $t$ with a
vertex of  $\sd_2 \lk_K (s_0)$ is at least $2$. Therefore the minimum
length of a chain of edges connecting $t$ with a vertex that does
not belong to  $\sd_2 \st_K (s_0) $ is at least  $3$. While the minimum
 length of a chain of edges connecting $t$ with a vertex of
$a \cap b$ is at most $2$.

\smallskip
It follows from the second observation that all vertices of  $a \cap b$ belong
to
 $\sd_2 \st_K (s_0)$.  Therefore, by the first observation,
all simplices of
  $a \cap b$ belong to $\sd_2 \st_K (s_0)$.
\end{proof}
It is easy to verify the following assertion.
\begin{prop}\label{prop49}
If $i \geq 2$, then the star of a simplex from  $\sd_i K$ either
is itself a ball consistent with $K$, or is contained in a star of $\sd_{i-1}K$
 consistent
with $K$.
\end{prop}

Consider the set $\bs{B} (X)$ of all closed  $\PL$  balls embedded
into $X$. Denote by $\bs{B}^\infty (X)$ the set of all finite
subsets in $\bs{B} (X)$.
There are two partial orders on  $\bs{B}^\infty (X)$:
\begin{itemize}
\item the order  \bem{by inclusion}:
$$\{B_i\}_i \sqsubseteq \{C_j\}_j
\Leftrightarrow  \{B_i\}_i \underset{\bs{B} (X)}{\subseteq} \{C_j\}_j,$$
\item the order by \bem{covering}:
$$\{B_i\}_i \preccurlyeq \{C_j\}_j \Leftrightarrow
\cup_i B_i \underset{X}{\subseteq} \cup_j C_j.$$
\end{itemize}

Consider the infinite sequence
$$K={\sd}_0 K \trianglerighteq {\sd}_1 K \trianglerighteq {\sd}_2 K \trianglerighteq {\ldots} .$$
Set $\bs{S}_l K=\cup_{i=l}^\infty \bs{\st} \sd_i K \subset
\bs{B}(|K|)$.
\begin{prop}[On separation of small consistent balls] \label{balls2}
Let us be given  $m$ $K$-consistent balls from $\bs{S}_{i} K$,
$i\geq 3m-2$. Then  they can be separated  by
 a collection of at most
$m$ disjoint $K$-consistent balls from $\bs{S}_{i-3(m-1)}
K$.
\end{prop}
\begin{proof}
Induction by $m$. If $m=1$, then the assertion is trivial.  Assume
that it is true for
 $m=k-1$, $k\geq 2$.

Let us be given a collection $A=\{a_1,{\ldots} ,a_k\}$ of $K$-consistent balls
from $\bs{S}_i K$.
If the balls from $A$ are disjoint, then the assertion is true.
Otherwise  $a_i\cap a_j \neq \emptyset$ for some pair $i,j\in
\{1,{\ldots} ,k\}$, $i\neq j$.  Applying Propositions~\ref{balls0} and~\ref{balls1}
and using the fact that $\sd_{i} K = \sd_2\sd_{(i-2)}K$,
we can find a ball $c' \in \bs{S}_{i-2}K$ such that $a_i \cup a_j
\subset c'$. By assumption, $i \geq 3k-2$ and $k \geq 2$. Therefore
the inequality $i-2\geq 3k-4 \geq 2$ holds, and
the conditions of Proposition~\ref{prop49} are satisfied for the
ball $c'$.
Therefore there exists a $K$-consistent ball $c \in
\bs{S}_{i-3}K $ that contains $c'$. Consider the new set of $K$-consistent
balls $B=\{c\}\cup (A \setminus \{a_i, a_j\})$.
By construction, the set
 $B$ covers $A$, consists of
$k-1$ elements, and belongs to $\bs{S}_{i-3}$. By the inductive
assumption, for the set of balls $B$ the assertion is true. By the
transitivity of $\preccurlyeq$, the inductive step is proved.
\end{proof}

Denote by  $\bs{S}^i_j K \subset \bs{B}^\infty (|K|)$ the set of all
unordered  collections of at most $i$ balls from
$\bs{S}_j K$. Obviously, the order  ${\sqsubseteq}$ on $\bs{B}^\infty
(|K|)$ is stronger than $\preccurlyeq$, i.e.,
 $(A\sqsubseteq B) \Rightarrow (A \preccurlyeq B)$.
Denote by
$\bs{B}^\infty (|K|)_\sqsubseteq \xar{\bs \beta}\bs{B}^\infty (|K|)_\preccurlyeq$
the morphism of weakening the order.
The poset $(\bs{S}^i_j K)_\preccurlyeq$ will be denoted by $\bs{T}^i_j K$, and
the poset $ (\bs{S}^i_j K)_\sqsubseteq$ will be denoted by the same symbol
 $\bs{S}^i_jK$.
Let $\bs{D}_j^i K \hra{} \bs{T}_j^i K$ be the subposet of $\bs{T}^i_j
K$ formed by all disjoint collections of $K$-consistent balls.

Obviously, $\bs{S}_{j_0}^{i_0} K \hra{} \bs{S}_{j_1}^{i_1}K$ and
$\bs{T}_{j_0}^{i_0} K \hra{} \bs{T}_{j_1}^{i_1}K$ if   $i_0 \leq
i_1$ and $j_0 \geq j_1$. That is, $\bs{S}_{\bullet}^{\bullet} K$ and
$\bs{T}_{\bullet}^{\bullet} K$ are double filtrations of
$\bs{S}_{0}^{\infty} K$ and $\bs{T}_{0}^{\infty} K$,
respectively.
\subsection{}\label{861}
We need a special notation: denote by $\wt{m}$ the number
$3(\sum_{j=1}^{j=m} l!-m)$. It satisfies the recurrence
$\wt{m}=\wt{m-1}+3(m!-1)$. The composite
$$\bs{S}_{i}^m K \xar{\mb{\beta}} \bs{T}_{i}^m K
 \hra{} \bs{T}_{i - \wt{m} -1}^{m!}
K $$ will be denoted by $\bs{\alpha}^m_i$.

\begin{prop}[On functorial separation of small balls]  \label{fuballsep}
For any
geometric simplicial manifold  $K$, any positive integer  $m$,  and any
$i\geq \wt m +2$, there exists a natural transformation $\theta^m_i$
of the poset morphism $\bs{S}_i^m K
\xar{\bs{\alpha}^m_i}\bs{T}_{i-\wt m -1}^{m!} K$ into a poset
morphism $\bs{S}_{i}^m K \xar{\bs{\gamma}^m_i}\bs{D}_{i-\wt
m-1}^{m!} K
 \hra{} \bs{T}_{i-\wt m -1}^{m!} K
$.
\end{prop}
\begin{proof}
Induction on $m$.

\smallskip
1. Let $m=1$. Then $m!=1$, $\wt m = 0$, and the poset structure on
$\bs{S}_{i}^m K$ is trivial. Set $\bs{S}_{i}^m
K\xar{\bs{\gamma}^1_i} \bs{D}_{i-1}^{1} K $. According to
Proposition~\ref{prop49},  for every ball $B \in \bs{S}_{i}^m$
 we can choose a larger $K$-consistent ball $\bs{\gamma}^1_i(B)$
 from $\bs{D}_{i-1}^{1} K$. We have
$B \subseteq \bs{\gamma}^1_i(B)$, i.e., $\bs{\alpha}^1_i (B)
\preccurlyeq \bs{\gamma}^1_i(B)$, and
$\theta^1_i$ is defined.

\smallskip
2. Assume that  $\theta_i^{m-1}, \bs{\gamma}_i^{m-1}$ are already
constructed. Note that $\bs{S}_j^i K, \bs{T}_j^i K$ are
truncated upper
 semilatices
(the supremum exists provided that there exists an upper bound). The semilatice
$\bs{S}_j^i K$ is free. Consider the diagram
\begin{equation}\label{balls3}
\bfig
\square/..>`<-_{)}`<-_{)}`>/<800,400>[\bs{S}_{i}^{m}K`\bs{T}_{i-\wt{m} -1}^{m!}K
`\bs{S}_{i}^{m-1}K`\bs{D}_{i-\wt{m-1} +1}^{(m-1)!}K.;\bs{\delta}_m```\bs{\gamma}_i^{m-1}]
\efig
\end{equation}
Since $\bs{S}_{i}^{m}K$ is free, the extension
 $\bs{\delta_m}$
is canonically defined on the set $\{A_i\}_{i=1}^m \in
\bs{S}_{i}^{m} K$ by the rule
$$ \bs \delta_m(\{A_i\}_{i=1}^m)=
\vee_{i=1}^m \bs{\gamma}_i^{m-1}(\{A_j\}_{j \in \{1,\ldots,m\}\setminus
\{i\}}).$$ Therefore
 $ \bs \delta_m(\{A_i\}_{i=1}^m)$ has at most
$m!$ balls. Since $\bs{S}_{i}^{m}$ is free, it follows that if
we choose, according to Proposition~\vref{balls2}, an arbitrary separating
cover  $ \bs{\delta}_m(\{A_i\}_{i=1}^m)$, then we will obtain a natural
transformation $\bs{\delta}_m \to^{\vartheta_m}_\bullet
\bs{\gamma}_m$ of the functor $\bs{\delta}_m$ into some functor
$$\bs{S}_{i}^{m}K\xar{\bs{\gamma}_m}\bs{D}^{m!}_{i - \wt m -1} K.$$
Now we can put $\bs \theta_m=\vartheta_m\psi_m$, where $\psi_m$ is
canonically defined by the commutativity of the square of natural
transformations
$$
\bfig
\square/..>`<-_{)}`<-_{)}`>/[\bs{\alpha}_i^m`\bs \delta_m`\bs{\alpha}_i^{m-1}`\bs{\gamma}_i^{m-1}.;%
\psi_m``(\ref{balls3})`\bs \theta{m-1}]
\efig
$$
Here the left arrow
is the canonical embedding as a subfunctor, and
the right arrow
is the embedding defined by diagram~(\ref{balls3}).
\end{proof}

\section{Lemma on prismatic fragmentation
of fiberwise homeomorphisms over the cube} \label{fragm} Here we
assemble the knowledge collected in Secs.~\ref{some_scheme}--\ref{balls_sep}
into a lemma which is a multidimensional generalization of
Hudson's isotopy fragmentation.
\begin{lemma}[On prismatic fragmentation
of fiberwise homeomorphisms over the cube]\label{lem_frag_on_qube}
Assume that the following data are fixed: $\ms Q \in \mb R (X)$, a
set of homeomorphisms $G_1,{\ldots} ,G_m,U \in \PL_{I^\mb{k}}(X)$, a
triangulation $\ms K \trianglelefteq \ms Q$, a positive integer $M \geq 1$, an
Alexandroff presheaf
$$(\partial \frac{1}{M} [I^k])^\op \xar{\xi}\mb D^\infty_c(K).$$
Assume that the homeomorphisms $G_1,{\ldots} , G_m$ are $\delta_{\ms
Q}$-prismatic and the homeo\-morphisms $G_1\lfloor_{\partial
I^k},{\ldots} ,G_m\lfloor_{\partial I^k}, U\lfloor_{\partial I^k}$ are
$\xi$-prismatic.

Then there exist
\begin{itemize}
\item homeomorphisms
$\wt G_1,{\ldots} ,\wt G_m, \wt U \in \PL_{I^k \times I}(X)$ such that
$\wt G_i\lfloor_{I^k \times \{0\}}=G_i$, $\wt U\lfloor_{I^k \times \{0\}}=U$,
\item an $\mr{OPL}$ Alexandroff topology $\mc T$ on the cubic bucket
$\Xi=\partial I^k \times I \cup I^k \times \{1\}$ such that
$\mc T|_{\partial I^k \times \{0\}}= (\partial \frac{1}{M}[ I^k])^\op$,
\item an Alexandroff presheaf $\mr (\Xi, \mc T)\xar{\wt \xi}
 \mb D^\infty_c(K)$
extending $\xi$
such that $\wt G_1,{\ldots} ,\wt G_m$ are $\delta_{\ms
Q}$-prismatic and $\wt G_1\lfloor_{\Xi},{\ldots} ,\wt G_m\lfloor_{\Xi},\wt
U\lfloor_{\Xi}$ are $\wt \xi$-prismatic. If the homeomorphism $U$
is $\delta_{\ms Q}$-prismatic, then the homeomorphism $\wt{U}$ can
be chosen to be $\delta_{\ms Q}$-prismatic.
\end{itemize}
\end{lemma}
\begin{proof}
0. In the proof we construct the following:
\begin{itemize}
\item six  data sets $\la G^j_1,{\ldots} ,G^j_m,U^j,\mc T^j, \xi^j\ra$, where
$G^j_i,U^j\in PL_{I^k}$, $j=0,{\ldots} ,5$, $\mc T^j$ is an Alexandroff
topology on
 $\partial I^k$, $(\partial I^k,\mc T^j)\xar{\xi^j}\mb D^\infty_c(K)$;
\item five data sets $\la G^{j,j+1}_1,{\ldots} ,G^{j,j+1}_1,U^{j,j+1}, \mc
T^{j,j+1},\xi^{j,j+1} \ra$, where $G_i^{j,j+1},U^{j,j+1} \in
\PL_{I^k\times I}$, $ j=0,{\ldots} ,4$,
$(\partial I^k\times I,\mc T^{j,j+1})\xar{\xi^{j,j+1}}\mb D^\infty_c(K))$;
\item a topology $\mc E$ on $I^k$ and an Alexandroff presheaf $(I^k, \mc
E)\xar{\psi} \mb D^\infty_c(K)$
such that
\begin{itemize}
\item $G_i^{j,j+1}\lfloor_{I^k\times 0}=G^j $,
\item $G_i^{j,j+1}\lfloor_{I\times 1}=G^{j+1}$,
\item $U^{j,j+1}\lfloor_{I\times \{0\}}=U^j$,
\item $U_i^{j,j+1}\lfloor_{I\times \{1\}}=U^{j+1}$,
\item $h^0 \xi^{j,j+1}=\xi^i$, $h^1 \xi^{j,j+1}=\xi^{j+1}$,
\item $G^0_i=G_i$, $U^0_i=U$, $\xi^0=\xi$,
$\partial \psi = \xi^5$,
\item $G^j_i\lfloor_{\partial I^k }, U^j\lfloor_{\partial I^k }$ are
$\xi^j$-prismatic,
\item $G_i^5, U^5$ are $\psi$-prismatic, $G^{j,j+1}_i\lfloor_{\partial I^k
}, U^{j,j+1}\lfloor_{\partial I^k }$ are
$\xi^{j,j+1}$-prismatic,
\item all  $G^j_i, G_i^{j,j+1}$ are $Q$-prismatic,
\item  $U^j,U^{j,j+1}$ are $Q$-prismatic if $U$ is $Q$-prismatic.
\end{itemize}
\end{itemize}

Having all these data,  we can past them, according to Proposition~\vref{p17},
into homeomorphisms on $\PL_{I^k\times[0,5]}$ using
the scheme presented at the figure below.
$$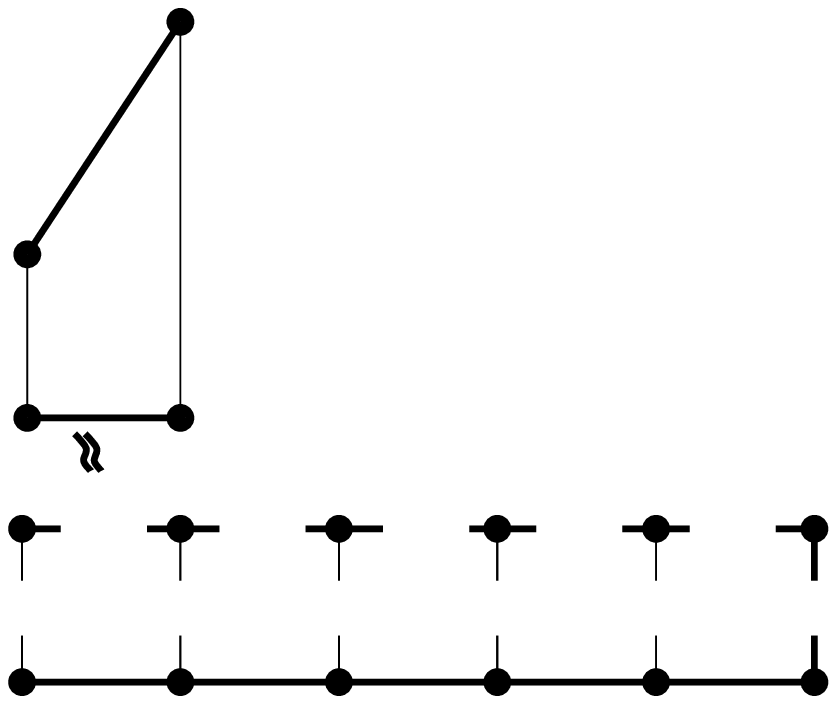$$
Then, making a linear change of coordinates in the base, we will obtain the
required homeomorphisms from $\PL_{I^k\times I}$.

\smallskip
1. Using Proposition~\vref{prop55}, build some Alexandroff
presheaf
$$(\partial I^k, \mc A (\frac{1}{M} \partial I^k))\xar{\eta}\mb D^\infty_c(K)$$
such that $\eta > \xi$.

\smallskip
2. Choose a simplicial subdivision $K' \trianglelefteq K$ that simultaneously
subdivides all  balls in all configurations $\eta(y)$, $y\in\partial
I^\mb{k}$, and in all configurations
 $\xi(y)$, $y\in\partial I^\mb{k}$. There are finitely many such balls,
therefore  $K'$ exists.

\smallskip
3. Using Lemma~\vref{prop48}, construct a number  $N$ and a Hudson graph
system $\mc H=\mc H(\{G_1,{\ldots} ,G_m,U\},\sd_{\wt k} K',M) \in
\mb{HGF}(X,k,N,M)$,
where $\sd_{\wt k} $ is defined in Sec.~\vref{861}.

\smallskip
4. Construct $G^{0,1}$, $U^{0,1}$, $\mc T^{0,1}$, $\xi^{0,1}$.

We apply Lemma~\vref{lem8}. Consider the following  fiberwise homeomorphisms of
the trivial bundle $X \times T^k \xar{\pi_2}T^k$: $\ol
G^{0,1}_i=F_\mc{H}\rtimes G_i$, $i=1,{\ldots} ,m$, and
 $\ol U^{0,1}=F_\mc{H}\rtimes U$.
On $T^\mb{k}$ we have the
Alexandroff topology $\Theta(M)$.
Proposition~\vref{45}
defines a $\PL$ homeomorphism
$\varphi$, which is also a homeomorphism of Alexandroff spaces:
$$
\bfig \Vtriangle[(T^\mb{k},\Theta)`(\frac{1}{M}
{[I^k]})^\op\times I^\mr{triv}` {I^\mr{triv};};\varphi``] \efig $$
besides, $\varphi$ has the following properties:
$\varphi|_{H^k_0}=\id$, $\varphi({H^k_1})=I^k\times\{1\}$,
 $\varphi(W^k)=(\partial I^k) \times I$.
Put $G^{0,1}_i=\ol G^{0,1}_i\circ \varphi^{-1}$, $i=1,{\ldots} ,m$, and
$U^{0,1}=\ol U^{0,1}\circ \varphi^{-1}$.

According to assertion~(b) of Lemma~\ref{lem8}, the
homeomorphisms
 $G_i^{0,1}$
are $\ms Q$-prismatic and the homeomorphism $U^{0,1}$ is $\ms Q$-prismatic if
$U$ is  $\ms Q$-prismatic.

Define $\mc T^{0,1}$ as $\varphi^{-1} \Theta = \mc
A^c(\partial\frac{1}{N}[I^k] \times I^{\mr{triv}})$. Define
$\ol\xi^{0,1}$ according to (\ref{eq78}) and put
$\xi^{0,1}=\ol\xi^{0,1}\circ\varphi^{-1}: (\partial I^k \times I,
\mc T^{0,1})\xar{}\mb D^\infty_c(K)$. By assertion~(c) of
Lemma~\vref{lem8}, the homeomorphisms
 $G_i^{0,1}\lfloor_{(\partial I^k) \times I}$ and
 $U^{0,1}\lfloor_{(\partial I^k) \times I}$
are $\xi^{0,1}$-prismatic.

 \smallskip
5. Define $G_i^j$, $U^j$, $G^{j,j+1}_i$, $U^{j,j+1}$ for $j\geq 1$.

Put $G_i^1=G_i^{0,1}\lfloor_{I^k\times \{1\}}$ and $U^1=
U^{0,1}\lfloor_{I^k\times \{1\}}$. For $j \geq 1$, put $G_i^j \equiv
G^1_i$ and $G^{j,j+1}_i \equiv G^1_i \times \id$. Put $U^j \equiv
U^1$ and $U^{j,j+1} \equiv U^1 \times \id$.

\smallskip
6. Define $\mc E$ and $\psi$.

\smallskip
6.1. According to assertion~(a) of Lemma \ref{lem8}, the
homeomorphisms  $G^{1}_i$ are $\mb J'(=\mb J_\mc{H}\circ
\varphi^{-1})$-prismatic, where $(I^k, \mc E =\varphi^{-1} \mc
E(N))\xar{\mb J'}\mb{Al}(X)$ and the topology $ \mc E|_{\partial
I^k}$ is inscribed into the topology $\mc A(\frac{1}{M}\partial
[I^k]) $ on $I^k$.

\smallskip
6.2. We can apply  Proposition~\vref{fuballsep} to $\mb  J'$
and obtain an Alexandroff presheaf $(I^k, \mc E
)\xar{\psi} \mb D^{k!}_c(\sd_1 K') $ such that $G^1_i$ and $U^1$
are $\psi$-prismatic.

\smallskip
7. Define  $\mc T^{1,2}$, $\xi^{1,2}$, $\mc T^2$, $\xi^2$.

Note that the topology $\mc A(\frac{1}{M}\partial [I^k])$ is
dense (see Sec.~\vref{dense}), therefore, by  Proposition~\vref{pro15}, there
exists a morphism $\mf B \partial (I^k, \mc E)\xar{\phi}\mf B
(\frac{1}{M}\partial [I^k])^\op $ such that for every $E \in \mf B
(I^k, \mc E)$ the inclusion $E \subset \phi(E)$ holds. Now we can apply
Proposition~\vref{p22}
to $\xi^1=\xi$ and $\phi$. Put $\mc T^2=\partial \mc E $,
 $\xi^2=\xi^1\circ \phi\circ o :(\partial I^k,  \mc E^2) \xar{}
\mb D^\infty_c (K)$. Put $\mc T^{1,2}=\wt{\mr{Cyl}}_\phi$ and
$\xi^{1,2}=\wt{\mr{Cyl}}_\phi(\xi^1,\xi^2)$. Proposition~\vref{p22}
states that $G_i^{1,2},U^{1,2}$ are prismatic with respect to
$\xi^{1,2}$.

Put $\mc T^5=\mc T^4 = \mc T^3 =\mc T^2=\partial \mc E $.

\smallskip
8. Define $\xi^{2,3}$.

Put $\xi^3 = \eta \circ \phi\circ o :\partial (I^k, \mc E) \xar{}
\mb D^\infty_c (K)$. We have $\eta > \xi$, therefore $\xi^2 <
\xi^3$. Apply Proposition~\vref{p21} for $\alpha = \xi^2$,
$\beta=\xi_3$, $w=\id$. Define $\mc T^{2,3}=\underline{\mr{Cyl}}_\id$,
$\xi^{2,3}=\underline{\mr{Cyl}}_\id(\xi^3,\xi^2)$. By Proposition~\vref{p21},
the homeomorphisms  $G_i^{2,3},U^{2,3}$ are prismatic
with respect to $\xi^{2,3}$.

\smallskip
9. Define $\xi^{3,4}$, $\xi^{4,5}$.

Build a morphism $\partial (I^k, \mc E)
\xar{\xi^4}\mb D^{k!}_c \sd_1 K'$. The morphism $\xi^4$ is defined
as follows: $\xi^4(b)=$\{all the balls from  $\psi(b)$ having a
common point with a ball from $\xi^2(b)$\}. The homeomorphisms $G_i^3,
U^3$ are $\xi^2$-prismatic and $\psi$-prismatic, therefore they are
$\xi^4$-prismatic. By construction,
 $\psi|_{\partial I^k}=\xi^5>\xi^4 < \xi^3$.
Applying Proposition~\vref{p20} to $\xi^3,\xi^4, \id$, we obtain
$\xi^{3,4}=\ol{\mr{Cyl}}_\id (\xi^3,\xi^4)$. Applying Proposition~\vref{p21}
to $\xi^4,\xi^5, \id$, we obtain
$\xi^{4,5}=\underline{\mr{Cyl}}_\id (\xi^4,\xi^5)$.
\end{proof}

\section{The proof of the lemma on a common $\mb R(X)$-triangulation
of fiberwise homeomorphisms} \label{proof} Here we will prove
Lemmas~\ref{loc_lem} and~\ref{lem13}. The latter is the $\pmb
\frown$-lemma for the pair $(\underline{\mr{Prism}}^{m}(X)\times
\PL(X), \underline{\mr{Prism}}^{m+1}(X))$. This will complete the
proof of Theorem~A.
\subsection{Local lemma}
Here we prove that some data related to the prismaticity of a family
of fiberwise homeomorphisms over a simplicial bucket can be extended
to similar data over the
filling of the bucket.

Fix the following data:
\begin{itemize}
\item{} $Q \in \mc N_k \mb R (X)$,
$Q=(\ms Q_0\trianglelefteq \ms Q_1\trianglelefteq...\trianglelefteq \ms Q_k)$,
\item{} the standard  $Q$-coloring of the
ordered simplicial complex $[\Delta^k]$ and the
Alexandroff presheaf $(\Delta^k, \mc A^c([\Delta^k]))
\xar{\mr{Max} Q} \mb R(X)$  generated by the
coloring (see~(\ref{e33})),
\item{} the prism $\Delta^k \times I$ with the structure of the ball
complex $[\Delta^k]\times [0,1]$,
\item{} the Alexandroff
presheaf
$$(\Delta^k \times I, \mc A^c([\Delta^k]\times [0,1]))\xar{\alpha (Q)=
\mr{Max} Q\}\lfloor_\pi} \mb R (X)$$
 induced by the projection
$[\Delta^k]\times [0,1] \xar{\pi} [\Delta^k]$,
\item{} the ``lower simplicial bucket''
$$\underline\Gamma^k= \Delta^k\times\{0\}\cup \partial \Delta^k \times I
\subset \Delta^k \times I $$
with the ball structure $[\underline\Gamma^k]$ induced by
the embedding,
\item{} an $\OPL$ Alexandroff topology
$\mc T$ on the rim  $$\partial \Delta^k \hra{\id \times \{1\}}
 \underline \Gamma^k$$ of the bucket,
\item{} a common triangulation $K\trianglelefteq \ms Q_i$, $i=0,{\ldots} ,k$,
\item{} an Alexandroff presheaf
$(\partial \Delta^k\times \{1\}, \mc T)\xar{\xi} \mb
D^\infty_c(K)$,
\item{} a set of fiberwise homeomorphisms over the bucket:
$$G_1,{\ldots} ,G_m, U \in \PL_{\underline\Gamma^k}(X).$$
\end{itemize}
\begin{lemma} \label{loc_lem}
Assume that our data satisfy the following conditions:
\begin{enumerate}
\item{} the homeomorphisms
$G_1,{\ldots} , G_m$ are
$\alpha(Q)\lfloor_{\underline\Gamma^k}$-prismatic,
\item{} the homeomorphisms $G_1\lfloor_{\partial \Delta^k \times \{1\}},{\ldots} ,
G_m\lfloor_{\partial \Delta^k \times \{1\}}, U\lfloor_{\partial
\Delta^k \times \{1\}}$ are $\xi$-prismatic,
\item{} if the homeomorphism $U\lfloor_{d_i \Delta^k \times \{0\}}$
is $\alpha(Q)\lfloor_{d_i \Delta^k \times \{0\}}$-prismatic, then
its extension $U\lfloor_{d_i \Delta^k
\times I}$ to the wall of the bucket is $\alpha(Q)\lfloor_{d_i \Delta^k \times
I}$-prismatic.
\end{enumerate}
Then there exist
\begin{itemize}
\item
an $\OPL$ Alexandroff topology $\wt{\mc T}$ on  $\Delta^k \times
 \{1\}$,
\item
an Alexandroff presheaf $(\Delta^k\times \{1\}, \wt{\mc
T})\xar{\wt{\xi}} \mb D_c^\infty(K)$,
\item a set of fiberwise homeomorphisms
 $\wt{G}_1,{\ldots} ,\wt{G}_m,\wt{U} \in \PL_{\Delta^k \times I}(X)$
\end{itemize}
such that
\begin{enumerate}
\item $\wt{\mc T}\lfloor_{\partial \Delta^k \times \{1\}}=\mc T$,
 $\wt{\xi}\lfloor_{\partial \Delta^k \times \{1\}}=\xi$,
\item
$\wt{G}_i\lfloor_{\underline{\Gamma}^k}=G_i$, $i=1,{\ldots} ,m$,
$\wt{U}\lfloor_{\underline{\Gamma}^k}=U$,
\item
$\wt{G}_i\lfloor_{\Delta^k \times\{1\}}$, $i=1,{\ldots} ,m$, and
 $\wt{U}\lfloor_{\Delta^k\times\{1\}}$ are
$\wt{\xi}$-prismatic,
\item
$\wt{G}_i$, $i=1,{\ldots} ,m$, are $\alpha(Q)$-prismatic; if
$U\lfloor_{\Delta^k \times \{0\}}$ is $Q$-prismatic, then one can
choose $\wt{U}$ to be $\alpha(Q)$-prismatic.
\end{enumerate}

\end{lemma}
\begin{figure}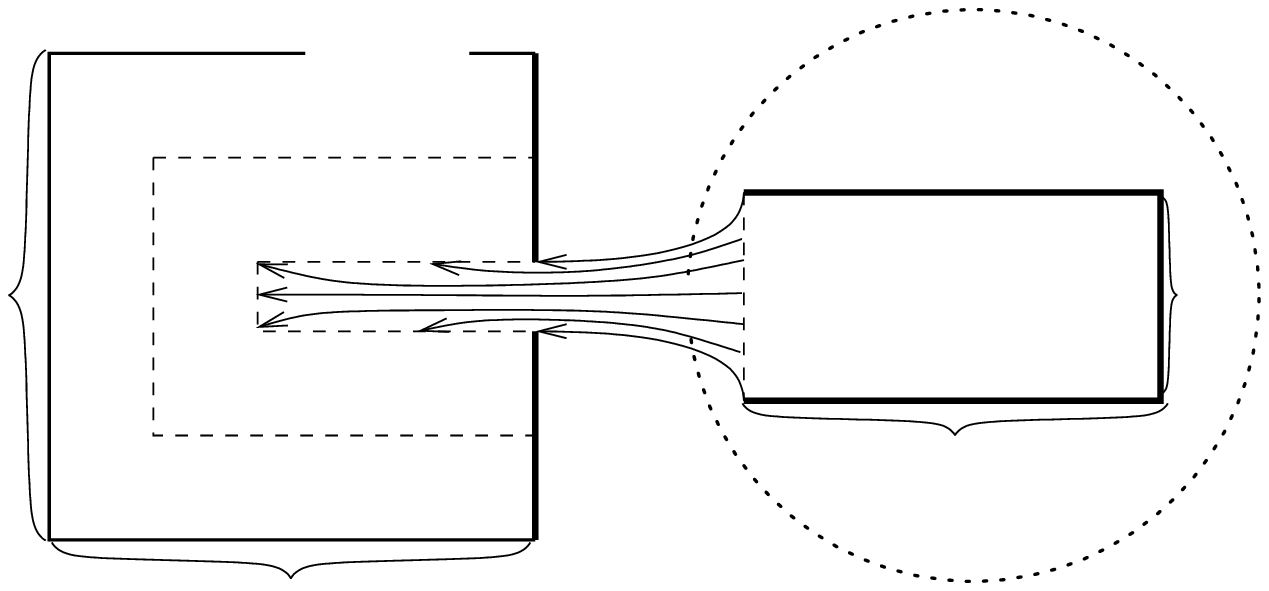\caption{\label{lem13fig}}\end{figure}
\begin{proof} We construct a filling
 $\wt{\mc T}$, $\wt{\xi}$, $\wt{G_i}$, $i=1,{\ldots} ,m$, $\wt{U}$
 of the bucket from three pieces (see Fig.~\ref{lem13fig}). The
last, principal, piece will be obtained from Lemma~\ref{lem_frag_on_qube}.
The two first pieces are needed to adjust our
data to take the form required for applying Lemma~\ref{lem_frag_on_qube}.

\smallskip
1. We need the following data:
\begin{itemize}
\item a homeomorphism $I^k \xar{\varphi} \Delta^k$ of the cube onto the simplex,
\item $G_i^{j,j+1}\in \PL_{\Gamma^k \times I}(X)$, $i=1,{\ldots} ,m$, $j=0,1$,
$U^j\in \PL_{\Delta^k \times I}$, $j=0,1$,
\item
 $\wh{G}_i\in \PL_{I^k \times I} (X)$, $i=1,{\ldots} ,m$,
 $\wh{U} \in \PL_{I^k \times I} (X)$,
\item Alexandroff topologies
 $\mc T^{j,j+1}$ on $\partial \Delta^k \times I$ and an Alexandroff
topology $\wh{\mc T} $ on the upper cubic bucket $\ol{\Xi}^k =
\partial I^k \times I \cup I^k \times \{1\}$,
\item morphisms
$(\partial \Delta^k \times I, \mc T^{j,j+1} \xar{\xi^{j,j+1}}
\mb D^\infty_c (K))$  for $j=0,1$ and a morphism
 $(\ol{\Xi}^k, \wh{\mc T}) \xar{\wh{\xi}} \mb D^\infty_c (K)$.
\end{itemize}
Assume that the following conditions are satisfied:
\begin{itemize}
\item[(a)] $G^{0,1}_i\lfloor_{\ol{\Gamma}^k\times \{0\}} = G_i$, $i=1,{\ldots} ,m$;
$U^{0,1}\lfloor_{\ol{\Gamma}^k\times \{0\}} = U$;
$\xi^{0,1}\lfloor{\partial \Delta \times \{0\}}=\xi$;
\item[(a1)] $G^{0,1}_i$, $i=1,{\ldots} ,m$, are $\delta_{\ms{Q}_k}$-prismatic;
$U^{0,1}$ is $\delta_{\ms{Q}_k}$-prismatic if $U$ is
$\alpha(Q)$-prismatic;
\item[(a2)] $G^{0,1}_i\lfloor_{\partial \Delta^k \times \{1\} \times I  }$
and $U^{0,1}_i\lfloor_{\partial \Delta^k \times \{1\} \times I
}$ are
 $\xi^{0,1}$-prismatic;
\item[(b)]  $G^{1,2}_i\lfloor_{\ol{\Gamma}^k\times \{0\}}
= G^{0,1}_i\lfloor_{\ol{\Gamma}^k\times \{1\}}=G_i^1$,
$i=1,{\ldots} ,m$; $U^{1,2}\lfloor_{\ol{\Gamma}^k\times \{0\}} =
U^{0,1}\lfloor_{\ol{\Gamma}^k\times \{1\}}=U^1$;
$\xi^{1,2}\lfloor{\partial \Delta \times \{0\}}
=\xi^{0,1}\lfloor{\partial \Delta \times \{1\}}$;
\item[(b1)] $G^{1,2}_i$, $i=1,{\ldots} ,m$, are $\delta_{\ms{Q}_k}$-prismatic;
 $U^{1,2}$ is $\delta_{\ms{Q}_k}$-prismatic if
$U^{0,1}\lfloor_{\Gamma^k \times \{1\}}$ is
 $\delta_{\ms{Q}_k}$-prismatic;
\item[(b2)] $G^{1,2}_i\lfloor_{\partial \Delta^k \times \{1\} \times I  }$
and $U^{1,2}_i\lfloor_{\partial \Delta^k \times \{1\} \times I
}$ are
 $\xi^{1,2}$-prismatic;
\item[(c)] $\wh{G}_i\lfloor_{I^k \times \{0\}}=
(G^{1,2}\lfloor_{\Gamma^k \times \{1\}})\lfloor_{\varphi}=G^2_i$,
$i=1,{\ldots} ,m$; $\wh{U}\lfloor_{I^k\times \{0\}}=
(U^{1,2}\lfloor_{\ol{\Gamma}^k \times
\{1\}})\lfloor_{\varphi}=U^2$; $\wh{\xi}\lfloor_{I^k\times
\{0\}}= (\xi^{1,2}\lfloor_{\Gamma^k \times
\{1\}})\lfloor_{\varphi}$;
\item[(c1)] $\wh{G}_i$, $i=1,{\ldots} ,k$, are $\delta_{\ms{Q}_k}$-prismatic;
 $\wh{U}$ is $\delta_{\ms{Q}_k}$-prismatic if
$\wh{U}\lfloor_{I^k\times \{0\}}$ is
$\delta_{\ms{Q}_k}$-prismatic;
\item[(c2)] $\wh{G}_i\lfloor_{\ol{\Xi}^k}$, $i=1,{\ldots} ,k$, and
$\wh{U}\lfloor_{\ol{\Xi}^k}$ are $\wh{\xi}$-prismatic.
\end{itemize}
If we have data satisfying all these conditions, then we
construct a $(k+1)$-dimensional $\PL$ ball $P$ as the colimit of the
following diagram of solid arrows:
$$
\xymatrix{& & P & &  \\
          \Gamma^k\times I \ar@{-->}[urr] & & \Gamma^k\times I \ar@{-->}[u]&  & I^{k}\times I. \ar@{-->}[ull] \\
         & \Gamma^k \ar[ur]^{\id \times \{0\}} \ar[ul]^{\id\times \{1\}} &
          & I^k \ar[ul]_{\varphi\times \{1\}}\ar[ur]_{\id \times \{0\}} }
$$
The boundary of $P$ is assembled from two
 $k$-dimensional balls $\Gamma^k$ and
$S$:
\begin{equation}\label{ee73}
\xymatrix{& & S & &  \\
          \partial \Delta^k \times I \ar@{-->}[urr] &
          & \partial \Delta^k \times I \ar@{-->}[u]&  & \ol{\Xi}^k. \ar@{-->}[ull] \\
         & \partial \Delta^k \ar[ur]^{\id \times \{0\}} \ar[ul]^{\id\times \{1\}}
         & & \partial I^k \ar[ul]_{\varphi\times \{1\}}\ar[ur]_{\id \times \{0\}} }
\end{equation}
We define $\wt{G}'_i, \wt{U}' \in \PL_P(X)$ as pastings (see Sec.~\vref{ss611}):
$$
\xymatrix{& & \wt{G}'_i, \wt{U}' & &  \\
G_i^{0,1}, U^{0,1} \ar@{-->}[urr] & & G_i^{1,2}, U^{1,2} \ar@{-->}[u]&
 & \wh{G}_i, \wh{U}. \ar@{-->}[ull] \\
         & G^1_i, U^1 \ar[ur]^{\lfloor_{\id \times \{0\}}} \ar[ul]^{\lfloor_{\id\times \{1\}}} &
          & G^2_i,U^2 \ar[ul]_{\lfloor_{\varphi\times \{1\}}}\ar[ur]_{\lfloor_{\id \times \{0\}}}
          }
$$
We define a topology $\wt{\mc T}'$ on the ball $S$ and an Alexandroff
presheaf
$(S,\wt{\mc T}' )\xar{\wt{\xi}'}\mb{D}^\infty_c(K)$ as the natural
pasting of the topologies and presheaves $\mc T^{0,1}, \xi^{0,1}, \mc
T^{1,2},\xi^{1,2}$ and $\wh{\mc T}, \wh{\xi}$ according to diagram~(\ref{ee73}).
Pick a homeomorphism $ \Delta^k \times I \xar{\psi} P
$ that is identical on $\Gamma^k$ and sends
 $\Delta^k \times \{1\}$ to $S$.

Now, if we put $\wt{G}_i=\wt{G}'_i\lfloor_\psi$,
$\wt{U}=\wt{U}'\lfloor_\psi$, $\wt{\mc T}=\wt{\mc T}'\lfloor_{\psi}$,
and $\wt{\xi}=\wt{\xi}'\lfloor_{\psi}$, then we will
obtain data that satisfy
 conditions~(1)--(4) due to  conditions~(a)--(c2) on the ingredients of our pasting.

\smallskip
2. Now we will present the ingredients required in Step~1.

\smallskip
2.A. Put
 $G^{0,1}_i = G_i \times \id$, $U^{0,1}=U \times \id \in
 \PL(X)_{\Gamma^k \times I}$.
Build a dense Alexandroff topology $\mc T^1$ on
 $\partial \Delta^k$ that strengthens the topology $\mc T$.
Such a strengthening always exists by Proposition~\vref{p24}. Let
$(\partial \Delta^k, \mc T^1) \xar{w}(\partial \Delta^k, \mc T)$ be
the morphism of weakening the topology.
Let $I\xar{\rm{inv}}I$ be the linear
homeomorphism ${\rm inv}(t)=-t+1$, so that
 ${\rm inv}(0)=1$, ${\rm inv}(1)=0 $.
Using Proposition~\vref{p20}, build a topology $\ol{\mr{Cyl}}_w$ on
$\partial \Delta^k\times I$. Put $\mc T^{0,1}=\ol{\mr{Cyl}}_w
\lfloor_{\mr{inv}}$ (i.e., we change the orientation of the parameter
in the construction from Proposition~\ref{p20}). Put $\xi^1=\xi\circ
w$ and $\xi^{0,1}= \ol{\mr{Cyl}}_w(\xi^1\circ
w,\xi)\lfloor_{\mr{inv}}$. By  Proposition~\ref{p20},
conditions~(a), (a1), (a2) are satisfied for
$G_i^{0,1}, U^{0,1}, \xi^{0,1}$.

\smallskip
2.B. Put $G^{1,2}_i = G_i \times \id$, $U^{1,2}= U \times \id \in
\PL(X)_{\Gamma^k \times I}$. The Alexandroff topology $\mc
T^1\lfloor_\varphi$ is an  $\OPL$ topology on $\partial I^k$. By
the Lebesgue lemma (see \cite{Ke}), there is a positive integer $M$ such that the
minimal base of the topology $\mc A^c(\partial \frac{1}{M}[I^k])$ is
inscribed into the minimal base of the topology $\mc
T^1\lfloor_\varphi$. Therefore $\mf B (A^c(\partial
\frac{1}{M}[I^k])\lfloor_{\varphi^{-1}})$ is inscribed into
 $\mf B (\mc T^1)$.
Since the topology $\mc T^1$ is dense, it follows by Proposition~ \vref{pro15} that
there is a morphism $\mf B (A^c(\partial
\frac{1}{M}[I^k])\lfloor_{\varphi^{-1}}) \xar{\phi}\mf B (\mc T^1)$
such that $$\mbox{for every }U\in \mf B ( A^c(\partial \frac{1}{M}[I^k])
\lfloor_{\varphi^{-1}})\mbox{ the inclusion } U \subset \phi (U)
\mbox{ holds}.$$ Put $\mc T^2=
A^c(\partial \frac{1}{M}[I^k])\lfloor_{\varphi^{-1}}$, $\mc
\xi^\cdot=\xi'^1\circ \phi \circ o $, $\mc
T^{1,2}=\wt{\mr{Cyl}}_{\phi}\lfloor_{\mr{inv}}$, and $
\xi^{1,2}=\wt{\mr{Cyl}}_\phi(\xi^1,\xi^\cdot)\lfloor_{\mr{inv}}$
(see Proposition~\vref{p32}). Then,
 by Proposition~\ref{p32},  conditions~(b), (b1), (b2) are
satisfied  for
$ G_i^{1,2}, U^{1,2}, \xi^{1,2}$.

\smallskip\noindent
2.C. Consider the homeomorphisms induced by $\varphi$: $G^2 = G_i
\lfloor_\varphi$, $U^2 = U\lfloor_\varphi\in \PL_{I^k}$. We are in the
conditions of Lemma~\vref{lem_frag_on_qube}: the homeomorphisms
$\partial G^2_i, \partial U^2 $ are prismatic with respect to
the Alexandroff presheaf
$$\partial \frac{1}{M}[I^k]^\op
 \xar{\xi^2= \xi^{\cdot}\lfloor_\varphi} \mb D^\infty_c (K)$$
and  $(G^2)_i$ is $\delta_{\ms Q_k}$-prismatic. Therefore there
exists a topology $\wh{\mc T}$ on the upper cubic bucket $\ol{\Xi}^k$
that extends $\partial \mc A^c \frac{1}{M}[I^k] $, there exists an
Alexandroff presheaf
 $(\ol{\Xi}^k, \wh {\mc T} )\xar{\wh{\xi}} \mb D^\infty_c(K)$
that extends
 $\xi^2$, there exist
$\delta_{\mc Q_k}$-prismatic homeomorphisms $\wh G_i, \wh U$ such
that $\wh G_i\lfloor_{\ol\Xi^k}, \wh U\lfloor_{\ol\Xi^k} $ are $(\ol
\xi')^2$-prismatic. If  $U^2$ is $\delta_{\mc Q_k}$-prismatic, then
$\wh U$ is $\delta_{Q_k}$-prismatic.

Therefore for $\wh G_i, \wh U, \wh \xi$ conditions~(c), (c1), (c2)
are satisfied.
\end{proof}

\subsection{The proof of the lemma on common $\mb R(X)$-triangulations of
fiberwise homeomorphisms}

\begin{lemma}[The $\pmb \frown$-lemma for
$(\underline{\mr{Prism}}^{m}(X)\times \PL(X),
 \underline{\mr{Prism}}^{m+1}(X))$]\label{lem13}
Let $\ms B$ be  a finite simplicial ball, $B=|\ms B|$,
$\ms S=\partial \ms B$, $S=|\ms S|$.

Let $G_1,{\ldots} ,G_m, U\in \PL_Y(X)$ be fixed. Let $\mc Q$ be an $\mb
R(X)$-coloring of $\ms B$
and  $\mc Q_{\ms S}$ be the induced coloring of $\ms S$.

Assume that the coloring  $\mc Q$ is an $\mb R(X)$-triangulation of
the homeomorphisms $G_1,{\ldots} ,G_m$. Assume that the coloring
 $\mc Q_\ms{S}$ is an  $\mb R(X)$-triangulation of
$\partial U$.

Then there exist
\begin{itemize}
\item{} a triangulation  $\ms T$ of $B\times I$
such that $h^0 \ms T=\ms B$;
\item{} an $\mb R(X)$-coloring
$\wt{ \mc Q}$ of the triangulation $\ms T$ such that $h^0 \wt{
\mc Q} = \mc Q$;
\item{} homeomorphisms $\wt G_1,{\ldots} ,\wt G_m, \wt U$
such that $h^0 \wt G_i=G_i$, $h^0 \wt U=U$, $\wt{\mc Q}$
triangulates $\wt G_1,{\ldots} ,\wt G_m$, $\wt{\mc Q}\lfloor_{\ms
T_\Lambda}$ triangulates $\wt U\lfloor_\Lambda$.
\end{itemize}
Here $\Lambda=S\times I \cup B\times \{1\}$ is the upper bucket and
$\ms T_\Lambda$ is the induced triangulation of  $\Lambda$.
\end{lemma}

\begin{proof}
1. Fix a subdivision $\ms K \trianglelefteq \mc Q$. Consider the ball complex
$\Xi=\ms B \times [I]$. Define an Alexandroff presheaf $( |\Xi|, \mc
A (\Xi))\xar{\mc Q'} \mb R(X)$ by
$\mc Q'=\mr{Max}\mc Q \circ \pi_1$.

We will define
\begin{itemize}
\item homeomorphisms
 $G'_1,{\ldots} ,G'_m,U'\in \PL_{|\Xi|}$;
\item an $\OPL$-topology $\mc T$
on $B\times \{1\}$;
\item an Alexandroff presheaf
 $(B, \mc T)\xar{\xi'}\mb D_c^\infty(K)$
\end{itemize}
such that
\begin{itemize}
 \item $G'_1,{\ldots} ,G'_m,U'$ extend
 $G_1,{\ldots} ,G_m$;

 \item $G'_1,{\ldots} ,G'_m$ are $\mc Q'$-prismatic;

\item $ U'\lfloor _{S\times I}$ are
 $\mc Q'_\ms{S\times I}$-prismatic;

\item $h_1 G'_i, h_1 U'$ are $\xi$-prismatic.
\end{itemize}

\smallskip
2. Consider the following ball subcomplexes of  $\Xi$: $\underline \Xi_i=\ms B_i
\cup \ms B_{i-1}\times [I]$ and $\ol \Xi_i=\ms B_i \cup \ms B_{i}\times
[I]$, where $\ms B_i$ is the $i$-skeleton. We can define  $G'_i$ and
$\xi_i$ inductively on $\underline \Xi_i$. On $\underline \Xi_0$ the
construction is trivial. Then we can fill the buckets of $\underline
\Xi_1$ using Lemma~\ref{loc_lem} and obtain required data on
$\underline{\Xi}_1$; continuing in this way, we construct the data
announced in Step~1 of the proof.
$$\begin{picture}(0,0)%
\includegraphics{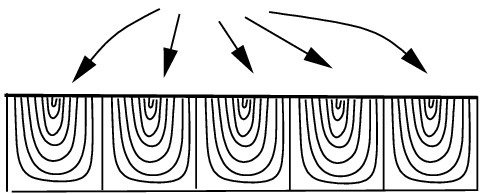}%
\end{picture}%
\setlength{\unitlength}{3947sp}%
\begingroup\makeatletter\ifx\SetFigFont\undefined%
\gdef\SetFigFont#1#2#3#4#5{%
  \reset@font\fontsize{#1}{#2pt}%
  \fontfamily{#3}\fontseries{#4}\fontshape{#5}%
  \selectfont}%
\fi\endgroup%
\begin{picture}(2314,1148)(859,-6530)
\put(1476,-5526){\makebox(0,0)[lb]{\smash{{\SetFigFont{12}{14.4}{\rmdefault}{\mddefault}{\updefault}{\color[rgb]{0,0,0}Lemma \ref{loc_lem}}%
}}}}
\end{picture}%
 $$

\smallskip
3. As a result of the constructions from Steps~1 amd~2 of the proof,
on $h^1
|\Xi|$ we have exactly the situation of Lemma~\vref{lemma7}. Therefore
there is a triangulation  $\ms V \trianglelefteq h^1 \Xi$ such that the
associated common triangulation of the homeomorphisms $\mc V$ allows
pasting the construction from Lemma~8. To complete
the proof, it remains to use Proposition~\vref{pr64} for extending the
triangulation $\mc V$ to a triangulation of $G'_1,{\ldots} ,G'_m,U'$
that does not change the triangulation $\mc Q$.
$$\begin{picture}(0,0)%
\includegraphics{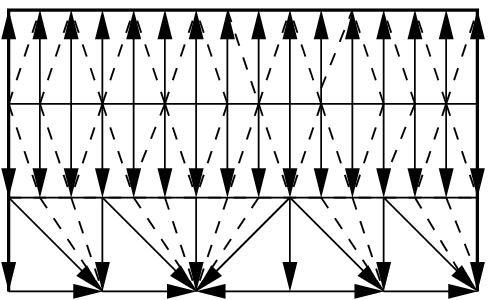}%
\end{picture}%
\setlength{\unitlength}{3947sp}%
\begingroup\makeatletter\ifx\SetFigFont\undefined%
\gdef\SetFigFont#1#2#3#4#5{%
  \reset@font\fontsize{#1}{#2pt}%
  \fontfamily{#3}\fontseries{#4}\fontshape{#5}%
  \selectfont}%
\fi\endgroup%
\begin{picture}(2334,1414)(859,-8503)
\end{picture}%
 $$
\end{proof}

\section{The case of combinatorial manifolds, Theorem~B}
The proof of Theorem~B comes from the fact that all our arguments
remain valid if we replace the poset $\mb R(X)$ by $\mb T(X)$ and
the category $\bg R(X)$ by $\bg T(X)$. Only one place in the
arguments should be specially tuned: the reference to  Lemma~\ref{lem9}
in Step~3 (p. \pageref{comb}) of the proof of Lemma~\ref{lemma7}
should be replaced by the reference to Lemma~\ref{le10}.
That is,  we should use the procedure of inscribing small
balls into combinatorial manifolds, which is slightly more complicated
than inscribing small balls into general ball complexes.

\section{Combinatorial models for $\PL_n$ fiber bundles,
the tangent bundle and Gauss map of combinatorial manifolds} Milnor
\cite{Milnor:1961} defined the group $\PL_n$ as the simplicial group
of germs at the zero section of $\PL$ homeomorphisms of $\R^n$. In this
section, we use the  Kuiper--Lashof models for
 $\PL_n$ fiber bundles
 \cite{KLI, KLII}.
Denote by $\PL((\R^n,0))$
the simplicial group of origin-preserving
$\PL$ homeomorphisms of $\R^n$. Denote by $\PL(S^n, 0,\infty)$ the
simplicial group of $\PL$-homeomorphisms of the sphere $S^n$ preserving
two different points ``$0$'' and ``$\infty$.'' Denote by
$\PL(S^n,\infty)$ the simplicial group of $\PL$-homeomorphisms of
$S^n$ preserving a single point ``$\infty$.'' It is proved in
\cite{KLII} (Lemma~1.6 at p.~248 and Lemma~1.8 at p.~249)
that in the row of natural homomorphisms of simplicial groups
\begin{equation} \label{e73} \PL((\R^n,0))\xar{\alpha}\PL_n\xleftarrow{\beta}
\PL(S^n, 0,\infty)\xar{\gamma}\PL((S^n,\infty)), \end{equation}
all homomorphisms are simplicial homotopy equivalences. In
(\ref{e73}), the homomorphism $\alpha$ sends a homeomorphism of
$(\R^n,0)$ to its germ at $0$, the homomorphism $\beta$ sends a homeomorphism of
 $(S^n,0,\infty)$ to its germ at $0$,
and  $\gamma$ is the embedding. The simplicial groups of homeomorphisms
$$\PL((\R^n,0)),\quad \PL(S^n, 0,\infty),\quad \PL((S^n,\infty))$$
are the structure groups of $\PL$ fiber bundles with marked
sections. The group $\PL( \R^n,0)$ corresponds to fiber bundles with
fiber  $\R^n$ and zero section, the group $\PL((S^n,\infty))$
corresponds to fiber bundles with fiber  $S^n$ and marked
$\infty$-section, and the group $\PL((S^n,\infty))$ corresponds to
fiber bundles with fiber   $S^n$ and two
sections $0$ and $\infty$ that have no common points.
The group  $\PL_n$ corresponds to germs
of $\R^n$-bundles near the zero section. The chain~(\ref{e73}) of simplicial
homotopy equivalences generates a chain of homotopy
equivalences of classifying spaces
\begin{equation}\label{e74}B\PL((\R^n,0))\approx B\PL_n \approx  B\PL(S^n, 0,\infty)
\approx B\PL((S^n,\infty)).\end{equation}
\begin{thm}[Kuiper--Lashof theorem on models of $\PL_n$ bundles]
There exist  functorial one-to-one correspondences between
isomorphism classes of Milnor $\PL$ microbundles and isomorphism
classes of the following $\PL$ fiber bundles:
\begin{itemize}
\item[{\rm(i)}] with fiber $\R^n$ and zero section,
\item[{\rm(ii)}] with fiber  $S^n$ and a marked section,
\item[{\rm(iii)}] with fiber $S^n$ and two marked
sections that have no common points.
\end{itemize}
\end{thm}
\subsection{Proof of Theorem C}
By the Kuiper--Lashof theory, it suffices to prove that
\begin{equation}\label{eq75}B\bg R_n \approx B\PL(S^n, \infty).\end{equation}
We will indicate how we should tweak our general constructions to obtain
the proof of (\ref{eq75}).

Consider the
$n$-dimensional $\PL$ sphere $S^n$ with a fixed point
marked by ``$\infty$.'' Consider the subposet $\mb R_n$ of the poset
$\mb R(S^n)$ (see Sec.~\vref{ss26}) formed by all ball complexes
having
 $\infty$ in the interior of a maximal ball.
This ball will be  called ``marked.'' Consider the
functor  $\mb R_n
\xar{\mb P} \bg R_n$ sending the marked ball
 to the marked combinatorial ball. Consider an $m$-simplex
$$\ms Q = \ms Q_0 \trianglelefteq \ms Q_1 \trianglelefteq{\ldots} \trianglelefteq \ms Q_n$$
of the simplicial set $\mc N \mb R_n$. Build (see Sec.~\vref{prism}) the
structure  of the
cellular bundle $\ms T(\ms Q)\xar{\ms e(\ms Q)}[\Delta^m]$
on the trivial fibration
$S^n\times \Delta^m
\xar{\pi_2} \Delta^m$ and, additionally, mark the constant section
$\infty_m = \{\infty \} \times \Delta^m \subset S^n\times \Delta^m
$. Denote by $\ms e^\infty(\ms Q)$ the pair
 $\la \ms e(\ms Q), \infty_m\ra $\label{e_infty}.
There exists one  ``marked'' prism $\ms T(\ms Q)$ (see
Sec.~\vref{prizms}) that contains the section $\infty_m$. A
$\ms Q$-prismatic homeomorphism $G\in \PL(S^n, \infty)$ (see Sec.~\vref{ss322})
is naturally defined. By construction, a
$\ms Q$-prismatic
homeomorphism preserves $\infty_m$ and sends the marked prism to the
marked prism.
We obtain the groupoid $\mr{Prism}(S^n,\infty)$ of prismatic
homeomorphisms, which is a  subgroupoid of $\mr{Prism}(S^n)$
(see Sec.~\ref{ss323}). For $\mr{Prism}(S^n,\infty)$-homeomorphisms, Lemma~\ref{Alex}
is valid, with an additional observation that we can force the
parametric Alexander trick to respect one fixed section (this is
equivalent to the fact that $\PL(D^{n},0,\partial) $ is contractible, see
\cite[Lemma~1.4, p.~248]{KLII}). Such a version of Lemma~\ref{Alex}
is sufficient for the validity of all arguments concerning the  $\mc
W$-construction  for $\mr{Prism}(S^n,\infty)$ and $\PL(S^n,\infty)$
from Secs.~\ref{wprism}--\ref{Aproof} and Sec.~\ref{cm}. Therefore the
proof of (\ref{eq75}) is reduced to the following wording
of the geometric Lemma~\vref{CombFrag1}:

\medskip\noindent
{\em the pair
$$|(\underline{\mr{Prism}}^{N-1}(S^n,\infty)\times \PL(S^n, \infty)|,
|\underline{\mr{Prism}}^N(S^n,\infty))|$$ is homotopy trivial.}

\medskip
This lemma is automatically proved by all our constructions of
fragmentation and surgery, since the scheme of fragmentation of
fiberwise homeomorphisms from Secs.~8--14 and all surgery of
generalized prismatic homeomorphisms from Sec.~6 respect all fixed
sections. We should make two remarks on the choice of
triangulations. In Lemma~\vref{lem_frag_on_qube} we should make sure
that the triangulation $\ms K \trianglelefteq \ms Q\in \mb R(S^n,\infty)$
contains the point $\infty$ as a vertex. Then the configurations
of small closed balls from Sec.~13 will involve either  balls containing
 the point $\infty$ in the interior or  balls that do not
touch this point at all. Therefore the construction of Lemma~\vref{lem9}
will not lead us out of $\mb R_n$ provided that we take care
of the following:
 all common triangulations that appear in the proof of
Lemma~\ref{lem9} should be performed inside $\mb R_n$, which is always possible.

\subsection{Theorem D}
\subsubsection{Milnor's tangent microbundle}
Traditionally, the tangent bundle of a $\PL$ manifold $M^n$ is
defined as the tangent microbundle \cite{Milnor:1961}. The tangent
microbundle is constructed as follows.

Consider the square  $M^n\times M^n$ of the manifold, the projection
 $M^n\times M^n \xar{t} M^n$
to the first argument, and the section $M\xar{0_t} M^n\times M^n$:
$0_t(x)=(x,x)$. The \bem{tangent microbundle of $M^n$} is the germ of
$t$ at the section $0_t$. According to the Kuiper--Lashof theory, there
is a unique, up to isomorphism, $(S^n, 0, \infty)$-bundle
$t_M^{0,\infty}$ on $M^n$ such that its germ at the zero section is
isomorphic to $t$. The fiber bundle $t_M^{0,\infty}$ will be called the
tangent $(S^n, 0, \infty)$-bundle of $M^n$. For the tangent $(S^n, 0,
\infty)$-bundle, there is a unique, up to isomorphism, $(S^n,
\infty)$-bundle $t^\infty_M$ (just forget the $0$-section). This
fiber bundle will be called the tangent $(S^n, \infty)$-bundle of
$M^n$. These correspondences are one-to-one correspondences of
isomorphism classes. This means that the Gauss maps
\begin{equation}M^n\xar{G} B\PL_n, \quad M^n\xar{G^{0,\infty}} B\PL(S^n,0,\infty),\quad
M^n \xar{G^\infty} B\PL(S^n,\infty) \end{equation}
of the tangent bundles $t_M, t_M^{0,\infty},t_M^{\infty}$
coincide up to homotopy after the identification~(\ref{e74}).

\subsubsection{The proof of Theorem D}
According to  (\ref{eq75}), $B\bg R_n \approx B\PL(S^n
\infty)$. Therefore for any locally ordered simplicial  complex $K$ with an
$\bg R_n$-coloring $\bs i K \xar{\mc Q} \mb d \mc N \bg R_n$ (see
Sec.~\vref{color}), $|K|\xar{|\mc Q|} B\bg R_n$ is the
Gauss map for the  $(S^n,\infty)$-bundle $\ms e^\infty (\mc Q)$,
where $\ms e^\infty (\mc Q)$ is pasted (Sec.~\vref{ss342}) from
prismatic trivializations of $\ms e^\infty(\mc Q)$ over the
simplices of $K$ (Sec.~\vref{e_infty}).

To prove Theorem~D, it suffices to associate with the functor  $\mb
G$  an $(S^n, 0, \infty)$-bundle $e^{0,\infty}(\mb G)$
such that
\begin{itemize}
\item[(i)] forgetting the   $0$-section $e^{0,\infty}(\mb G)$ yields the
bundle $\ms e^\infty (\mb G)$;
\item[(ii)] the germ of $ e^{0,\infty}(\mb G)$ at the zero section
is isomorphic to $t_M$.
\end{itemize}
This can be achieved by the obvious special choice of
prismatic trivializations  associated with
Milnor's diagonal construction during the
construction of $\ms e^{\infty}(\mb G)$.

\def\cprime{$'$} \def\cprime{$'$}

\end{document}